\documentclass{article}
\usepackage{amsmath}
\usepackage{latexsym}
\usepackage{amssymb}
\newcommand{\qed}{\Box}
\newtheorem{thm}{Theorem}[section]
\newtheorem{la}[thm]{Lemma}
\newtheorem{Defn}[thm]{Definition}
\newtheorem{Remark}[thm]{Remark}
\newtheorem{prop}[thm]{Proposition}
\newtheorem{cor}[thm]{Corollary}
\newtheorem{Number}[thm]{\!\!}
\newenvironment{defn}{\begin{Defn}\rm}{\end{Defn}}
\newenvironment{rem}{\begin{Remark}\rm}{\end{Remark}}
\newenvironment{numba}{\begin{Number}\rm}{\end{Number}}
\newenvironment{proof}{{\noindent\bf Proof.}}%
                  {\nopagebreak\hspace*{\fill}$\qed$\medskip\par}   
\newcommand{\Punkt}{\nopagebreak\hspace*{\fill}$\qed$}
\newcommand{\wb}{\overline}

\newcommand{\mto}{\mapsto}
\newcommand{\ve}{\varepsilon}

\newcommand{\N}{{\mathbb N}}
\newcommand{\R}{{\mathbb R}}
\newcommand{\bO}{{\mathbb O}}
\newcommand{\Q}{{\mathbb Q}}
\newcommand{\Z}{{\mathbb Z}}
\newcommand{\C}{{\mathbb C}}
\newcommand{\K}{{\mathbb K}}

\newcommand{\cC}{{\cal C}}

\newcommand{\cE}{{\cal E}}

\newcommand{\cL}{{\cal L}}
\DeclareMathOperator{\Iso}{Iso}

\newcommand{\sub}{\subseteq}
\DeclareMathOperator{\GL}{GL}

\DeclareMathOperator{\car}{char}

\DeclareMathOperator{\id}{id}
\DeclareMathOperator{\Lip}{Lip}
\newcommand{\sbull}{{\scriptscriptstyle \bullet}}
\DeclareMathOperator{\Diff}{Diff}
\DeclareMathOperator{\Hom}{Hom}
\DeclareMathOperator{\graph}{graph}

\begin{document}
\begin{center}
{\bf\Large Finite order differentiability properties,\\[1mm]
fixed points and implicit functions\\[2.5mm]
over valued fields}\vspace{5.8 mm}\\
{\bf Helge Gl\"{o}ckner}\vspace{4mm}
\end{center}
\begin{abstract}\vspace{.5mm}
{\protect\noindent}We prove an implicit function theorem
for $C^k$-maps from arbitrary topological
vector spaces over valued fields to Banach spaces
(for $k\geq 2$).
As a
tool, we show the $C^k$-dependence
of fixed points on parameters
for suitable families of contractions
of a Banach space. Similar results are obtained
for $k$ times
strictly differentiable maps, and for
$k$ times Lipschitz differentiable maps.
In the real case, our results
subsume an implicit function theorem
for Keller $C^k_c$-maps from arbitrary
topological vector spaces to Banach spaces.
\end{abstract}
\section*{Introduction}
Generalizations of the implicit function theorem
for mappings from suitable real or complex topological
vector spaces to Banach spaces have been
obtained in various settings of analysis.
Hiltunen (\cite{Hil},\cite{Hi2})
studied implicit functions
in the framework
of Keller's $C^k_\Pi$-theory~\cite{Kel}
(cf.\ also \cite{Hi3} for recent generalizations
to non-Banach range spaces).
Teichmann~\cite{Tei} proved an implicit function
theorem in the ``convenient setting'' of analysis,
under more restrictive conditions.
In the framework of Keller's $C^k_c$-theory,
implicit functions from real and complex
topological vector spaces to Banach spaces
were discussed in~\cite{IMP}.\\[3mm]
The notion of a Keller $C^k_c$-map can
be generalized
to a notion of $C^k$-map between open subsets
of topological vector spaces over an arbitrary
(non-discrete) topological field~\cite{Ber}
(see also~\cite{SUR} for a survey).
For mappings between open subsets
of an ultrametric field, these $C^k$-maps
coincide with those usually considered
in Non-Archimedian Analysis (as in \cite{Sch}).
On the basis of~\cite{Ber},
the paper \cite{IMP} also provided an implicit function
theorem for $C^k$-maps
from metrizable topological vector spaces over complete
valued fields to Banach spaces,
with a possible loss of one order of
differentiability in the case of an infinite-dimensional
range space.\\[3mm]
In the present paper, we discuss implicit functions
from topological vector spaces over valued fields
to Banach spaces by a different method, which enables
us to remove the metrizability condition.
We can also avoid
the former loss of one order
of differentiability.
Thus,
local solutions
to $C^k$-equations are always $C^k$, as they should be
(if $k\geq 2$).
Stimulated by
discussions of mappings between
real Banach spaces in~\cite{Ir2},
our new strategy of proof now is
to discuss, in a first step,
the existence of fixed points for families of contractions
of a Banach space and their differentiable
dependence on parameters.
Next, we prove a suitable
Lipschitz Inverse Function Theorem
(which varies a result from \cite{Wel}).
The desired implicit function theorem
is then an immediate
consequence.
In fact, the Lipschitz Inverse Function
Theorem entails that, locally,
an implicit function $\lambda$
exists and that its value
$\lambda(x)$ at~$x$
can be constructed as a fixed point
of a suitable contraction~$g_x$ (if $k\geq 2$, say).
By the $C^k$-dependence of fixed points
on parameters, $\lambda(x)$ is a $C^k$-function of~$x$.\\[3mm]
It useful to work with several
variants of $C^k$-maps,
because the mere $C^k$-property
(notably, the $C^1$-property)
is slightly too weak for some
of our purposes, and envisaged
applications.
We therefore introduce
so-called \emph{$k$ times strictly differentiable mappings}
($SC^k$-maps) between open subsets
of topological $\K$-vector spaces
(and subsets with dense interior),
as well as \emph{$k$ times Lipschitz differentiable
maps}
($LC^k$-maps).
Our definition of $SC^k$-maps
generalizes an earlier definition
for mappings from normed spaces to
polynormed vector spaces from~\cite{IMP}.
Furthermore, a mapping between
open subsets of Banach spaces
is once strictly differentiable
in our sense if and only if it is strictly
differentiable at each point in the sense
of Bourbaki~\cite[1.2.2]{Bou} (cf.\ also \cite{Lea}
and \cite{Car}).
We also mention the ``$r$-Lipschitz maps''
$\bO\to\K$ studied by Barsky~\cite{Bar},
where $\K$ is a local field and $\bO$
its maximal compact subring.
The differentiability properties just described
are related as follows:
\[
C^{k+1}\; \Longrightarrow\;
LC^k\;
\Longrightarrow\;
SC^k\;
\Longrightarrow
\;
C^k\,.
\]
Thus
$C^\infty$-maps, $SC^\infty$-maps
and $LC^\infty$-maps all coincide,
but we have a certain range
of finite-order differentiability
properties.\\[3mm]
Among our main results is the following
generalization of the Implicit Function Theorem
(Theorem~\ref{genimp}).\\[3mm]
{\bf Generalized Implicit Function Theorem.}
\emph{Let $\K$ be a valued field,
$E$ be a topological $\K$-vector space,
$F$ be a Banach space over~$\K$,
and $f\colon U\times V\to F$ be a
map, where $U\sub E$ is a subset with dense
interior and $V\sub F$ is open.
Given $x\in U$, abbreviate $f_x:=f(x,\sbull)
\colon V\to F$. Assume that $\K$,
$F$,
$k\in \N\cup\{\infty\}$ and~$f$
have properties as shown
in the following table.
Furthermore, assume
that $f(x_0,y_0)=0$ for some $(x_0,y_0)\in U\times V$
and $f_{x_0}'(y_0)\in \GL(F)$.
Then there exists an open neighborhood
$U_0\sub U$ of~$x_0$, an open neighborhood
$V_0\sub V$ of $y_0$, and a
map
$\lambda \colon U_0\to V_0$ such that
\[
\{(x,y)\in U_0\times V_0\colon f(x,y)= 0 \}\;=\;
\graph\,\lambda \,,
\]
where $\lambda$ has the differentiability
property shown in the table}:\vspace{1mm}
\begin{center}
\begin{tabular}{||c|c|c|c|c||}\hline\hline
$\K$ & $F$ & $k$ &  $f$ & $\lambda $ \\ \hline\hline
arbitrary & arbitrary  & arbitrary & $LC^k$ & $LC^k$ \\ \hline
arb. & arb. & arb. & $SC^k$ & $SC^k$\\ \hline
arb. & arb. & $\geq 2$ & $C^k$ & $C^k$ \\ \hline
locally compact & $\dim\, F<\infty$ & arb. & $C^k $ & $C^k$\\ \hline\hline
\end{tabular}
\vspace{3mm}
\end{center}
The generalized implicit function
theorem is deduced from
a suitable ``Inverse Function Theorem with Parameters''
(Theorem~\ref{advif}),
dealing with families of local diffeomorphisms.
This theorem is our
actual main result.
As a technical tool, in Theorem~\ref{fixpardep}
we prove
differentiable dependence
of fixed points on
parameters,
for uniform families of contractions
in the sense of Definition~\ref{defunicon}:\\[3mm]
{\bf Theorem (on the Parameter-Dependence of Fixed Points).}
\emph{Let $E$ be a topological $\K$-vector space,
$F$ be a Banach space over~$\K$
and
$f\colon P\times U\to F$
be a map, where
$P\sub E$ is a subset with dense interior
and $U\sub F$ is open.
Given $p\in P$, abbreviate $f_p:=f(p,\sbull)\colon
U\to F$.
Assume that $(f_p)_{p\in P}$
is a uniform family of contractions
and assume that $f$ is $C^k$, $SC^k$, resp., $LC^k$
for some $k\in \N_0\cup\{\infty\}$.
Let $Q$ be the set of all $p\in P$ such that
$f_p$ has a fixed point $x_p$.
Then the following holds:}
\begin{itemize}
\item[(a)]
\emph{$Q$ is open in~$P$};
\item[(b)]
\emph{The map $\phi\colon Q\to F$, $p\mto x_p$
is $C^k$, $SC^k$, resp., $LC^k$.}\vspace{1mm}
\end{itemize}
{\bf Applications.} It is clear that generalizations
of such basic and central results as
the inverse- and implicit function theorems have
immediate implications.
\begin{itemize}
\item
In \cite{ARB}, the ultrametric
inverse function theorem with parameters
in a Fr\'{e}chet space
is used to prove that the inversion map
$\Diff(M)\to \Diff(M)$, $\gamma\mto \gamma^{-1}$
of the diffeomorphism group of a paracompact,
finite-dimensional smooth manifold over a
local field is smooth.
Also composition being smooth,
$\Diff(M)$ is a Lie group
(see also \cite{IMP} for an outline
of the proof).\footnote{Cf.\
also \cite{Lu1}, \cite{Lu8} for certain diffeomorphism
groups for $\car \,\K=0$.}
\item
Irwin~\cite{Ir1} explained how to construct
stable manifolds around hyperbolic fixed points
for discrete dynamical systems
modeled on real Banach spaces,
with the help of the implicit function theorem
(see also \cite{Wel}).
In \cite{INV}, the
$C^k$-, $SC^k$- and $LC^k$-versions
of our inverse function theorem
are used to adapt Irwin's method to
dynamical systems modeled on Banach
spaces over valued fields.
Stable manifolds
around hyperbolic fixed points
of the respective differentiability class
are constructed (and also analytic ones).
Furthermore,
the $LC^k$-inverse function theorem with parameters
is used in \cite{INV} to study the dependence
of the stable manifolds on the non-linearity
(see also \cite{IMP} and \cite{SUR} for an outline of
the method).
\item
Adapting Irwin's discussion of pseudo-stable
manifolds (\cite{Ir2}; see also \cite{LaW}),
one can also construct pseudo-stable manifolds
around hyperbolic fixed points
for dynamical systems modeled on Banach
spaces over valued fields
(see \cite{IN2}).
At the heart of these studies
is a refined analysis
of the dependence of fixed
points on parameters (in specialized situations).
\item
Varying a classical idea by Robbin
(see also Chow and Hale~\cite{CaH}
and \cite{Fre}),
$C^{k+n}$-solutions to (systems of)
$p$-adic differential
equations of the form $y^{(k)}=f(x,y,y',\ldots, y^{(k-1)})$
for $f$ an $C^n$-map and $k\in \N$, $n\in \N_0$
(with an extra local Lipschitz condition
if $n=0$)
can be constructed
using our inverse function theorems (with and without
parameters),
which depend on initial conditions and parameters
in a controlled way
(work in progress;
cf.\ \cite[\S\,65]{Sch} for $C^1$-solutions to scalar-valued
first order equations).
\item
Various applications of the Inverse Function Theorem
with Parameters in the structure theory of
infinite-dimensional
real Lie groups can be found in \cite{GaN}.
For example, as shown by Neeb,
it can be used to prove
that locally compact subgroups of locally exponential
Lie groups are Lie subgroups (as in the
classical case of Banach-Lie groups).\footnote{A smooth Lie group
$G$ modeled on a real locally
convex space is called \emph{locally exponential}
if it has an exponential function
$\exp_G\colon L(G)\to G$
which induces a local $C^\infty$-diffeomorphism
from an open $0$-neighborhood
in $L(G)$ onto an open identity neighborhood
in~$G$.}
\end{itemize}
The stable, unstable and center-stable
manifolds available through the preceding constructions
are useful for the theory of finite-dimensional
Lie groups over local fields.
In~\cite{SC2},
they are used to generalize structure theorems
for $p$-adic Lie groups and their automorphisms
(as in \cite{Wan} and \cite{SC1})
to the case of Lie groups
over local fields of positive characteristic
(under suitable hypotheses).
Further applications of the generalized
implicit function theorems provided here
(and their precursors from~\cite{IMP})
are summarized in the introduction of~\cite{IMP}.\\[3mm]
{\bf Structure of the article.}
In Section~\ref{secprel}
(as well as Appendix~\ref{appC0}),
we compile definitions,
notation and basic facts
concerning the differential calculus
of $C^k$-maps over topological fields.
We also introduce a certain concept of
a ``gauge'' as a convenient substitute for
seminorms when dealing with arbitrary
(not necessarily polynormed) topological
vector spaces over valued fields.
In Section~\ref{seclip}, we
discuss Lipschitz continuous maps
between subsets of topological vector spaces
over valued fields and Lipschitz differentiable
maps between subsets with dense interior.
Next, we define and discuss
strictly differentiable maps
and clarify their relations
to $C^k$-maps and $LC^k$-maps
(Section~\ref{secstrict}).
Fixed points of contractions and their
dependence on parameters are
discussed in Section~\ref{secfp},
and the results obtained are then used
in Section~\ref{sectimpl} to deduce
implicit function theorems as well as
inverse function theorems
with and without parameters.
Appendix~\ref{apptechni}
records a variant of a technical lemma
(not needed in the main text).
Appendix~\ref{newppx} was added in 2007.
It describes how some of our results
involving $C^1$-maps on finite-dimensional
vector spaces over locally compact fields
can be generalized to
$C^1$-maps on finite-dimensional
vector spaces over
complete valued fields.
%
%
%
%
%
\section{Preliminaries and basic facts}\label{secprel}
The general framework
for our studies is the
differential calculus
of $C^k$-maps between open subsets
of topological vector spaces
over non-discrete topological fields
developed in~\cite{Ber}.
In this section, we first recall 
basic definitions and facts from this
theory, which applies to arbitrary
topological ground fields,
and slightly extend them by replacing open sets
by sets with dense interior.
We then introduce various
concepts which are important
when dealing with \emph{valued} fields.
In particular, we shall define
certain generalizations of seminorms,
which we call ``gauges''.
With the help of
these gauges, we can treat general
topological vector spaces as though they were
locally convex.
Many definitions,
results and proofs will look exactly
as in the locally convex (or polynormed)
case, except that continuous seminorms
have been replaced with gauges.\footnote{A similar
idea is also implicit in Lang's definition of
total differentiability in general real topological
vector spaces \cite[I, \S\,3]{Lan}.}\\[3mm]
All topological fields are assumed
Hausdorff and non-discrete;
all topological vector spaces are assumed Hausdorff.
In Definition~\ref{def1}--Lemma~\ref{symmetries}, we shall
assume that $\K$ is a topological field;
in Definition~\ref{defvalued}--Lemma~\ref{seeRooj},
we assume that $(\K,|.|)$
is a (non-trivial) valued field.
\subsection*{{\boldmath $C^k$}-maps over topological fields}
Let $E$ and $F$ be topological $\K$-vector
spaces, and $f\colon  U\to F$ be a mapping,
defined on a subset $U\sub E$
with dense interior.
Then the directional difference quotient
\[
f^{]1[}(x,y,t)\; :=\;
\frac{f(x+ty)-f(x)}{t}
\]
makes sense for all $(x,y,t)$ in the subset
\[
U^{]1[}\; :=\; \{(x,y,t)\in U\times E\times \K^\times\colon x+ty\in U\}
\]
of $E\times E\times \K$.
To define directional derivatives,
we must enlarge this set by allowing also
the value $t=0$.
Hence, we consider now
\[
U^{[1]}\;:=\; \{(x,y,t)\in U\times E\times \K
\colon x+ty\in U\}\,.
\]
Thus $U^{[1]}=U^{]1[}\cup (U\times E\times \{0\})$, as a disjoint union.
%
%
\begin{defn}\label{def1}
$f\colon U\to F$ is called \emph{continuously
differentiable} (or $C^1$)
if $f$ is continuous ($C^0$)
and there exists a continuous map $f^{[1]}\colon U^{[1]}\to F$
which extends
$f^{]1[}\colon U^{]1[}\to F$.
\end{defn}
Thus, we assume the existence of a continuous map
$f^{[1]}\colon U^{[1]}\to F$ such that
\[
f^{[1]}(x,y,t)\; =\;
\frac{f(x+ty)-f(x)}{t}\quad\mbox{for all $(x,y,t)\in U^{[1]}$
such that $t\not=0$.}
\]
If it exists, then $f^{[1]}$ is unique,
by part\,(a) of the following lemma:
%
%
\begin{la}\label{U1dense}
Let $E$ be a topological $\K$-vector
space, and $U\sub E$ be a subset with dense interior~$U^0$.
Then the following holds:
\begin{itemize}
\item[\rm (a)]
$U^{]1[}$ is a dense open subset of $U^{[1]}$;
\item[\rm (b)]
$U^{[1]}$ is a subset of $E\times E\times \K$
with dense interior.
\item[\rm (c)]
$W^{]1[}$ is open in $E\times E\times \K$
and dense in $U^{[1]}$,
for each dense subset $W\sub U$ which is open in~$E$
$($e.g., for $W:=U^0)$.
\end{itemize}
\end{la}
\begin{proof}
Clearly $U^{]1[} =U^{[1]}\cap (E\times E\times \K^\times)$
is open in $U^{[1]}$.
Furthermore, $W^{]1[}$ (and hence $(U^0)^{]1[}$)
is open in $E\times E\times\K$,
as a consequence of the continuity of the maps $(x,y,t)\mto x$
and $(x,y,t)\mto x+ty$ on $E\times E\times\K$.
All other assertions will follow readily
if we can show that $W^{]1[}$ is dense
in $U^{[1]}$.
To this end, let $(x,y,t)\in U^{[1]}$
and $X \sub E$, $Y \sub E$ and $S \sub \K$
be open neighborhoods of $x$, $y$ and $t$, respectively.

If $t\not=0$,
then $x+tY$ is a neighborhood
of $x+ty\in U\sub \wb{W}$, whence
$(x+tY)\cap W\not=\emptyset$.
We therefore find $y'\in Y$ such that $x+ty'\in W$.
Then $(W-ty')\cap X$ is a neighborhood
of $x$ and hence has non-empty intersection with~$W$.
We pick $x'\in W\cap (W -ty')\cap X$.
Then $x'\in W$ and $x'+ty'\in W$ and thus
$(x',y',t)\in W^{]1[}$.
Furthermore, $(x',y',t)\in X\times Y\times S$.

If $t=0$, we pick an element $x'\in X\cap W$.
By openness of $W$ and continuity
of scalar multiplication,
after shrinking~$Y$ and~$S$ we may assume
that $x'+SY\sub W$.
Thus $\{x'\}\times Y \times S \sub W^{]1[}$.
We choose $y'\in Y$ and $s\in S\setminus\{0\}$.
Then $(x',y',s)\in W^{]1[}\cap (X\times Y\times S)$.

We have shown that $W^{]1[}$ is dense in $U^{[1]}$,
which completes the proof.
\end{proof}
Lemma~\ref{U1dense}\,(b)
facilitates to define $C^k$-maps by recursion.
\begin{defn}\label{defnCr}
Let $E$ and $F$ be topological $\K$-vector
spaces and $f\colon  U\to F$ be a map,
defined on a subset $U\sub E$
with dense interior.
Given $k\in \N$ with $k\geq 2$,
we say that $f$ is \emph{$k$ times continuously differentiable}
(or $C^k$) if $f$ is $C^1$ and $f^{[1]}\colon U^{[1]}\to F$
is $C^{k-1}$.
The map~$f$ is called $C^\infty$
(or \emph{smooth})
if it is $C^k$ for all $k\in \N_0$.
\end{defn}
\begin{numba}
For example, every continuous linear
map
$\lambda\colon  E\to F$ is smooth~\cite{Ber},
with $\lambda^{[1]}(x,y,t)=\lambda(y)$
for all $(x,y,t)\in E\times E\times \K$
(thus $\lambda^{[1]}$ is again continuous linear).
\end{numba}
\begin{numba}
The domain of the mapping $f^{[2]}:=(f^{[1]})^{[1]}$
is the set
\[
U^{[2]}\; :=(U^{[1]})^{[1]}\; \sub \; E\times E\times \K\times
E\times E\times \K\times \K\,.
\]
We write $U^{[k]}:=(U^{[1]})^{[k-1]}$
for the domain of $f^{[k]}$,
and set $U^{]k[}:=(U^{]1[})^{]k-1[}$.
\end{numba}
%
%
\begin{rem}\label{usefulsimp}
Applying Lemma~\ref{U1dense}\,(c) to $U^{[1]}$
and its dense open subset
$W:=(U^0)^{]1[}$,
we find that $(U^0)^{]2[}$ is dense
in $U^{[2]}$, and inductively that
$(U^0)^{]k[}$ is dense in $U^{[k]}$,
for each $k\in \N$.
Hence $(f|_{U^0})^{[k]}$ uniquely determines
$f^{[k]}$ in particular (if it exists).
This enables us to extend
all relevant results and proofs from \cite{Ber}
(where only open domains were considered)
to the case of non-open domains,
by trivial arguments based on continuous
extension.
We shall therefore cite (and apply)
results from \cite{Ber} freely
also for mappings on non-open domains.
We remark that, if one is working exclusively
with Hausdorff topological vector spaces
over Hausdorff topological fields
(in contrast to \cite{Ber}, where
the scope was wider),
it can be convenient to replace
the somewhat technical axioms from
\cite{Ber} by simpler ones.
Such variants are suggested in
Appendix~\ref{appC0}.
\end{rem}
%
%
\begin{rem}\label{remonesev}
A trivial induction shows that $U^{[k]}=(U^{[k-1]})^{[1]}$,
and that $f$ is $C^k$ if and only if $f$ is $C^{k-1}$
and $f^{[k-1]}$ is $C^1$
(cf.\ \cite[Rem.\,4.2]{Ber}).
In this article, we shall
avoid to use $f^{[k]}$ for higher~$k$
as far as possible. Usually, we only need~$f^{[1]}$.
\end{rem}
\begin{numba}\label{differentials}
Given a $C^1$-map $f\colon  U\to F$ as before, we define
its \emph{directional derivative} at $x\in U$
in the direction
$y\in E$ via
\[
(D_yf)(x)\;:=\; df(x,y)\;:=\; f^{[1]}(x,y,0)\,.
\]
If $x\in U^0$, then
\[
df(x,y)\; =\; f^{[1]}(x,y,0)\; =\;
\lim_{0\not=t\to 0} f^{[1]}(x,y,t)\;
=\; \lim_{0\not=t\to 0} {\textstyle
\frac{1}{t}(f(x+ty)-f(x))}\,,
\]
and thus $df(x,y)$ can be interpreted as a limit
of difference quotients.
The map $df\colon U\times E\to F$ is continuous, being a partial map
of~$f^{[1]}$, and it can be shown that
the ``differential'' $f'(x):=
df(x,\sbull)\colon  E\to F$
of $f$ at~$x$ is a continuous $\K$-linear map,
for each $x\in U$ (cf.\ \cite[Proposition~2.2]{Ber}).
If $f$ is $C^2$,
we define a continuous map $d^2f\colon  U\times E^2\to F$ via
$d_2f(x,y_1,y_2):= (D_{y_2}(D_{y_1}f))(x)$. Thus
\[
d^2f(x,y_1,y_2):=\lim_{t\to 0}
{\textstyle
\frac{1}{t}(df(x+ty_2,y_1)-df(x,y_1))}
=f^{[2]}((x,y_1,0),\, (y_2,0,0),\, 0)\, .
\]
Similarly, if $f$ is $C^k$,
we obtain continuous
maps $d^jf\colon  U\times E^j\to F$
for all $j\in \N_0$ such that $j\leq k$
via $d^jf(x,y_1,\ldots, y_j):=(D_{y_j}\cdots D_{y_1}f)(x)$.
It can be shown that
$d^jf(x,\sbull)\colon E^j\to F$
is a symmetric $j$-linear map (cf.\ \cite[Lemma~4.8]{Ber}).
\end{numba}
\begin{rem}
We remark that,
in the real locally convex case,
our $C^k$-maps coincide with Keller's $C^k_c$-maps
(as used, e.g., in \cite{RES}, \cite{GaN}, \cite{Mil}).
More precisely, let
$E$ be a real topological vector
space, $F$ be a locally convex real topological vector space,
$U\sub E$ be open,
$f\colon  U\to F$ be a map, and $k\in \N_0\cup\{\infty\}$
Then $f$ is $C^k$ if and only if
$f$ is continuous, the
limits $d^jf(x, y_1, \ldots, y_j):=
(D_{y_j}\cdots D_{y_1}f)(x)$
described in {\rm {\bf \ref{differentials}}}
exist for all $j\in \N$ with $j\leq k$,
and the maps
$d^jf\colon U\times E`j\to F$
so obtained are continuous
(see \cite[Proposition~7.4]{Ber}).
We mention that if also $E$ is locally convex,
then this characterization
remains valid for mappings on
locally convex subsets $U\sub E$
with dense interior (as considered
in~\cite{GaN}).
\end{rem}
\begin{numba}\label{chainr}
(Chain Rule).
If $E$, $F$, and $H$ are topological
$\K$-vector spaces, $U\sub E$ and $V\sub F$ are
subsets with dense interior, and $f\colon  U\to V\sub F$,
$g\colon  V\to H$ are $C^k$-maps,
then also the composition
$g\circ f\colon  U\to H$ is~$C^k$.
If $k\geq 1$, we have
$(f(x),f^{[1]}(x,y,t),t)\in V^{[1]}$
for all $(x,y,t)\in U^{[1]}$, and
\begin{equation}\label{formchain}
(g\circ f)^{[1]}(x,y,t)=g^{[1]}(f(x),f^{[1]}(x,y,t),t)\,.
\end{equation}
In particular, $d(g\circ f)(x,y)=dg(f(x),df(x,y))$
for all $(x,y)\in U\times E$ (cf.\ \cite[Prop.\,3.1 and 4.5]{Ber}).
\end{numba}
We recall from \cite[La.\,4.9]{Ber} that
being $C^k$ is a local property.
\begin{la}\label{Crlocal}
Let $E$ and $F$ be topological $\K$-vector spaces,
and $f\colon  U\to F$ be a map, defined
on an open subset $U$ of~$E$.
Let $k\in \N_0\cup\{\infty\}$. 
If there is an open cover $(U_i)_{i\in I}$
of~$U$ such that
$f|_{U_i}\colon  U_i\to F$ is
$C^k$ for each $i\in I$,
then $f$ is~$C^k$.\Punkt
\end{la}
\subsection*{Symmetry properties of {\boldmath $f^{[1]}$}}
The map $f^{[1]}$ and
also the higher different quotient maps
satisfy various identities,
which we need to exploit occasionally.
We now describe some properties
of $f^{[1]}$.
A symmetry property of $f^{[2]}§$
is described in Appendix~\ref{apptechni}.
Slightly less explicit (and more complicated)
results concerning $f^{[k]}$ for arbitrary~$k$
can be found in \cite[La.\,6.8]{ARB}.
They are also essential for the construction
of invariant manifolds in \cite{INV} and \cite{IN2}.
%
%
\begin{la}\label{firsiden}
Let $E$, $F$ be topological $\K$-vector spaces,
$U\sub E$ be a subset with dense interior
and $f\colon U\to F$ be $C^1$.
Then $(x+t y_2,y_1-y_2,t)\in U^{[1]}$
for all $x\in U$, $y_1,y_2\in E$ and $t\in \K$
such that $(x,y_1,t),(x,y_2,t)\in U^{[1]}$,
and
%
%
\begin{equation}\label{eqfirsid}
f^{[1]}(x,y_1,t)-f^{[1]}(x,y_2,t)\;=\;
f^{[1]}(x+ty_2,y_1-y_2, t)\,.
\end{equation}
\end{la}
\begin{proof}
The first assertion is clear since
$x+ty_2\in U$ (because $(x,y_2,t)\in U^{[1]}$)
and $x+ty_2+t(y_1-y_2)=x+ty_1\in U$
(because $(x,y_1,t)\in U^{[1]}$).
If $t\not=0$, then
\begin{eqnarray*}
f^{[1]}(x,y_1,t)-f^{[1]}(x,y_2,t) &=&
\frac{f(x+ty_1)-f(x)}{t}
-
\frac{f(x+ty_2)-f(x)}{t}\\
&=& \frac{f(x+ty_1)-f(x+ty_2)}{t}\\
&= &\frac{f(x+ty_2+t(y_1-y_2))-f(x+ty_2)}{t}\\
&=&
f^{[1]}(x+ty_2,y_1-y_2, t)\,,
\end{eqnarray*}
as desired.
If $t=0$, then (\ref{eqfirsid})
turns into the identity
$df(x,y_1)-df(x,y_2)=
df(x,y_1-y_2)$,
which is valid by linearity of~$df(x,\sbull)$.
\end{proof}
Also the following slightly more complicated
symmetry properties will be used.
They shall enable us to shrink certain
entries of $f^{[1]}$
while inflating others in a controlled way.
%
%
\begin{la}\label{symmetries}
Let $E$ and $F$ be topological vector
spaces over a topological field~$\K$, and $f\colon  U\to F$
be a $C^1$-map, defined on a subset $U\sub E$
with dense interior.
If $f$ is $C^1$, $t\in \K^\times$,
and $(x,y,s)\in E\times E\times \K$
such that $(x,y,ts)\in U^{[1]}$,
then also $(x,ty,s)\in U^{[1]}$, and
%
\begin{equation}\label{bracket1}
t\, f^{[1]}(x,y,ts)= f^{[1]}(x,ty,s)\,.
\end{equation}
\end{la}
\begin{proof}
Since $x+(ts)y=x+s(ty)$, it is obvious
that $(x,ty,s)\in U^{[1]}$ if and only if $(x,y,ts)\in U^{[1]}$.
In this case, we have
\[
tf^{[1]}(x,y,ts)={\textstyle \frac{1}{s}}(f(x+tsy)-f(x))
=f^{[1]}(x,ty,s)
\]
provided $s\not=0$;
if $s=0$, then
$f^{[1]}(x,ty,s)=f^{[1]}(x,ty,0)=df(x,ty)=tdf(x,y)=tf^{[1]}(x,y,0)
=tf^{[1]}(x,y,ts)$.
Thus (\ref{bracket1}) holds.
\end{proof}
\subsection*{Vector spaces
over valued fields, seminorms and gauges}
We now fix terminology and notation
concerning topological vector spaces over
valued fields, and compile some simple
observations which are
analogous to well-known facts concerning
locally convex spaces.
%
%
\begin{defn}\label{defvalued}
A \emph{valued field} is a field~$\K$, together
with an absolute value $|.|\colon  \K\to [0,\infty[$
(see \cite{Wie}); we require
furthermore that the absolute value be
non-trivial (meaning that it gives rise
to a non-discrete topology on~$\K$).
An \emph{ultrametric field}
is a valued field $(\K,|.|)$
whose absolute value satisfies the ultrametric inequality
\[
|x+y|\leq \max\{|x|,|y|\}\quad \mbox{for all $x,y\in \K$.}
\]
Locally compact, totally disconnected,
non-discrete topological fields will be
referred to as \emph{local fields}.
It is well known that every local field~$\K$ admits
an ultrametric absolute value
defining its topology~\cite{Wei}.
Fixing such an absolute value on~$\K$, we can consider
$\K$ as an ultrametric field.
\end{defn}
\begin{rem}
Note that we do not require that
valued fields (nor ultrametric fields)
be complete (with respect to
the metric induced by the absolute value).
Of course, we are interested exclusively
in complete valued fields --
but none of our results
will depend on completeness of~$\K$.
\end{rem}
Given a valued field $(\K,|.|)$,
we can speak of \emph{seminorms}
on $\K$-vector spaces
(as in \cite[Ch.\,II, \S\,1, no.\,1, Defn.\,1]{BTV}).
If $(\K,|.|)$ is an ultrametric field,
we can also speak of \emph{ultrametric seminorms}
(Bourbaki's ``ultra-semi-norms''),
satisfying the ultrametric inequality.
\begin{numba}
Recall that a topological vector
space $E$ over an ultrametric field~$\K$
is called \emph{locally convex}
if every $0$-neighborhood of~$E$
contains an open $\bO$-submodule
of~$E$, where $\bO:=\{t\in \K\colon  |t|\leq 1\}$
is the valuation ring of~$\K$.
Equivalently, $E$ is locally convex if and only if
its vector topology is defined by a family
of ultrametric continuous seminorms $\gamma\colon  E\to [0,\infty[$ on~$E$
(cf.\ \cite{Mon} for more information, or also
the discussions of Minkowski functionals
given below).
Let $\K$ be a valued field.
We call a topological $\K$-vector space
\emph{polynormed} if its vector topology
is defined by a family of continuous seminorms
(which need not be ultrametric seminorms
if $\K$ is an ultrametric field).
\end{numba}
\begin{numba}
A \emph{Banach space} over a valued field~$\K$
is a normed $\K$-vector space $(E,\|.\|)$
(see \cite[Ch.\,I, \S1, no.\,2]{BTV})
which is complete in the metric associated with~$\|.\|$.
\end{numba}
\begin{numba}
We shall not presume that normed spaces (nor
Banach spaces) over ultrametric fields be ultrametric,
unless saying so explicitly.
For example, $\ell^1(\Q_p)$
is a non-ultrametric (and non-locally convex)
Banach space over~$\Q_p$.
\end{numba}
We are mainly interested in mappings
between polynormed (or even locally convex)
spaces, but all of our results can be proved
just as well for mappings on
arbitrary topological vector spaces
over valued fields, without
(or with very little)
additional effort.
To make the general proofs
look like those for polynormed
vector spaces,
we now introduce a convenient
replacement for
continuous seminorms,
namely the more general concept of a \emph{gauge}.
\begin{defn}
Let $E$ be a topological
vector space over a valued field
$(\K,|.|)$. A \emph{gauge on~$E$}
is an upper semicontinuous map $\gamma \colon E\to [0,\infty[$
(also written $\|.\|_\gamma :=\gamma$)
which satisfies
$\gamma(t x) = |t|\gamma(x)$ for all $t\in \K$
and $x\in E$.
\end{defn}
\begin{rem}
(a)
The upper semicontinuity of a gauge $\gamma\colon E\to [0,\infty[$
means
that $\gamma^{-1}([0,r[)$ is open in~$E$,
for each $r>0$.
This is equivalent to the following condition:
For each
$x\in E$ and net
$(x_\alpha)$ in~$E$ that converges to~$x$, we have
\begin{equation}\label{eqlimsup}
{\textstyle \limsup_\alpha} \, \|x_\alpha\|_\gamma \; \leq \; \|x\|_\gamma\,.
\end{equation}

(b)
Every gauge is continuous at~$0$.
Indeed: Since
$\gamma(0)=0$ and $\gamma\geq 0$ pointwise,
this follows from the upper semicontinuity.

(c)
Sums of gauges and non-negative multiples $r\gamma$ of gauges
are gauges.
\end{rem}
\begin{rem}
Typical examples
of gauges are Minkowski functionals
of balanced, open $0$-neighborhoods.
We recall: If $E$ is a topological
vector space over a valued field $\K$,
then a subset $U\sub E$ is called
\emph{balanced} if $tU\sub U$
for all $t\in \K$ such that $|t|\leq 1$.
It is easy to see that
the filter of $0$-neighborhoods of~$E$
has a basis of balanced, open $0$-neighborhoods.
For any balanced, open $0$-neighborhood~$U$,
we define its \emph{Minkowski functional}
$\mu_U\colon E\to [0,\infty[$ via
\[
\mu_U(x)\;:=\; \inf\{|t|\colon \mbox{$t\in \K^\times$
such that $x\in tU$}\}\,.
\]
Clearly $\mu_U(sx)=|s|\mu_U(x)$ for all $x\in E$ and $s\in \K$.
To see that $\mu_U$ is a gauge,
it only remains to check its upper semicontinuity.
To this end, let
$x\in E$ and $\ve>0$. There is $t\in \K^\times $
such that $x\in tU$ and $|t|\leq \mu_U(x)+\ve$.
Then $tU$ is a neighborhood
of~$x$ such that $\mu_U(y)\leq |t|\leq \mu_U(x)+\ve$
for all $y\in tU$.

For later use, we observe that
%
\begin{equation}\label{bllmink}
\{x\in E\colon  \mu_U(x)< 1\}\;\sub\;
U\; \sub \; \{x\in E\colon  \mu_U(x)\leq 1\}
\; \sub \; t U
\end{equation}
for each $t\in \K$ such that $|t|>1$.
Note that, in contrast to the familiar real
case, $tU$ on the right hand side of (\ref{bllmink})
cannot be replaced by the closure
$\wb{U}$ in general,
as the example $E:=\K:=\C_p$,
$U:=\{x\in \C_p \colon |x|<1\}$ shows.
In this case, $U$ is closed,
and it is a proper subset of
$\{x\in \C_p\colon |x|\leq 1\}$.
\end{rem}
%
%
%
%
\begin{defn}\label{defblls}
Given a topological $\K$-vector space~$E$,
$x\in E$, $r>0$ and a gauge $\gamma$ on~$E$,
we set
$B_r^\gamma(x):=\{y\in E\colon \|y-x\|_\gamma<r\}$
and $\wb{B}_r^\gamma(x):=\{y\in E\colon $
$\|y-x\|_\gamma\leq r\}$.
If $(E,\|.\|)$ is a normed space over~$\K$
and the norm $\gamma:=\|.\|$ is understood,
we simply write
$B_r^E(x):=B_r^\gamma(x)$
and
$\wb{B}_r^E(x):=\wb{B}_r^\gamma(x)$.
We abbreviate $B_r(x):=B_r^E(x)$
when no confusion is possible,
and $B_r(x):=B_r^E(x)$.
\end{defn}
Given a balanced, open $0$-neighborhood $V\sub E$
and $t\in \K$ such that $|t|>1$,
the set $U:=t^{-1}V$ is a balanced,
open $0$-neighborhood.
The right hand side of (\ref{bllmink}) being $V$,
we deduce that
the balls $\wb{B}^\gamma_1(0)$
form a basis of $0$-neighborhoods for~$E$,
if $\gamma$ ranges through the set of all gauges on~$E$.
\begin{defn}
Let $E$ be a topological vector space over
a valued field~$\K$.
A set $\Gamma$ of gauges on~$E$ is called
a \emph{fundamental system of gauges}
if
\[
\{B_r^\gamma(0)\colon \gamma\in \Gamma,r\in \;]0,\infty[\}
\]
is a basis for the filter of $0$-neighborhoods
in~$E$.
\end{defn}
It is useful to single out a simple argument which will
be used repeatedly.
%
\begin{la}\label{observi}
Let $E$ be a
vector space over a valued field~$\K$
and $\gamma,\eta \colon E\to [0,\infty[$
be mappings such that $\gamma(tx)=|t|\gamma(x)$
and $\eta(tx)=|t|\eta(x)$
for all $x\in E$, $t\in \K$.
We assume that there are $r,s>0$
such that $B_s^\eta(0)\sub B_r^\gamma(0)$,
using notation as in {\rm Definition~\ref{defblls}}.
Then
\[
\gamma\; \leq \; rs^{-1}|a|^{-1}\eta
\]
for each $a\in \K^\times$ such that $|a|<1$.
\end{la}
\begin{proof}
Let $x\in E$.
If $\eta(x)\not=0$, pick $k\in \Z$
such that $|a|^{k+1}\leq s^{-1}\eta(x) <|a|^k$.
If $\eta(x) =0$, pick any $k\in \Z$.
Then $\eta(a^{-k}x)<s$
and thus $|a|^{-k}\gamma(x)=
\gamma(a^{-k}x)<r$,
whence
\begin{equation}\label{ggg3}
\gamma(x) \; \leq \; r|a|^k\,.
\end{equation}
If $\eta(x)>0$,
then the right hand side of (\ref{ggg3})
is $\leq r s^{-1}|a|^{-1}\eta(x)$,
as required.
If $\eta(x)=0$, we can choose $k$ arbitrarily large,
and thus (\ref{ggg3})
entails that $\gamma(x)=0\leq rs^{-1}|a|^{-1}\eta(x)$
also in this case.
\end{proof}
%
%
\begin{rem}\label{onlymin}
If $\gamma$ is a gauge on~$E$,
$a\in \K^\times$ such that $|a|<1$,
and $U:=B_{|a|}^\gamma(0)$,
then $B^{\mu_U}_1(0)\sub U$ by (\ref{bllmink}),
whence
$\gamma  \leq \mu_U$
by Lemma~\ref{observi}.
Therefore Minkowski functionals
form a fundamental
system of gauges in particular.
\end{rem}
%
%
\begin{rem}\label{morefund}
If $\Gamma$ is a fundamental
system of gauges for~$E$,
and $\gamma$ is a gauge on~$E$,
then there exists a gauge $\eta\in \Gamma$
and $c>0$ such that $\gamma\leq c\eta$.
In fact, there exists $\eta\in \Gamma$
and $r>0$ such that $B^\eta_r(0)\sub B^\gamma_1(0)$.
Then $\gamma\leq r^{-1}|a|^{-1}\eta$,
by Lemma~\ref{observi}.
\end{rem}
%
%
\begin{la}\label{morebasstr}
Let $E$ be a topological vector space over
a valued field~$\K$,
and $\Gamma$ be a fundamental
system of gauges on~$E$.
Then the following holds:
\begin{itemize}
\item[\rm (a)]
The set of balls
$\{B_r^\gamma(x)
\colon \mbox{$\gamma\in \Gamma$, $r\in \;]0,\infty[$}\}$
is a basis for the filter of neighborhoods
of~$x$ in~$E$, and so is
$\{\wb{B}_r^\gamma(x)
\colon \mbox{$\gamma\in \Gamma$, $r\in \;]0,\infty[$}\}$.
If $]0,\infty[\,\Gamma\sub \Gamma$,
then also $\{B_1^\gamma(x)\colon \gamma\in \Gamma\}$
and $\{\wb{B}_1^\gamma(x)\colon \gamma\in \Gamma\}$
are bases.
\item[\rm (b)]
A map $f\colon X\to E$ from
a topological space to~$E$ is continuous
at $x\in X$ if and only if,
for each gauge $\gamma\in \Gamma$ and $\ve>0$,
there exists a neighborhood $U$ of~$x$ in~$X$
such that $\|f(y)-f(x)\|_\gamma<\ve$ for all
$y\in U$.
Or equivalently:
$\|f(x_\alpha)-f(x)\|_\gamma\to 0$,
for each $\gamma\in \Gamma$ and each
net $(x_\alpha)$ in~$X$ that converges to~$x$.
\item[\rm (c)]
Suppose that also $E_1$ is a topological $\K$-vector space,
$\Gamma_1$ a fundamental set of gauges on~$E_1$,
and $\alpha\colon E\to E_1$ a linear map.
Then $\alpha$ is continuous if and only if,
for each $\gamma \in \Gamma_1$,
there exists $\eta\in \Gamma$ and a constant
$c\in [0,\infty[$ such that
\[
\|\alpha(x)\|_\gamma\;\leq\; c \|x\|_\eta\qquad
\mbox{for all $\, x\in E$.}
\]
If $[0,\infty[\,\Gamma\sub \Gamma$
and $\alpha$ is continuous,
then $\eta$ and~$c$ can always be chosen
such that $c=1$.
\end{itemize}
\end{la}
\begin{proof}
(a) Given $x\in E$, the map $E\to E$, $y\mto x+y$ is a
homeomorphism. Hence,
if $\Gamma$ is a fundamental system of gauges
on~$E$, then
$\{B_r^\gamma(x)=x+B_r^\gamma(0)
\colon \mbox{$\gamma\in \Gamma$, $r\in \;]0,\infty[$}\}$
is a basis for the filter of neighborhoods
of~$x$ in~$E$.
Similarly, so is
$\{\wb{B}_r^\gamma(x)=x+B_r^\gamma(0)
\colon \mbox{$\gamma\in \Gamma$, $r\in \;]0,\infty[$}\}$.
Since
$B_1^{\gamma/r}(x)=B_r^\gamma(x)$ and
$\wb{B}_1^{\gamma/r}(x)=\wb{B}_r(x)$
for each $\gamma\in \Gamma$ and $r>0$,
the final assertions follow.

(b) As a consequence of~(a),
the map $f$ is continuous
at~$x$ if and only if $f^{-1}(B_r^\gamma(f(x)))=\{y\in X\colon
\|f(y)-f(x)\|_\gamma<r\}$
is a neighborhood of~$x$,
for each $\gamma\in \Gamma$ and $r>0$.
Hence the first assertion holds.
It also follows from (a)
that $f(x_\alpha)\to f(x)$
if and only if $f(x_\alpha)\in B_r^\gamma(f(x))$
eventually for all $\gamma\in \Gamma$
and $r>0$.
Since $f(x_\alpha)\in B_r^\gamma(f(x))$
if and only if $\|f(x_\alpha)-f(x)\|_\gamma<r$,
the second assertion follows.

(c) Suppose that,
for each $\gamma$, a gauge $\eta$ can be chosen
as described. Then
$\alpha(B_{c^{-1}r}^\eta(0))\sub B_r^\gamma(0)$
for each $r>0$,
entailing that $\alpha$ is continuous
at~$0$ and hence continuous.
If $\alpha$ is continuous and $\gamma$ is a gauge on~$E_1$,
then $U:=\alpha^{-1}(B_1^\gamma(0))=B_1^{\gamma\circ\alpha}(0)$
is a balanced open $0$-neighborhood
in~$E$. There exists $\eta\in \Gamma$
and $r>0$ such that
$B_r^\eta(0)\sub U$.
Applying Lemma~\ref{observi}
to $\gamma\circ \alpha$ and $\eta$,
we obtain $\gamma\circ \alpha \leq
r^{-1}|a|^{-1}\eta=c\eta$
with $c:=r^{-1}|a|^{-1}$.
\end{proof}
\begin{defn}\label{likeopnorm}
In the situation of
Lemma~\ref{morebasstr}\,(c), we set
\[
\|\alpha\|_{\gamma,\eta}:=\min\{c\geq 0\colon
\mbox{$\|\alpha(x)\|_\gamma \leq c \|x\|_\eta$
for all $x\in E$}\}\,.
\]
\end{defn}
Note that the triangle
inequality need not hold for gauges.
The following lemma
provides a certain substitute.
%
%
\begin{la}\label{substitut}
If $E$ is a topological vector
space over a valued field~$\K$
and $U,V\sub E$ are balanced open $0$-neighborhoods
such that $V+V\sub U$, then
\begin{equation}\label{repltri1}
\mu_U(x+y)\;\leq\; \max\{\mu_V(x),\mu_V(y)\}\quad
\mbox{for all $x,y\in E$.}
\end{equation}
As a consequence,
for each gauge $\|.\|_\gamma$
on~$E$, there is a gauge $\|.\|_\eta$
on~$E$
such~that
%
\begin{equation}\label{repltri2}
\|x+y\|_\gamma \;\leq\; \max\{\|x\|_\eta , \|y\|_\eta\}\quad
\mbox{for all $x,y\in E$.}
\end{equation}
\end{la}
\begin{proof}
Let $x,y\in E$.
Given $\ve>0$,
there exists $t\in \K^\times$
such that $|t|\leq \mu_V(x)+\ve$
and $x\in tV$,
and $s\in \K^\times$
such that $|s|\leq \mu_V(y)+\ve$
and $y\in sV$.
Assume that $|s|\leq |t|$
(the case $|s|>|t|$ is similar).
Then $x,y\in tV$
and thus $x+y\in tV+tV=t(V+V)\sub tU$,
entailing that $\mu_U(x+y)\leq |t|\leq
\max\{\mu_V(x),\mu_V(y)\}+\ve$.
As $\ve>0$ was arbitrary,
(\ref{repltri1}) follows.

Given $\gamma$, by Remark~\ref{onlymin}
there exists $U$ such that $\gamma\leq \mu_U$.
Choosing~$V$ as before
we then have (\ref{repltri2})
with $\eta:=\mu_V$.
\end{proof}
Recall that a subset $B\sub E$ of a topological vector space~$E$
over a valued field~$\K$ is called \emph{bounded}
if, for each $0$-neighborhood $U\sub E$,
there exists $t\in \K^\times$ such that $B\sub tU$.
\begin{defn}\label{deftopops}
If $E$ and~$F$ are topological vector
spaces over a valued field~$\K$,
we equip
the space $\cL(E,F)$ of continuous $\K$-linear maps $E\to F$
with the topology of
uniform convergence on bounded subsets of~$E$.
Thus, we equip $\cL(E,F)$
with the unique vector topology which has
the sets
\[
[B,U]\, :=\, \{\alpha\in \cL(E,F)\colon \alpha(B)\sub U\}
\]
as a filter basis of $0$-neighborhoods,
where $B$ ranges through the bounded
subsets of~$E$ and $U$ through
the $0$-neighborhoods
of~$F$ (it easily follows from \cite[Ch.\,I, \S\,1, no.\,5,
Prop.\,4]{BTV} that such a vector
topology exists).
As usual, we abbreviate $\cL(E):=\cL(E,E)$
and set
\[
\GL(E)\, :=\, \cL(E)^\times\, :=\, \{\alpha\in \cL(E)\colon
\mbox{$\exists \beta\in \cL(E)$ s.t.\ $\alpha\circ \beta=\beta\circ\alpha
=\id_E$}\}\,.
\]
\end{defn}
%
%
\begin{rem}\label{defgop}
Note that the gauges
\[
\|.\|_{\gamma, B}\colon \cL(E,F)\to [0,\infty[\,,\quad
\|\alpha\|_{\gamma,B}:=
\sup\{\|\alpha(x)\|_\gamma\colon x\in B\}
\]
define the vector topology on~$\cL(E,F)$,
for $B$ and $\gamma$
ranging through the bounded
subsets of~$E$ and the gauges
on~$F$, respectively.
If $F$ is polynormed (resp.,
locally convex),
then also $\cL(E,F)$ is
polynormed (resp., locally convex),
because $\|.\|_{\gamma,B}$
is a seminorm (resp., an ultrametric
seminorm) if so is~$\gamma$.
If $(E,\|.\|_E)$ is normed and $F$ is polynormed,
then the vector topology on~$\cL(E,F)$
arises from the family of continuous seminorms
$\|.\|_\gamma\colon \cL(E,F)\to [0,\infty[$
defined for $\alpha \in \cL(E,F)$ via
\begin{equation}\label{defopnorm}
\|\alpha \|_\gamma:=\sup\{ \|\alpha(v)\|_\gamma
\cdot \|v\|^{-1}_E \colon \,
0\not=v\in E\}\in [0,\infty[
\end{equation}
where $\gamma$ ranges through
the continuous seminorms on~$F$
(see Lemma~\ref{seeRooj} below; cf.\ also \cite[note on p.\,59]{Roo}).
By definition, we have
%
\begin{equation}\label{hencegood}
\|\alpha(x)\|_\gamma\;\leq\; \|\alpha\|_\gamma
\|x\|_E\qquad\mbox{for all $x\in E$.}
\end{equation}
If both $(E,\|.\|_E)$ and $(F,\|.\|_F)$ are normed, then
$\cL(E,F)$ is normable; its vector topology arises
from the \emph{operator norm} $\|.\|:=\|.\|_\gamma$ defined
in (\ref{defopnorm}), with
$\gamma :=\|.\|_F$.
It is easy to see that the operator norm $\|.\|$
is ultrametric if the norm $\|.\|_F$
is ultrametric.
If $F$ is a Banach space here,
then also $\cL(E,F)$ is complete
and hence a Banach space
(this can be shown as in the real case).
It easily follows from the definition of
the operator norm on $\cL(E)$
that $\|\alpha\circ\beta\|\leq\|\alpha\|\cdot\|\beta\|$
for all $\alpha,\beta\in \cL(E)$.
\end{rem}
The following observation is occasionally useful.
%
\begin{la}\label{seeRooj}
Let $\alpha\colon E\to F$ be a linear map,
$\gamma$ be a gauge on~$F$ and $\eta$ be a gauge on~$E$
such that
\[
\|\alpha\|_{\gamma,B}
\, :=\, 
\sup\{\|\alpha.v\|_\gamma\colon v\in B\}\,<\infty\,,
\]
where $B:=B^\eta_1(0)$.
Then
%
\begin{equation}\label{eqgaug}
\|\alpha\|_{\gamma, B}\; \leq\;
\|\alpha\|_{\gamma,\eta}
\; \leq \; |a|^{-1}\|\alpha\|_{\gamma, B}
\end{equation}
for each $a\in \K^\times$ such that $|a|<1$.
If $|\K^\times|$ is dense in $[0,\infty[$,
then $\|\alpha\|_{\gamma,B}=\|\alpha\|_{\gamma,\eta}$.
\end{la}
\begin{proof}
For each $x\in B$, we have
$\|\alpha(x)\|_\gamma\leq \|\alpha\|_{\gamma,\eta}\|x\|_\eta\leq
\|\alpha\|_{\gamma,\eta}$
and thus $\|\alpha\|_{\gamma, B}:=\sup\{\alpha(x)\|_\gamma\colon
x\in B\}\leq \|\alpha\|_{\gamma,\eta}$.
Hence the first half of (\ref{eqgaug}) holds.
To prove the second,
let $r>\| \alpha\|_{\gamma,B}$.
Then $\alpha(B)\sub B_r^\gamma(0)$
and thus $B\sub B_r^{\gamma\circ\alpha}(0)$.
Then $\gamma\circ \alpha\leq r^{-1}|a|^{-1}\eta$,
by Lemma~\ref{observi}, entailing that
$\|\alpha\|_{\gamma,\eta}\leq r |a|^{-1}$.
Letting $r\to \|\alpha\|_{\gamma,B}$,
the second half of (\ref{eqgaug}) follows.
The other assertions are now immediate.
\end{proof}
Note that the notation $\|.\|_{\gamma,B}$
is abused in Lemma~\ref{seeRooj}
for a set $B$ which need not be bounded
(unless $E$ is normed).\\[2.5mm]
If $E$ is a Banach space,
then $\GL(E)$ is open in $\cL(E)$
and is a topological group.
%
%
\begin{prop}\label{glislie}
If $(E,\|.\|)$ is a Banach space,
then $\GL(E)$ is open in $\cL(E)$
and the inversion map $\iota\colon \GL(E)\to\GL(E)$,
$\iota(\alpha):=\alpha^{-1}$ is continuous.
For each $\alpha\in \cL(E)$ with operator norm $\|\alpha\|<1$,
we have $\id_E-\alpha\in \GL(E)$,
%
\begin{equation}\label{neumann}
(\id_E-\alpha)^{-1}\, =\;\;  \sum_{k=0}^\infty \, \alpha^k \,,
\end{equation}
and
%
\begin{equation}\label{geomse}
\|(\id_E-\alpha)^{-1}\| \;\leq\;
\frac{1}{1-\|\alpha\|}\,.
\end{equation}
If $\K$ is an ultrametric here and the norm $\|.\|$
on~$E$ is ultrametric, then the set
\[
\Omega\;:=\;
\{\id_E-\alpha\colon \mbox{$\alpha\in \cL(E)$ such that $\|\alpha\|<1$}\}
\]
is an open subgroup of $\GL(E)$
and each $\alpha\in \Omega$ is an isometry.
As a consequence, in the ultrametric case also the set
\[
\Iso(E)\; :=\; \{\alpha\in \GL(E)\colon
\mbox{$(\forall u\in E)$ $\|\alpha(u)\|=\|u\|$}\}
\]
of surjective linear isometries is an open subgroup
of $\GL(E)$.
\end{prop}
\begin{proof}
Since $\|\alpha^k\|\leq \|\alpha\|^k$
and $\cL(E)$ is complete,
the Neumann series $\sum_{k=0}^\infty \alpha^k$
converges for all $\alpha\in \cL(E)$ such that
$\|\alpha\|<1$.
Then $(\id_E-\alpha)\sum_{k=0}^\infty \alpha^k=\id_E
=\big(\sum_{k=0}^\infty \alpha^k\big)(\id_E-\alpha)$,
showing that (\ref{neumann}) holds.
The convergence of the Neumann series being uniform
on the set $\{\alpha\in \cL(E)\colon \|\alpha\|<\frac{1}{2}\}$,
we see that $\iota$ is continuous on an identity
neighborhood and hence continuous
(since $\GL(E)$ is a topological monoid).
The calculation
\[
\|(\id_E-\alpha)^{-1}\|\, =\, \Big\|\sum_{k=0}^\infty\alpha^k\Big\|
\, \leq\, \sum_{k=0}^\infty \|\alpha\|^k \, =\, \frac{1}{1-\|\alpha\|}
\]
establishes (\ref{eqgaug}).
In the ultrametric case,
given $\alpha\in \cL(E)$ with $\|\alpha\|<1$
we have $\|\alpha(x)\|<\|x\|$ and hence $\|(\id_E-\alpha).x\|=\|x\|$
for all $0\not=x\in E$, whence $\id_E-\alpha$
is an isometry.
Furthermore, $\id_E -\alpha$ is invertible by the preceding,
with inverse $(\id_E-\alpha)^{-1}= \id_E -(-\sum_{k=1}^\infty \alpha^k)
\in \Omega$
as $\|-\sum_{k=1}^\infty \alpha^k\|\leq\max\{\|\alpha^k\|\colon k\in \N\}
<1$.
Given $\id_E-\alpha, \id_E-\beta \in \Omega$,
we have $(\id_E-\alpha)\circ
(\id_E-\beta)=\id_E-(\alpha+\beta-\alpha\circ\beta)
\in \Omega$. Thus $\Omega$ is an open subgroup
of $\GL(E)$. Since $\Omega\sub \Iso(E)$,
also $\Iso(E)$ is an open subgroup.
\end{proof}
Actually, $\GL(E)$ is
a Lie group, by \cite[Proposition~2.2]{ARB}.
%
%
\begin{la}\label{operatorpullback}
Let $E_1$, $E_2$ and $F$ be  topological
vector spaces over a valued field
and $\lambda\colon E_1\to E_2$
be a continuous linear map.
Then the mappings
\[
\cL(\lambda, F)\colon \cL(E_2,F)\to\cL(E_1, F)\,,\quad
\alpha\mto \alpha\circ \lambda
\]
and
\[
\cL(F,\lambda) \colon \cL(F,E_1)\to\cL(F, E_2)\,,\quad
\alpha\mto \lambda\circ \alpha
\]
are continuous and linear.
\end{la}
\begin{proof}
Clearly both maps are linear.
If $B\sub E_1$ is a bounded set
and $U\sub F$ a $0$-neighborhood,
then $\lambda(B)\sub E_2$ is bounded
and $\cL(\lambda,F)(\lfloor \lambda(B),U\rfloor )\sub
\lfloor B,U\rfloor$.
Hence $\cL(\lambda, F)$ is continuous at~$0$
and hence continuous.

If $B\sub F$ is bounded and $U\sub E_2$ a $0$-neighborhood,
then $\lambda^{-1}(U)$ is a $0$-neighborhood in~$E_1$
and $\cL(F,\lambda)(\lfloor B,\lambda^{-1}(U)\rfloor)\sub
\lfloor B, U\rfloor$.
Hence $\cL(F,\lambda)$ is continuous at~$0$
and hence continuous.
\end{proof}
\section{\hspace*{-2.6mm}Lipschitz continuous and Lipschitz
differentiable
\hspace*{-2.6mm}mappings}\label{seclip}
%
%
In this section,
we set up our terminology and
prove basic facts
concerning
Lipschitz conditions, Lipschitz continuity
and Lipschitz differentiability.
\subsection*{Lipschitz conditions and Lipschitz continuity}
We first consider functions between
normed spaces that
satisfy a global Lipschitz condition.
\begin{defn}
Let $E$ and $F$ be normed spaces
over a valued field~$\K$
and $U\sub E$ be a subset.
We say that a map $f\colon U\to F$
is \emph{Lipschitz} if
there exists $L\in [0,\infty[$
such that
%
\begin{equation}\label{eqlipglob}
\|f(x)-f(y)\|\,\leq\, L\,\|x-y\|\qquad
\mbox{for all $x,y\in U$.}
\end{equation}
In this case, we define
%
\begin{equation}\label{defLipf}
\Lip(f)\, :=\,\sup\left\{
\frac{\|f(x)-f(y)\|}{\|x-y\|}\colon \mbox{$x,y\in U$, $x\not=y$}
\right\}
\,\in \, [0,\infty[\, .
\end{equation}
Thus $\Lip(f)$ is the smallest possible choice
for the Lipschitz constant~$L$.
\end{defn}
%
%
\begin{la}\label{sweela}
Let $E$ and $F$ be normed spaces
over a valued field~$\K$, $U\sub E$ be a subset
with dense interior,
and $f\colon U\to E$ be a mapping
which is $C^1$ and Lipschitz.
Then $\|f'(x)\|\leq \Lip(f)$
for each $x\in U$.
\end{la} 
\begin{proof}
Given $x\in U^0$ and $y\in E$,
we have
\begin{eqnarray*}
\|f'(x).y\| &= &\lim_{t\to 0}\|t^{-1}(f(x+ty)-f(x))\|
\; =\;  \lim_{t\to 0}|t|^{-1}\|f(x+ty)-f(x)\|\\
&\leq &
\Lip(f)\|y\|
\end{eqnarray*}
because $\|f(x+ty)-f(x)\|\leq \Lip(f)\|ty\|=|t|\Lip(f)\|y\|$.
Hence $\|f'(x).y\|=\|df(x,y)\|\leq \Lip(f)\|y\|$
for all $x\in U$, because $df$
is continuous and $U^0$ is dense
in~$U$.
As a consequence,
$\|f'(x)\|=\sup\{\|f'(x).y\|/\|y\|\colon 0\not=
y\in E\}\leq \Lip(f)$.
\end{proof}
While the Lipschitz maps just
introduced satisfy a quite restrictive
condition, we shall 
use the term ``Lipschitz continuity''
for a much weaker property,
which amounts to a local Lipschitz
condition in the case of mappings between
normed spaces.
%
%
\begin{defn}\label{deflipcts}
Let $E$ and $F$
be topological vector spaces over a valued field~$\K$,
and $U\sub E$ be a subset.
A map $f\colon U\to F$
is called \emph{Lipschitz continuous}
if, for every $x_0\in U$ and gauge $\gamma$
on $F$, there exists a gauge $\eta$
on $E$ and $\delta>0$ such that
\[
\|f(y)-f(x)\|_\gamma \; \leq \; \|y-x\|_\eta \quad
\mbox{for all $\, x,y\in B_\delta^\eta(x_0)\cap U$.}
\]
\end{defn}
After replacing $\eta$ with a suitable multiple,
we may always assume that $\delta=1$.
If~$U$ is open, we may also assume that $B^\eta_\delta(x_0)\sub U$,
whenever this is convenient.
\begin{rem}
For example,
every continuous linear map is
Lipschitz continuous, by Lemma~\ref{morebasstr}\,(c).
\end{rem}
%
%
\begin{la}\label{baseHglob}
Let $E$, $F$ and $H$ be topological
vector spaces over a valued field~$\K$,
$U\sub E$ and $V\sub F$ be subsets,
and $f\colon U\to V\sub F$
and $g\colon V\to H$ be mappings.
Then the following holds:
\begin{itemize}
\item[\rm (a)]
If $f$ is Lipschitz continuous
then $f$ is continuous.
\item[\rm (b)]
If $f$ and $g$ are Lipschitz continuous,
then $g\circ f$ is Lipschitz continuous.
\item[\rm (c)]
If $U$ has dense interior
and $f$ is $C^1$, then $f$ is Lipschitz continuous.
\item[\rm (d)]
If $f$ is Lipschitz continuous
and both $E$ and $F$ are Banach spaces,
then every point $x \in U$ has a neighborhood
$W$ in~$U$ such that $f|_W$ is Lipschitz.
\end{itemize}
\end{la}
\begin{proof}
(a) Given $x_0\in U$ and a gauge $\gamma$ on~$F$,
choose a gauge $\eta$ on~$E$ and $\delta\in \;]0,1]$ as in
Definition~\ref{deflipcts}.
Then $f(U\cap B^\eta_\delta(x_0))\sub \wb{B}_\delta^\gamma(f(x_0))\sub
\wb{B}_1^\gamma(f(x_0))$,
whence $f$ is continuous at $x_0$ by
Lemma~\ref{morebasstr}\,(a).

(b) Let $x\in U$. Given a gauge $\zeta$ on~$H$,
there exists a gauge $\gamma$ on~$F$ and
$\theta> 0$
such that $\|g(z)-g(y)\|_\zeta\leq\|z-y\|_\gamma$
for all $y,z\in V\cap B^\gamma_\theta(f(x))$.
There exists a gauge $\eta$ on~$E$ and $\delta\in
\;]0,\theta]$
such that $\|f(z)-f(y)\|_\gamma\leq \|z-y\|_\eta$
for all $z,y\in U\cap B_\delta^\eta(x)$.
Then $f(y)\in V\cap B^\gamma_\theta(f(x))$
for all $y\in B^\eta_\delta(x)$
and thus $\|g(f(z))-g(f(y))\|_\zeta\leq \|f(z)-f(y)\|_\gamma
\leq\|z-y\|_\eta$ for all $z,y\in B^\eta_\delta(x)$.

(c) We use the first order Taylor expansion
\[
f(x+ty)-f(x)\; =\; tdf(x,y)+t R_1(x,y,t)
\]
of the $C^1$-map
$f\colon E\supseteq U\to F$ (cf.\ \cite[Theorem~5.1]{Ber}).
Here $R_1\colon U^{[1]}\to F$ is a continuous map
and
\[
R_1(x,y,1)=tR_1(x,t^{-1}y,t)\;\, \mbox{for $t\in \K^\times$
and $(x,y)\in U\times E$
such that $x+y\in U$.}
\]
Fix $x_0\in U$. Let $\gamma$ be a gauge on~$F$.
Pick $a\in \K^\times$ such that $|a|<1$,
and a gauge $\eta$ on~$F$ such that
$\|u+v\|_\gamma\leq\max\{\|u\|_\eta,\|v\|_\eta\}$
for all $u,v\in F$.
Since $df(x_0,0)=0$, using the continuity of $df$
we find a gauge $\zeta$ on~$E$
such that
$\|df(x,y)\|_\eta \leq |a|$
for all $x\in B^\zeta_1(x_0)\cap U$ and $y\in B^\zeta_1(0)$,
whence $\|df(x,y)\|_\eta \leq \|y\|_\zeta$
for all $x\in B^\zeta_1(x_0)\cap U$ and $y\in E$
(cf.\ Lemma~\ref{seeRooj}).
Since $R_1(x_0,0,0)=0$, we find a
gauge $\sigma $ on~$E$ and $r \in \,]0,1]$ such that
$\|R_1(x,y,t)\|_\eta  \leq 1$
for all $x\in B^\sigma_r(x_0)$,
$y\in B^\sigma_r(0)$ and $t\in B_r(0)\sub \K$
such that $(x,y,t)\in U^{[1]}$;
we may assume that $\sigma \geq \zeta$.
Let $\tau$ be a gauge on~$E$
such that $\tau\geq |a|^{-1}r^{-1}\sigma$
and $\|u+v\|_\sigma\leq\max\{\|u\|_\tau,\|v\|_\tau\}$
for all $u,v\in E$.
Define $\delta:=\frac{1}{2}r^2|a|$.
Given $x,y\in B^\tau_\delta(x_0)\cap U$,
set $z:=y-x$.
If $\|z\|_\sigma>0$, there is $k\in \Z$ such that
$|a|^{k+1}\leq r^{-1}\|z\|_\sigma <|a|^k$.
Then $\|a^{-k}z\|_\sigma <r$
and $|a^k|\leq |a|^{-1}r^{-1}\|z\|_\sigma \leq |a|^{-1}r^{-1}\delta
< r$,
whence $\|R_1(x,z,1)\|_\eta=|a^k|\, \|R_1(x,a^{-k}z,a^k)\|_\eta
\leq |a^k|\leq |a|^{-1}r^{-1}\|z\|_\sigma\leq \|z\|_\tau$
and thus
$\|f(y)-f(x)\|_\gamma=
\|f(x+z)-f(x)\|_\gamma =\|df(x,z)+R_1(x,z,1)\|_\gamma
\leq \max\big\{\|df(x,z)\|_\eta,
\|R_1(x,z,1)\|_\eta\big\}
\leq \|z\|_\tau$. Hence
\begin{equation}\label{fn2}
\|f(y)-f(x)\|_\gamma \; \leq \; \|y-x\|_\tau \,.
\end{equation}
If $\|z\|_\sigma =0$, given $\ve>0$ pick $t\in \K^\times$
such that $|t|<\min\{r,\ve\}$.
Then $\|df(x,z)\|_\eta=0$
and $\|R_1(x,z,1)\|_\eta=|t|\, \| R_1(x,t^{-1}z,t)\|_\eta\leq |t|\leq \ve$,
whence $\|R_1(x,z,1)\|_\eta =0$ (as $\ve$ was arbitrary).
Thus (\ref{fn2}) also holds if $\|z\|_\sigma=0$.

(d) Let $\|.\|_E$ and $\|.\|_F$ be the norms
on the normed spaces
$E$ and $F$, respectively.
Given $x_0\in U$, by Lipschitz continuity
there exists a neighborhood
$W\sub U$ of $x_0$ and a gauge $\gamma$ on~$E$
such that $\|f(z)-f(y)\|_F\leq \|z-y\|_\gamma$
for all $z,y\in W$.
Since $\{\|.\|_E\}$ is a fundamental system
of gauges for~$E$,
Remark~\ref{morefund}
provides $L>0$ such that $\gamma\leq L\|.\|_E$.
Thus $\|f(z)-f(y)\|_F\leq L\|z-y\|_E$
for all $z,y\in W$,
showing that $f|_W$ is Lipschitz.
\end{proof}
\subsection*{Lipschitz differentiable maps}
We now strengthen the $C^k$-property
by imposing Lipschitz continuity
of the extended difference quotient
maps. The $LC^k$-maps
so obtained are
valuable, for example, in the context
of $p$-adic differential equations.
They are also used in \cite{INV}
to study the parameter dependence
of stable manifolds.
\begin{defn}
Let $\K$ be a valued field,
$E$ and $F$ be topological $\K$-vector spaces,
$U\sub E$ be a subset
with dense interior,
and $f\colon U\to F$ be a mapping.
Let $k\in \N_0\cup\{\infty\}$.
We say that $f$ is \emph{$k$ times Lipschitz differentiable}
(or an \emph{$LC^k$-map})
if $f$ is $C^k$
and $f^{[j]}\colon U^{[j]}\to F$
is Lipschitz continuous
for all $j\in \N_0$
such that $j\leq k$
(where $f^{[0]}:=f$).
\end{defn}
Note that the Lipschitz continuity of $f^{[j]}$
is automatic for $j<k$,
by Lemma~\ref{baseHglob}\,(c).
In particular,
$f$ is $LC^\infty$ if and only if
$f$ is $C^\infty$.
Also note that
$LC^0$-maps are precisely
the Lipschitz continuous maps
(on subsets with dense interior).
%
%
\begin{prop}\label{chainLip}
Let $\K$ be a valued field,
$E$, $F$ and $H$ be topological $\K$-vector spaces,
$U\sub E$ and $V\sub F$ be subsets with dense interior,
and $f\colon U\to V$ as well as $g\colon V\to H$ be
$LC^k$-maps.
Then also $g\circ f\colon U\to H$
is $LC^k$.
\end{prop}
\begin{proof}
In \cite{Ber},
so-called ``$\cC^0$-concepts''
were introduced,
dealing (in particular)
with mappings between
open subsets of Hausdorff topological vector spaces.
As mentioned before, this idea can directly
be adapted to mappings between subsets
with dense interior
(cf.\ also Appendix~\ref{appC0}).
It is clear that Lipschitz continuous maps
between subsets with dense interior of Hausdorff topological
$\K$-vector spaces define a $\cC^0$-concept,
and that
$LC^k$-maps are precisely the $\cC^k$-maps
with respect to this $\cC^0$-concept.
We now use that compositions of $\cC^k$-maps
are always $\cC^k$
(cf.\ \cite[Proposition~4.5]{Ber}).
\end{proof}
%
%
\begin{la}\label{abllip}
Let
$E$ and $F$ be topological vector spaces over a complete
valued field~$\K$
and $f\colon U\to F$ be an $LC^1$-map on a subset $U\sub E$
with dense interior.
Then the following holds:
\begin{itemize}
\item[\rm (a)]
The map $f'\colon U\to \cL(E,F)$, $x\mto f'(x)=df(x,\sbull)$
is Lipschitz continuous and hence continuous.
\item[\rm (b)]
For each gauge $\gamma$ on~$F$
and $x\in U$,
there exists a gauge $\xi$ on~$E$
such that
%
\begin{equation}\label{stronga}
\|f'(z)-f'(y)\|_{\gamma,\xi}\;\leq\; \|z-y\|_\xi \quad
\mbox{for all $z,y\in B_1^\xi(x)\cap U$.}
\end{equation}
\end{itemize}
\end{la}
\begin{proof}
(a) Given $x\in U$, let
$\gamma$ be a gauge on~$F$ and $B\sub E$ be a bounded
set. Since $f^{[1]}$ is Lipschitz continuous,
also its partial map $df\colon U\times E\to F$
is Lipschitz continuous.
As a consequence, there exists a gauge $\eta$ on~$E$
such that
%
\begin{equation}\label{ggaug}
\|df(z,v)-df(y,u)\|_\gamma\leq
\max\{\|z-y\|_\eta,\|v-u\|_\eta\}
\end{equation}
for all $(y,u)$, $(z,v)\in (U\cap B_1^\eta(x))
\times B_1^\eta(0)$.
By boundedness of~$B$,
there exists $t\in \K^\times$ such that $tB\sub B_1^\eta(0)$.
We may assume that $|t|\leq 1$.
Then also $\zeta:=|t|^{-1}\eta$ is a gauge
on~$E$.
Given $z,y\in U\cap B_1^\zeta(x)\sub B_1^\eta(x)$
and $u\in B$, we then have
$\|f'(z).u-f'(y).u\|_\gamma
=|t^{-1}|\,\|f'(z).tu-f'(y).tu\|_\gamma\leq
|t^{-1}|\,\|z-y\|_\eta=\|z-y\|_\zeta$.
Hence
$\|f'(z)-f'(y)\|_{\gamma, B}\leq \|z-y\|_\zeta$
for all $z,y\in U\cap B^\zeta_1(x)$.
As the gauges $\|.\|_{\gamma,B}$
form a fundamental system,
we deduce that $f'$ is Lipschitz 
continuous.

(b)
Taking $u=v$ in (\ref{ggaug}),
we see that
%
\begin{equation}\label{dfbr}
\|df(z,u)-df(y,u)\|_\gamma\; \leq \; \|z-y\|_\eta
\end{equation}
for all $z,y\in U\cap B^\eta_1(x)$ and
all $u\in B_1^\eta(0)$.
Pick $a \in \K^\times$ such that $|a|<1$.
Then (\ref{stronga}) is satisfied
with $\xi:=|a|^{-1}\eta$.
This follows from (\ref{dfbr}) and Lemma~\ref{seeRooj}.
\end{proof}
The following technical lemma will
be needed
when we consider families of contractions
of a Banach space
and study differentiable dependence
of fixed points on parameters.
%
%
\begin{la}\label{monstrtech}
Let $\K$ be a valued field,
$E$ and $F$ be topological
vector spaces over~$\K$ and
$f\colon U\to F$ be an $LC^1$-map
on a subset $U\sub E$ with dense
interior.
Let $x_0\in U$, $y_0\in E$
and $\gamma$ be a gauge on~$F$.
Then there exists
a gauge $\eta$ on~$E$
with the following property:
For each $\ve>0$,
there is $\delta>0$ such that
%
\begin{equation}\label{eqmonstrt}
\|f^{[1]}(x,y_1,t)-f^{[1]}(x,y_2,t)-f'(x_0).(y_1-y_2)\|_\gamma
\;\leq\; \ve\, \|y_1-y_2\|_\eta
\end{equation}
for all elements $x\in B_\delta^\eta(x_0)\cap U$,
$y_1,y_2\in B_\delta^\eta(y_0)$ and $t\in B_\delta(0)\sub \K$
satisfying $(x,y_1,t), (x,y_2,t)\in U^{[1]}$.
\end{la}
\begin{proof}
We shall prove the following stronger
assertion, from which the Lemma readily
follows by taking $\delta:=\min\{\rho, \frac{\ve}{C}\}$:
\begin{numba}\label{evwonum}
\emph{There exists a gauge $\eta$
on~$E$ and $\rho, C>0$ such that}
%
\begin{eqnarray}
\lefteqn{\|f^{[1]}(x,y_1,t)-f^{[1]}(x,y_2,t)-f'(x_0).
(y_1-y_2)\|_\gamma}\notag \\
&\leq &
C\cdot
\max\big\{\|y_1-y_0\|_\eta, \|y_2-y_0\|_\eta,
\|x-x_0\|_\eta, |t|\big\}\cdot
\|y_1-y_2\|_\eta \label{evworse}
\end{eqnarray}
\emph{for all elements $x\in B_\rho^\eta(x_0)\cap U$,
$y_1,y_2\in B_\rho^\eta(y_0)$ and $t\in B_\rho(0)\sub \K$
satisfying $(x,y_1,t), (x,y_2,t)\in U^{[1]}$.}
\end{numba}
Let $\zeta$ be a gauge on~$F$ such that
$\|u+v\|_\gamma\leq \max\{\|u\|_\zeta,\|v\|_\zeta\}$.
Given any elements
$x\in U$,
$y_1,y_2 \in E$ and $t\in \K$
such that $(x,y_1,t), (x,y_2,t)\in U^{[1]}$,
we have
%
%
%
\begin{eqnarray}
\lefteqn{\|f^{[1]}(x,y_1,t)-f^{[1]}(x,y_2,t)
-f'(x_0).(y_1-y_2)\|_\gamma}\quad\notag\\
&=&
\|f^{[1]}(x+ty_2,y_1-y_2,t)-f'(x_0).(y_1-y_2)\|_\gamma\notag\\
&=&
\|f^{[1]}(x+ty_2,y_1-y_2,t)-f'(x+ty_2).(y_1-y_2)\notag\\
& & \qquad\qquad
+(f'(x+ty_2)-f'(x_0)).(y_1-y_2)\|_\gamma\notag\\
&\leq &\max\big\{
\|f^{[1]}(x+ty_2,y_1-y_2,t)-f^{[1]}(x+ty_2, y_1-y_2,0)\|_\zeta,\notag\\
& &\qquad\qquad
\|(f'(x+ty_2)-f'(x_0)).(y_1-y_2)\|_\zeta \big\}\,,\label{spltprb}
\end{eqnarray}
using Lemma~\ref{firsiden} to obtain the first
equality.\\[3mm]
By Lemma~\ref{abllip}\,(b),
there exists a gauge $\xi$ on~$E$ such that
\begin{equation}\label{sml1}
\|f'(x_1)-f'(x_2)\|_{\zeta,\xi}\;\leq\; \|x_1-x_2\|_\xi
\quad\mbox{for all $x_1,x_2\in B^\xi_1(x_0)\cap U$.}
\end{equation}
There exists $r\in \;]0,1]$
and a gauge $\kappa$ on~$E$
such that $\|u+v\|_\xi\leq \max\{\|u\|_\kappa,\|v\|_\kappa\}$
for all $u,v\in E$,
and such that
$a+sb\in B_1^\xi(x_0)$
for all $a\in B^\kappa_r(x_0)$, $b\in B_r^\kappa(y_0)$
and $s\in B^\K_r(0)$.
As $\kappa$ is upper semicontinuous, after
shrinking $r$ we may assume that
$\kappa(b)<\kappa(y_0)+1=:C_1$
for all $b\in B^\kappa_r(y_0)$.
Then
%
\begin{eqnarray}
\hspace*{-10mm}
\lefteqn{\|(f'(x+ty_2)-f'(x_0)).(y_1-y_2)\|_\zeta}\qquad\notag\\
&\leq & \|f'(x+ty_2)-f'(x_0)\|_{\zeta,\xi}\|y_1-y_2\|_\xi\notag\\
&\leq & \|x-x_0+ty_2\|_\xi\|y_1-y_2\|_\xi  \notag\\
&\leq & \max\{\|x-x_0\|_\kappa,\|y_2\|_\kappa|t|\}
\cdot\|y_1-y_2\|_\kappa\notag\\
&\leq & C_1 \max\{\|x-x_0\|_\kappa, |t|\} \cdot
\|y_1-y_2\|_\kappa \label{sml2}
\end{eqnarray}
whenever $x\in B_r^\kappa(x_0)$, $y_2\in B^\kappa_r(y_0)$
and $t\in B^\K_r(0)$
in the above situation.
Hence, we have established estimates of the desired
form for the second term in (\ref{spltprb}).\\[3mm]
By Lipschitz continuity of $f^{[1]}$,
there exists a gauge $\theta\geq \kappa$ on~$E$,
$L>0$
and $\sigma \in \;]0,r]$ such that
\[
\|f^{[1]}(x_1,u_1,t_1)-f^{[1]}(x_2,u_2,t_2)\|_\zeta
\;\leq\;
L \max\{\|x_1-x_2\|_\theta,\|u_1-u_2\|_\theta,|t_1-t_2|\}
\]
for all $x_1,x_2\in B_\sigma^\theta(x_0)\cap U$,
$u_1,u_2\in B^\theta_\sigma(0)$, and $t_1,t_2\in B^\K_\sigma(0)$.
We choose $a\in \K^\times$ such that $|a|<1$.
There exists
a gauge $\eta$ on~$E$ and $\rho\in \;]0,\sigma]$
such that
$\|u+v\|_\theta\leq\max\{\|u\|_\eta,\|v\|_\eta\}$
for all $u,v\in E$
and
%
\begin{equation}\label{nwcond}
B_\rho^\eta(x_0)+B^\K_\rho(0)B_\rho^\eta(y_0)\; \sub\;
B_\sigma^\theta(x_0)\,.
\end{equation}
To see that $\eta$ and $\rho$ have
the desired properties,
let $x$, $y_1$, $y_2$ and $t$ be as
described in {\bf\ref{evwonum}}.
If $\|y_1-y_2\|_\theta\not=0$,
we let $k\in \N_0$ be the unique
element such that $|a|^{k+1}\leq \rho^{-1}\|y_1-y_2\|_\theta
<|a|^k$.
If $\|y_1-y_2\|_\theta=0$,
we let $k\in \N_0$ be arbitrary.
In either case, we abbreviate $s:=a^k$.
Then $|s|\leq 1$.
Using the difference quotient
identity (\ref{bracket1})
from Lemma~\ref{symmetries},
we obtain the following
estimates for the first term
in (\ref{spltprb}):
%
%
\begin{eqnarray}
\lefteqn{\|f^{[1]}(x+ty_2,y_1-y_2,t)-
f^{[1]}(x+ty_2, y_1-y_2,0)\|_\zeta}\notag\\
&=& |s|\cdot \|f^{[1]}(x+ty_2,s^{-1}(y_1-y_2),st)-
f^{[1]}(x+ty_2, s^{-1}(y_1-y_2),0)\|_\zeta\notag \\
&\leq & |s| L |st|\, =\, |s|^2 L|t|\, \leq \, |s|^2 L
\,=\, |a|^{2k}L \label{yeatg}
\end{eqnarray}
because $\|s^{-1}(y_1-y_2)\|_\theta <\rho$
by definition of~$s$,
and $x+ty_2\in B_\sigma^\theta(x_0)$ by (\ref{nwcond}). 
If $\|y_1-y_2\|_\theta\not=0$,
then
\begin{eqnarray*}
|a|^{2k}L & \leq &
L |a|^{-2}\rho^{-2}
\|y_1-y_2\|_\theta^2\\
&\leq&
L |a|^{-2}\rho^{-2}
\max\{\|y_1-y_0\|_\eta,\|y_2-y_0\|_\eta\}\cdot\|y_1-y_2\|_\theta
\end{eqnarray*}
and thus
%
%
\begin{eqnarray}
\hspace*{-1cm}\lefteqn{\|f^{[1]}(x+ty_2,y_1-y_2,t)-
f^{[1]}(x+ty_2, y_1-y_2,0)\|_\zeta}\qquad \notag\\
&\leq &
C_2
\max\{\|y_1-y_0\|_\eta,\|y_2-y_0\|_\eta\}\cdot\|y_1-y_2\|_\eta\label{grpsec}
\end{eqnarray}
with $C_2:=L |a|^{-2}\rho^{-2}$.
If $\|y_1-y_2\|_\theta=0$, then $k$
can be chosen arbitrarily large in
(\ref{yeatg}).
Hence
$\|f^{[1]}(x+ty_2,y_1-y_2,t)-
f^{[1]}(x+ty_2, y_1-y_2,0)\|_\zeta=0$, and thus
(\ref{grpsec}) also holds in this case.
Since $\rho\leq r$ and $\eta\geq \kappa$,
we deduce from (\ref{sml2})
and (\ref{grpsec}) that
(\ref{spltprb}) holds,
with $C:=\max\{C_1,C_2\}$.
\end{proof}
\section{Strictly differentiable mappings}\label{secstrict}
%
%
In this section, we discuss a second
class of maps which are $k$-times
differentiable in a stronger sense
then mere $C^k$-maps,
namely $k$ times strictly differentiable
maps ($SC^k$-maps).
These $SC^k$-maps resemble
to some extent the familiar
continuously Fr\'{e}chet differentiable
mappings between real Banach spaces
(for example, it is known that a map
between real Banach spaces is once strictly
differentiable if and only if it
is continuously Fr\'{e}chet
differentiable, see \cite{Car}).
The $SC^k$-property is weaker than
the $LC^k$-property,
but yet sufficiently
strong for many purposes.
For example, because mere $C^1$-maps
need not be approximated well
enough by their linearization
around a given point,
we shall not be able to prove inverse
function theorems for $C^1$-maps
in general (only for $C^k$-maps with $k\geq 2$,
or in the presence of locally compactness).
Strict differentiability, by contrast,
provides exactly the quality of
approximation needed to make
the construction of inverse functions work.
\subsection*{Definition of strictly differentiable maps}
Before we define strictly
differentiable maps
in general, let us
consider the simpler special case
of mappings between normed spaces.
\begin{defn}\label{defstrictspec}
Let $\K$ be a valued field,
$E$ and $F$ be normed $\K$-vector spaces,
$U\sub E$ be a subset with dense interior,
and $f\colon  U\to F$ be a map.
Given $x\in U$,
we say that $f$ is \emph{strictly differentiable
at $x$}
if there exists a continuous linear map
$f'(x)\in \cL(E,F)$ such that, for every $\ve>0$,
there exists $\delta>0$ such that
%
\begin{equation}\label{prestric}
\|f(z)-f(y)-f'(x).(z-y)\|\; < \; \ve \,\|z-y\|
\end{equation}
for all $y,z\in U$ such that $\|z-x\|<\delta$ and $\|y-x\|<\delta$.
The map $f$ is called \emph{strictly
differentiable} if it is strictly differentiable
at each $x\in U$.
\end{defn}
%
%
\begin{rem}\label{strmeanLip}
Clearly strict differentiability at~$x$
implies total differentiability at~$x$
in the conventional sense (fixing $z=x$).
But it is a stronger
condition, as we are even allowed
to let two elements $z$ and $y$ pass to~$x$
simultaneously.
\end{rem}
%
%
\begin{rem}\label{inpetlip}
It is illuminating to interpret
strict differentiability in terms of
Lipschitz conditions.
Writing $\tilde{f}\colon U\to F$,
$\tilde{f}(y):=f(y)-f(x)-f'(x).(y-x)$,
we have $f(y)=f(x)+f'(x).(y-x)+\tilde{f}(y)$,
i.e., $\tilde{f}$ is the remainder
term of the affine linear approximation
(first order Taylor expansion)
of~$f$ at~$x$.
Strict differentiability at~$x$ means
that, for each $\ve>0$,
we can find $\delta>0$ such that
$\tilde{f}|_{B_\delta(x)\cap U}$
is a Lipschitz map with
$\Lip(\tilde{f}|_{B_\delta(x)\cap U})\leq \ve$.
To see this, note that
the left hand side of (\ref{prestric})
can be written as $\|\tilde{f}(z)-\tilde{f}(y)\|$.
\end{rem}
We now state the
appropriate
generalization of Definition~\ref{defstrictspec}
for mappings between arbitrary
topological vector spaces over
valued fields.
%
%
\begin{defn}\label{defstrict}
Let $\K$ be a valued field,
$E$ and~$F$ be topological $\K$-vector spaces,
$U\sub E$ be a subset with dense interior,
and $f\colon  U\to F$ be a map.
Given $x\in U$,
we say that $f$ is \emph{strictly differentiable
at $x$}
if there exists a continuous linear map
$f'(x)\in \cL(E,F)$
such that, for each gauge $\|.\|_\gamma$ on~$F$
there exists a gauge $\|.\|_\eta$ on~$E$
with the following property:
For each
$\ve>0$, there exists $\delta>0$
such that
%
\begin{equation}\label{defeqstr}
\|f(z)-f(y)-f'(x).(z-y)\|_\gamma < \ve \,\|z-y\|_\eta
\end{equation}
for all $y,z\in U$ such that $\|z-x\|_\eta<\delta$ and
$\|y-x\|_\eta <\delta$.
The map $f$ is called \emph{strictly
differentiable} if it is strictly differentiable
at each $x\in U$.
\end{defn}
\begin{la}
$f'(x)$ is uniquely determined
in the preceding situation.
\end{la}
\begin{proof}
Suppose that also $\alpha \in \cL(E,F)$
satisfies the property of $f'(x)$.
If $\alpha \not=f'(x)$,
we find $u\in E$ such that
$f'(x).u\not=\alpha(u)$,
whence $\|\alpha(u)-f'(x).u\|_\zeta \not=0$
for some gauge $\|.\|_\zeta$
on~$F$.
By Lemma~\ref{substitut},
there exists a gauge
$\|.\|_\gamma$ on~$F$
such that $\|v+w\|_\zeta \leq\max\{\|v\|_\gamma ,\|w\|_\gamma \}$
for all $v,w\in F$.
Choose a gauge
$\|.\|_\eta $ on~$E$ as in Definition~\ref{defstrict},
which works for both $f'(x)$
and~$\alpha$.
Pick $\ve>0$ such that
$\ve\|u\|_\eta <\|\alpha(u)-f'(x).u\|_\zeta $
and let $\delta>0$ be
such that (\ref{defeqstr}) and its analog
with $\alpha$ in place of~$f'(x)$ hold.
Since $U^0$ is dense in~$E$,
we find $y\in U^0\cap B_\delta^\eta(0)$.
By openness of $U^0\cap B_\delta^\eta(0)$,
there exists $t \in \K^\times$
such that $y+tu\in U^0\cap B_\delta^\eta(0)$.
Then
\begin{eqnarray*}
\hspace*{-3mm}|t|\cdot \|\alpha(u)-f'(x).u\|_\zeta & =&
\|\alpha(t u)-f'(x).tu\|_\zeta\\
&\leq &
\max\{\|f(y+tu)-f(y)-f'(x).tu)\|_\gamma ,\\
& &\qquad
\|f(y+tu)-f(y)-\alpha(tu)\|_\gamma \}\\
&\leq &\ve\, \|tu\|_\eta
\, = \, |t|\, \ve\,\|u\|_\eta
\end{eqnarray*}
and thus $\|\alpha(u)-f'(x).u\|_\zeta \leq \ve\|u\|_\eta$,
contradicting our choice of~$\ve$.
\end{proof}
\begin{rem}
Of course, equivalently we can use
gauges in any given fundamental systems
$\Gamma_E$ and $\Gamma_F$ of gauges
for~$E$ and $F$ in Definition~\ref{defstrict}.
In particular, if $E$
(resp., $F$) is polynormed,
we may replace $\|.\|_\eta$ (resp., $\|.\|_\gamma$)
in the definition by a continuous
seminorm.
If $(E,\|.\|_E)$ is a normed space,
we can always take $\|.\|_\eta=\|.\|_E$,
and if $(F,\|.\|_F)$ is normed,
we only need to test the condition
for $\|.\|_\gamma =\|.\|_F$.
\end{rem}
\subsection*{Strictly differentiable maps are {\boldmath$C^1$}}
We now verify that strict differentiability
is a stronger differentiability property
than being~$C^1$.\\[2.5mm]
A simple lemma
by Bourbaki and Dieudonn\'{e}~\cite{Die}
(see also \cite[Exercise~3.2\,A\,(b)]{Eng})
will be useful:
%
\begin{la}\label{Dieudonn}
Let $X$ be a topological space,
$X_0\sub X$ be a dense subset
and $f\colon X_0\to Y$ be a continuous map
to a regular topological space~$Y$.
Then~$f$ has a continuous extension to~$X$
if and only if $f$ has a continuous
extension to $X_0\cup\{x\}$ for each $x\in X$.\Punkt
\end{la}
%
%
\begin{la}\label{strictC1}
Let $\K$ be a valued field,
$E$ and
$F$ be topological $\K$-vector spaces,
$U\sub E$ be a subset with dense interior,
and $f\colon  U\to F$ be a strictly
differentiable map.
Then $f$ is $C^1$,
its strict differential is given by
$f'(x)=df(x,\sbull)$ for all $x\in U$,
and the map
$f'\colon  U\to \cL(E,F)$, $x\mto f'(x)$
is continuous.
\end{la}
\begin{proof}
{\em $f'$ is continuous at each $x\in U$.}
To see this, given $x$ let $B\sub E$ be bounded
and $\|.\|_\zeta$ be a gauge on~$F$.
Choose a gauge $\|.\|_\gamma$ on~$F$
such that $\|u+v\|_\zeta\leq\max\{\|u\|_\gamma,\|v\|_\gamma\}$
for all $u,v\in F$.
Pick $\|.\|_\eta$ as in Definition~\ref{defstrict}.
The set $B$ being bounded,
$M_\eta:=\sup\eta(B)$ is finite.
Given $\ve'>0$, set $\ve:=\ve'/(1+M_\eta)$.
Choose $\delta >0 $ such that
(\ref{defeqstr}) holds.
Let $y\in B_\delta^\eta(x)\cap U$ be arbitrary;
we show that
%
\begin{equation}\label{showy}
\|f'(y)-f'(x)\|_{\zeta,B}\;\leq\; \ve'
\end{equation}
(with notation as in Remark~\ref{defgop}).
By strict differentiability of $f$ at~$y$,
there exists a gauge $\xi\geq \eta $ on~$E$
such that, for each $\ve''>0$,
there exists $\rho>0$ such that
\[
\|f(w)-f(v)-f'(y).(w-v)\|_\gamma\;\leq\; \ve''\|w-v\|_\xi
\]
for all $v,w\in B^\xi_\rho(y)\cap U$.
We set $M_\xi:=\sup\xi(B)<\infty$
and choose $\rho$ as before
for $\ve'':= \ve'/(M_\xi+1)$.
There exists $z\in B_\delta^\eta(x)\cap B^\xi_\rho(y)\cap U^0$,
and $t\in \K^\times$ such that
$z+t B\sub B_\delta^\eta(x)\cap B^\xi_\rho(y)\cap U^0$.
For $u\in B$, we then have
\begin{eqnarray*}
\|f'(y).u-f'(x).u\|_\zeta &=& t^{-1}\|f'(y).tu-f'(x).tu\|_\zeta\\
&\leq &t^{-1}\max\big\{
\|f(z+tu)-f(z)-f'(y).tu\|_\gamma,\\
& & \qquad \|f(z+tu)-f(z)-f'(x).tu\|_\gamma\big\}\\
&\leq & t^{-1}\max\big\{ \ve'' \|tu\|_\xi, \ve\|tu\|_\eta\big\}\\
&\leq &
\max\big\{ \ve'' M_\xi , \ve M_\eta\big\}
\; \leq\; \ve'\,.
\end{eqnarray*}
Hence (\ref{showy}) holds.
The continuity of~$f'$ at $x$ follows.\\[3mm]
{\em $f$ is $C^1_\K$.}
Note first that $f$ is continuous.
In fact, given $x\in U$ and a gauge
$\zeta$ on~$F$, we let $\gamma$
be a gauge on~$F$ such that $\|y+z\|_\zeta
\leq\max\{\|y\|_\gamma,\|z\|_\gamma\}$
for all $y,z\in F$.
Choose a gauge~$\|.\|_\eta$
as in Definition~\ref{defstrict}.
After replacing $\|.\|_\eta$ by a larger
gauge, we may assume that $\|f'(x).w\|_\gamma\leq\|w\|_\eta$
for all $w\in E$. Given $\ve>0$, there
is $\delta\in \;]0,\min\{\ve,1\}]$
such that (\ref{defeqstr}) holds.
Then
\begin{eqnarray*}
\|f(y)-f(x)\|_\zeta & \leq & \max\{\|f(y)-f(x)-f'(x).(y-x)\|_\gamma,
\|f'(x).(y-x)\|_\gamma\}\\
& \leq & \max\{\ve \|y-x\|_\eta, \|y-x\|_\eta\}\,
\leq \, \ve
\end{eqnarray*}
for all $y\in B_\delta^\eta(x)\cap U$.
We deduce that~$f$ is continuous.

Next, let $W:=\{(x,y,t)\in U^{[1]}\colon  t\not=0\}$.
Define $g\colon  U^{[1]}\to F$ via
$g(x,y,t):=\frac{1}{t}(f(x+ty)-f(x))$
for $(x,y,t)\in W$,
while we set $g(x,y,0):=f'(x).y$
for\linebreak
$(x,y)\in U\times E$.
Then $g|_W$ is continuous since
$f$ is continuous.
Hence, by Lemma~\ref{Dieudonn},
$g$ will be continuous if we can show that
$g(x_\alpha,y_\alpha,t_\alpha)\to g(x,y,0)$,
for each net $((x_\alpha,y_\alpha,t_\alpha))_{\alpha\in I}$
in~$W$ which converges to some $(x,y,0)\in U^{[1]}$.
To this end, given a gauge $\|.\|_\zeta$ on~$F$
let $\|.\|_\gamma$ be a gauge on~$F$
such that $\|y+z\|_\zeta
\leq \max\{\|y\|_\gamma,\|z\|_\gamma\}$.
Let $\|.\|_\eta$ be a gauge on~$E$ as
in Definition~\ref{defstrict}.
After replacing~$\|.\|_\eta$ by a larger
gauge if necessary, we may assume that
$\|f'(x).w\|_\gamma\leq \|w\|_\eta$
for all $w\in E$.
Given $\ve>0$, let $\delta\in \;]0,\ve]$
be such that (\ref{defeqstr}) holds.
Since $x_\alpha\to x$
and $x_\alpha+t_\alpha y_\alpha\to x$,
we have $x_\alpha \in U\cap B^\eta_\delta(x)$
and
$x_\alpha+t_\alpha y_\alpha\in U\cap B^\eta_\delta(x)$
eventually.
Furthermore, $\|y_\alpha\|_\eta\leq \|y\|_\eta+1$
and $\|y_\alpha -y\|_\eta \leq \ve$
eventually, as $y_\alpha\to y$.
For any such $\alpha$, we obtain
\begin{eqnarray*}
\lefteqn{\|g(x_\alpha,y_\alpha,t_\alpha)-g(x,y,0)\|_\zeta}\\
& = &
\left\|
\frac{f(x_\alpha + t_\alpha y_\alpha)-f(x_\alpha)}{t_\alpha}
-f'(x).y\right\|_\zeta\\
& \leq &
\max\Big\{
\frac{1}{|t_\alpha|}
\|f(x_\alpha + t_\alpha y_\alpha)-f(x_\alpha)
-f'(x).t_\alpha y_\alpha \|_\gamma,\,
\|f'(x).(y_\alpha-y)\|_\gamma\Big\}\\
&\leq &\max\{\ve \|y_\alpha\|_\eta,\,\|y_\alpha-y\|_\eta\}
\,\leq\, \ve(\|y\|_\eta+1)\,,
\end{eqnarray*}
which can be made arbitrarily small.
Thus $g(x_\alpha,y_\alpha,t_\alpha)\to g(x,y,0)$
in~$F$, which completes the proof.
\end{proof}
\subsection*{\boldmath{$LC^1$}-maps are strictly differentiable}
We now show that every Lipschitz differentiable
map is strictly
differentiable. As a consequence, every $C^2$-map
is strictly differentiable.
%
\begin{prop}\label{Lipstrict}
Let $E$ and $F$ be  topological vector spaces
over a valued field~$\K$ and $f\colon U\to F$
be an $LC^1$-map on a subset $U\sub E$
with dense interior.
Then $f$ is $SC^1$, with strict differential
$f'(x)=df(x,\sbull)$ at $x\in U$.
\end{prop}
\begin{proof}
Let $x_0\in U$ and $\|.\|_\gamma$ be a gauge on~$F$.
Abbreviate $f'(x_0):=df(x_0,\sbull)$.
Choose a gauge $\|.\|_\zeta$ on~$F$ such that
$\|v+w\|_\gamma\leq \max\{\|v\|_\zeta,\|w\|_\zeta \}$ for all
$v,w\in F$.
Then
%
\begin{eqnarray}
\|f(z)-f(y)-f'(x_0).(z-y)\|_\gamma
&\leq &
\max\big\{ \|f(z)-f(y)-f'(y).(z-y)\|_\zeta,\notag \\
& & \quad \|(f'(y)-f'(x_0)).(z-y)\|_\zeta\big\}. \label{firsec}
\end{eqnarray}
The first order Taylor remainder
$R_1\colon U^{[1]}\to F$ (as in the proof of Lemma~\ref{baseHglob}\,(c))
being Lipschitz continuous
(cf.\ proof of Proposition~\ref{chainLip}
and \cite[Theorem~5.1]{Ber}),
there exists an open neighborhood
$V\sub U^{[1]}$ of $(x_0,0,0)$,
a gauge $\|.\|_\xi$ on~$E$ and $L > 0$
such that
\begin{eqnarray}
\lefteqn{\|R_1(x_1,y_1,t_1)-R_1(x_2,y_2,t_2)\|_\zeta}\qquad\notag\\
&\leq &L\max\{\|x_1-x_2\|_\xi,\|y_1-y_2\|_\xi,
|t_1-t_2|\} \label{usetaylip}
\end{eqnarray}
for all $(x_1,y_1,t_1), (x_2,y_2,t_2)\in V$.
After replacing $\xi$ by a larger
gauge if necessary, there exists
$r>0$ such that
$U^{[1]}\cap (B_r^\xi(x_0)\times B^\xi_r(0)\times B^\K_r(0))
\sub V$.
After replacing $\xi$
with a larger gauge,
by Lemma~\ref{abllip}\,(b) we may assume that
\[
\|f'(z)-f'(y)\|_{\zeta,\xi}\;\leq\; \|z-y\|_\xi\qquad
\mbox{for all $z,y\in B_r^\xi(x_0)$.}
\]
There exists a gauge $\|.\|_\eta$ on~$E$
such that $\|u+v\|_\xi \leq\max\{\|u\|_\eta,\|v\|_\eta\}$
for all $u,v\in E$.
Pick $a\in \K^\times$ such that $|a|<1$.
Given $\ve>0$,
let $\rho :=\min\{1,r,\ve\}$.
Choose $\delta \in \;\big]
0, \min\{\rho^2|a|,  \frac{|a|^2\ve\rho^2}{L}\} \big[$.
For all $z,y \in B^\eta_\delta(x_0)\cap U$,
we then have
\begin{eqnarray*}
\|(f'(y)-f'(x_0)).(z-y)\|_\zeta
& \leq  & \|f'(y)-f'(x_0)\|_{\zeta,\xi}\|z-y\|_\xi\\
& \leq  & \|y-x_0\|_\xi\|z-y\|_\xi\; \leq \;
\ve\|z-y\|_\xi\,,
\end{eqnarray*}
whence the second term of on the right hand side of
(\ref{firsec})
is no larger than~$\ve$.
We have $\|z-y\|_\xi\leq\max\{\|z-x_0\|_\eta,\|y-x_0\|_\eta\}
<\delta\leq\rho$.
If $\|z-y\|_\xi \not=0$,
let $k\in \N_0$ be the unique
element such that $|a|^{k+1} \leq \frac{\|z-y\|_\xi}{\rho}<|a|^k$.
If $\|z-y\|_\xi=0$, let $k\in \N_0$ be arbitrary.
Set $t:=a^k$.
Then $|t|\leq |a|^{-1}\|z-y\|_\xi \, \rho^{-1}<|a|^{-1}\delta \rho^{-1}
<\rho$. Hence
\begin{eqnarray*}
\|f(z)-f(y)-f'(y).(z-y)\|_\zeta
&=& \|R_1(y,z-y,1)\|_\zeta \;=\;
|t|\cdot\big\|R_1\big(y,{\textstyle \frac{z-y}{t}}, t\big)\big\|_\zeta\\
&=&
|t|\cdot\big\|R_1\big(y,{\textstyle \frac{z-y}{t}}, t\big)
-R_1\big(y,{\textstyle \frac{z-y}{t}}, 0 \big) \big\|_\zeta\\
&\leq& L|t|^2\,,
\end{eqnarray*}
by (\ref{usetaylip}).
If $\|z-y\|_\xi\not=0$, then
\begin{eqnarray*}
L|t|^2 & \leq & L|a|^{-2} \rho^{-2}\|z-y\|_\xi^2
\; \leq \; L|a|^{-2}\rho^{-2}\|z-y\|_\xi\max\{\|z-x_0\|_\eta,\|y-x_0\|_\eta\}\\
& \leq & L|a|^{-2}\rho^{-2}\delta \|z-y\|_\xi \; \leq \;
\ve \|z-y\|_\eta\,.
\end{eqnarray*}
If $\|z-y\|_\xi=0$,
then $L|t|^2=L|a|^{2k}$, where $k$ can be chosen arbitrarily large,
and thus $\|f(z)-f(y)-f'(y).(z-y)\|_\zeta=0\leq \ve\|z-y\|_\eta$
also in this case.
Hence both terms on the right hand side of (\ref{firsec})
are $\leq \ve\|z-y\|_\eta$, and thus
\[
\|f(z)-f(y)-f'(x_0).(z-y)\|_\gamma\;\leq\; \ve\|z-y\|_\eta
\]
for all $y,z\in B^\eta_\delta(x_0)\cap U$.
Thus $f$ is strictly differentiable
at~$x_0$, with $f'(x_0)$ as before.
\end{proof}
%
%
\begin{cor}\label{C2strictnew}
Let $E$ and $F$ be  topological vector spaces
over a valued field and $f\colon U\to F$
be a $C^2$-map on a subset $U\sub E$
with dense interior.
Then $f$ is $LC^1$ and hence also $SC^1$.
\end{cor}
\begin{proof}
Applying Lemma~\ref{baseHglob}\,(c)
to $f$ and $f^{[1]}$, we find that
$f$ is $LC^1$ and hence
$SC^1$, by Proposition~\ref{Lipstrict}.
\end{proof}
\subsection*{Strictly differentiable maps on locally compact
domains}
For mappings on open subsets of
finite-dimensional topological vector spaces
over locally compact topological fields,
the preceding result can be
strengthened:
such a map is $C^1$ if and only
if it is strictly differentiable.
More generally, this conclusion remains
valid for mappings on locally compact domains.
%
%
%
%
%
%
\begin{la}\label{C1impllud}
Let $\,\K$ be a locally
compact field,
$E$ be a finite-dimensional
$\K$-vector space,
$F$ be a topological $\K$-vector space,
$U\sub E$ be a locally compact subset
with dense interior, and
$f\colon  U\to F$ be a map.
Then $f$ is~$C^1$
if and only if $f$ is strictly
differentiable.
\end{la}
\begin{proof}
We already know that every strictly
differentiable map is~$C^1$. Conversely, assume
that $f$ is~$C^1$.
Let $|.|$ be an absolute value on~$\K$ defining its topology,
$\|.\|$ be a norm on~$E$,
and $0\not= a \in \K$ such that $|a|<1$.
Given $x_0\in U$,
let $V\sub U$ be an open neighborhood of~$x_0$
with compact closure $\wb{V}\sub U$.
Define a map
$f'\colon  U\to \cL(E,F)$
via $f'(x):=df(x,\sbull)=f^{[1]}(x,\sbull,0)$.
Given a gauge $\gamma$ on~$F$,
choose gauges $\eta$ and $\zeta$ on~$F$
such that $\|u+v\|_\gamma\leq \|u\|_\eta+\|v\|_\eta$
and $\|u_1+\cdots+ u_n\|_\eta\leq\sum_{i=1}^n\|u_i\|_\zeta$
for all $u,v,u_1,\ldots, u_n\in F$,
where $n:=\dim_\K(E)$.
Given $\ve>0$, consider
the continuous function
\[
g\colon  U^{[1]}\to F,\;\;\;\;
g(x,y,t):=f^{[1]}(x,y,t)-f^{[1]}(x,y,0)\,.
\]
Then $\wb{V}\times \wb{B_{\frac{1}{|a|}}^E(0)}\times\{0\}
\sub U\times E\times \{0\}\sub U^{[1]}$
is a compact subset on which $g$ vanishes identically.
Using a compactness argument, we find $\sigma>0$
such that
$\|g(x,y,t)\|_\eta<\frac{\ve}{2}$
for all $(x,y,t)\in
U^{[1]}\cap
\big(
\wb{V}\times \wb{B_{\frac{1}{|a|}}^E(0)}\times B_\sigma^\K(0)\big)$.
Let $e_1,\ldots,e_n$ be a basis
of~$E$, and $e_1^*,\ldots,e_n^*\in E'$ be its
dual basis.\\[3mm]
Given $\alpha \in \cL(E,F)$,
for each $v\in E$ we have
$\|\alpha(v)\|_\eta
=\|\sum_{i=1}^ne_i^*(v)\alpha(e_i)\|_\eta
\leq
\sum_{i=1}^n |e_i^*(v)|\cdot \|\alpha(e_i)\|_\zeta\leq
\sum_{i=1}^n\|e_i^*\|\cdot \|\alpha(e_i)\|_\zeta\|v\|$.
Thus
\begin{equation}\label{almtriv}
\|\alpha\|_\eta:=\sup\{\|\alpha(v)\|_\eta/\|v\|\colon 0\not=v\in
E\}
\leq\sum_{i=1}^n \|e_i^*\|\cdot \|\alpha(e_i)\|_\zeta
\end{equation}
for all $\alpha\in \cL(E,F)$.
Let $i\in\{1,\ldots, n\}$.
The map $\wb{V}\to F$, $x\mto df(x,e_i)$
being uniformly
continuous, we find $\delta_i>0$
such that $\|df(y,e_i)-df(x,e_i)\|_\zeta<\frac{\ve}{2n\|e_i^*\|}$
for all $x,y\in \wb{V}$ such that $\|x-y\|<\delta_i$.
Define $\delta:=\frac{1}{2}\min\{\sigma,\delta_1,\ldots,\delta_n\}$.
By (\ref{almtriv}) and the choice of~$\delta_i$, we have
$\|df(y,\sbull)-df(x,\sbull)\|_\eta<\frac{\ve}{2}$
for all $x,y\in \wb{V}$ such that $\|x-y\|< 2\delta$.

Let $x,y,z\in V$ be given such that $y\not=z$,
$\|y-x\|<\delta$, and $\|z-x\|<\delta$.
There exists $k\in \Z$ such that
$|a|^{k+1}\leq \|z-y\|<|a|^k$.
We set $s:=a^{k+1}$.
Then $\|\frac{1}{s}(z-y)\|<\frac{1}{|a|}$,
$|s|=|a|^{k+1}\leq \|z-y\|<2\delta\leq\sigma$,
and $\|z-y\|<2\delta$. Thus
\begin{eqnarray*}
\lefteqn{\frac{\|f(z)-f(y)-f'(x).(z-y)\|_\gamma}{\|z-y\|}}\\
&\leq &
\frac{\|f(z)-f(y)-f'(y).(z-y)\|_\eta}{\|z-y\|}
+\frac{\|(f'(y)-f'(x)).(z-y)\|_\eta}{\|z-y\|}\\
& < &
{\textstyle \frac{|s|}{\|z-y\|}\cdot
\left\|\frac{1}{s}(f(z)-f(y))-f'(y).\frac{1}{s}(z-y)\right\|_\eta
+\frac{\ve}{2}}\\
& \leq &
{\textstyle\left\|f^{[1]}\left(y,\frac{1}{s}(z-y),s\right)-
f^{[1]}\left(y,\frac{1}{s}(z-y),0\right)\right\|_\eta
+\frac{\ve}{2}}\\
& = & {\textstyle\left\|g
\left(y,\frac{1}{s}(z-y),s\right)\right\|_\eta+\frac{\ve}{2}
\leq \ve\,.}
\end{eqnarray*}
Hence $f|_V$ is strictly differentiable
at each $x\in V$.
As the proof shows, given $\gamma$
and $\ve$
we can even choose $\delta$ independently of~$x\in V$.
\end{proof}
Also a variant of
Lemma~\ref{C1impllud} involving parameters
will be needed later.
%
%
\begin{la}\label{C1paraunif}
Let $\K$ be a locally
compact topological field
and $|.|$ be an absolute value on~$\K$
defining its topology.
Let $E$ be a finite-dimensional normed
$\K$-vector space, $U\sub E$ be a locally compact subset
with dense interior,
$F$ be a topological $\K$-vector space,
and $P$ be a topological space.
Let $f\colon  P\times U\to F$ be a continuous
map such that $f_p:=f(p,\sbull)\colon  U\to F$
is $C^1$ for all $p\in P$, and such
that the map
\[
P\times U^{[1]}\to F,\quad
(p,y)\mto (f_p)^{[1]}(y)
\]
is continuous.
Let $p\in P$ and $u\in U$ be given.
Then, for every $\ve>0$ and gauge
$\gamma$ on~$F$, there is a neighborhood
$Q$ of~$p$ in~$P$ and $\delta>0$ such that
\[
\|f_q(z)-f_q(y)-f_q'(u).(z-y)\|_\gamma<\ve \|z-y\|
\]
%
%
for all $q\in Q$ and $y,z\in B_\delta(u)\cap U$,
where $f_q'(u):=d(f_q)(u,\sbull)$.
\end{la}
\begin{proof}
Given $\ve>0$ and $\gamma$,
let $\eta$ and $\zeta$ be as in the preceding proof.
Pick $0\not= a \in \K$ such that $|a|<1$.
Let $V\sub U$ be an open neighborhood of~$u$
with compact closure $\wb{V}\sub U$.
Consider the continuous mapping
\[
g\colon  P\times U^{[1]}\to F,\;\;\;\;
g(q,x,y,t):=f_q^{[1]}(x,y,t)-f^{[1]}_q(x,y,0)\,.
\]
Then $\{p\}\times \wb{V}\times \wb{B_{\frac{1}{|a|}}^E(0)}\times\{0\}
\sub P\times U\times E\times \{0\}\sub P\times U^{[1]}$
is a compact subset on which $g$ vanishes identically.
Using a compactness argument,
we find $\sigma>0$ and a neighborhood $P_0$ of~$p$ in~$P$
such that
\[
\|g(q,x,y,t)\|_\eta <{\textstyle \frac{\ve}{2}}\quad
\mbox{for all $(q,x,y,t)\in
P_0\times \big(
U^{[1]}\cap(\wb{V}\times \wb{B_{\frac{1}{|a|}}^E(0)}\times
B_\sigma^\K(0))\big)$.}
\]
Let $e_1,\ldots,e_n$ be a basis
of~$E$, and $e_1^*,\ldots,e_n^*$ be its
dual basis.
Using the compactness of~$\wb{V}$,
we find a neighborhood $Q\sub P_0$ of~$p$
and $\kappa>0$
such that $\|df_q(z,e_i)-df_q(y,e_i)\|_\eta <\frac{\ve}{2n\|e_i^*\|}$
for all $q\in Q$, $i\in \{1,\ldots, n\}$,
and all $y,z\in \wb{V}$
such that $\|z-y\|<\kappa$.
Let $\delta:=\min\{\frac{\sigma}{2},\frac{\kappa}{2}\}$.
Re-using the estimates from the proof of
Lemma~\ref{C1impllud}, we see
that the current assertion
holds for~$Q$ and~$\delta$.
\end{proof}
Recall that on a finite-dimensional
vector space $F$ over a topological field~$\K$,
of dimension~$n$,
there is a unique Hausdorff vector topology
making $F$ isomorphic to the direct product
$\K^n$ as a topological vector space.
It is called the \emph{canonical vector topology
on~$F$}.\\[2.5mm]
In the context of our current discussions,
the following observation is useful.
%
\begin{la}\label{operctsfin}
Let $E$ and $H$ be topological
vector spaces over a valued
field~$\K$, and $F$ be a finite-dimensional
$\K$-vector space, equipped with its canonical
vector topology.
Let $U\sub E$ and $V\sub F$ be subsets
with dense interior and
$f\colon U\times V\to H$ be a $C^1$-map.
Then the map
\[
U\times V\to \cL(F,H)\,,\quad
(x,y)\mto f'_x(y):=df((x,y),(0,\sbull))
\]
is continuous.
\end{la}
\begin{proof}
Let $e_1,\ldots, e_n$ be a basis of~$F$
and $e_1^*,\ldots, e_n^*\in F'$ be its
dual basis, determined by $e_i^*(e_j)=\delta_{ij}$.
Then
\begin{equation}\label{prepafo1}
f'_x(y).w\, =\, \sum_{j=1}^n e_j^*(w) \, df(x,y,0, e_j)
\end{equation}
for all $x\in U$, $y\in V$ and $w\in F$.
The assertion can now easily be derived.
In fact, if $\|.\|_\gamma$ is a gauge on~$H$
and $B\sub F$ a bounded subset,
choose a gauge $\|.\|_\zeta$ on $H$
such that $\|u_1+\cdots+ u_n\|_\gamma\leq\max\{
\|u_j\|_\zeta\colon j=1,\ldots, n\}$.
Pick $C\in \;]0,\infty[$
such that
%
\begin{equation}\label{prepafo2}
|e_j^*(w)| \; \leq\; C \qquad
\mbox{for all $w\in B$ and $j=1,\ldots, n$.}
\end{equation}
By continuity,
for $(x_0,y_0)\in U\times V$
there exist neighborhoods
$U_0\sub U$ of~$x_0$ and $V_0\sub V$ of $y_0$
such that
%
\begin{equation}\label{prepafo3}
\|df(x ,y ,0,e_j)-df(x_0,y_0,0,e_j)\|_\zeta\;\leq\;\frac{1}{C}
\quad\mbox{for all $(x,y)\in U_0\times V_0$.}
\end{equation}
Combining (\ref{prepafo1}),
(\ref{prepafo2}) and
(\ref{prepafo3}), we see that
$\|(f'_x(y)-f'_{x_0}(y_0)).w\|_\zeta\leq 1$
for all $(x,y)\in U_0\times V_0$ and $w\in B$, and thus
$\|f'_x(y)-f'_{x_0}(y_0)\|_{\zeta, B}\leq 1$
(using the notation from Remark~\ref{defgop}).
As a consequence, the map under consideration
is continuous at $(x_0,y_0)$.
\end{proof}
%
%
%
\begin{rem}\label{varC1paraunif}
In the situation of Lemma~\ref{C1paraunif},
we can achieve that furthermore
\[
\|f_q(z)-f_q(y)-f_p'(u).(z-y)\|_\gamma<\ve \|z-y\|
\]
for all $q\in Q$ and $y,z\in B_\delta(u)\cap U$.\\[2.5mm]
Indeed,
given $\ve>0$ and a gauge $\|.\|_\gamma$ on~$F$,
let $\|.\|_\xi$ be a gauge on $F$ such that
$\|u+v\|_\gamma \leq\max\{\|u\|_\xi,\|v\|_\xi\}$
for all $u,v\in F$.
By Lemma~\ref{C1paraunif},
there is a neighborhood
$Q$ of~$p$ in~$P$ and $\delta>0$ such that
\[
\|f_q(z)-f_q(y)-f_q'(u).(z-y)\|_\xi <\ve \|z-y\|
\]
for all $q\in Q$ and $y,z\in B_\delta(u)\cap U$.
By Lemma~\ref{operctsfin},
after shrinking~$Q$ we may assume
that $\|f'_q(u)-f'_p(u)\|_{\xi,\nu}\leq\ve$
for all $q\in Q$,
where $\nu:=\|.\|$ is the norm on~$F$.
Hence
\begin{eqnarray*}
\lefteqn{\|f_q(z)-f_q(y)-f_p'(u).(z-y)\|_\gamma}\qquad\\
&\leq &\max\big\{\|f_q(z)-f_q(y)-f_q'(u).(z-y)\|_\xi,
\|(f'_q(u)-f'_p(u)).(z-y)\|_\xi\big\}\\
&\leq & \max\big\{\ve\|z-y\|, \|f'_q(u)-f'_p(u)\|_{\xi,\nu}\|z-y\|\big\}
\; \leq \; \ve \|z-y\|\,.
\end{eqnarray*}
\end{rem}
\subsection*{Strict differentiability of higher order}\label{higherSCk}
We now define and discuss
$k$ times strictly differentiable mappings
between subsets of topological
vector spaces over valued fields.
\begin{defn}
Let $\K$ be a valued field, $E$ and~$F$ be topological
$\K$-vector spaces,
and $U\sub E$ be a subset with dense interior.
A map $f\colon  U\to F$ is called an \emph{$SC^0$-map}
if it is continuous; it is called an \emph{$SC^1$-map}
is it is strictly differentiable
(and hence $C^1$ in particular).
Inductively, having defined
$SC^k$-maps for some $k\in \N$ (which are $C^k$ in particular),
we call $f$ an \emph{$SC^{k+1}$-map}
if it is an $SC^k$-map and the map
$f^{[k]}\colon  U^{[k]}\to F$ is $SC^1$.
The map $f$ is \emph{$SC^\infty$} if it is
an $SC^k$-map for all $k\in \N_0$.
\end{defn}
\begin{rem}\label{remk1}
In other words, $f$ is $SC^k$ if and only
if $f$ is $C^k$ and $f^{[j]}\colon  U^{[j]}\to F$
is strictly differentiable for all $j\in \N_0$
such that $j<k$.
It follows from this and Remark~\ref{remonesev}
that $f$ is $SC^k$ if
and only if $f$ is $SC^1$ and $f^{[1]}$ is $SC^{k-1}$.
\end{rem}
\begin{rem}\label{remk2}
If $f\colon  E\supseteq
U\to F$ is $C^{k+1}$
in the preceding situation,
then $f$ is an $SC^k$-map.
In fact, for every $j\in \N_0$ such that $j<k$,
the map $f^{[j]}$ is $C^{k+1-j}$,
where $k+1-j\geq 2$.
Thus $f^{[j]}$ is strictly
differentiable, by Corollary~\ref{C2strictnew}.
It is also clear from the definitions
that every $LC^k$ map is $SC^k$, since
every $LC^1$-map is $SC^1$.
Hence, the relations
between the various differentiability
properties can be summarized as follows:
\[
C^{k+1}\; \Longrightarrow\;
LC^k\;
\Longrightarrow\;
SC^k\;
\Longrightarrow
\;
C^k\,.
\]
\end{rem}
\begin{rem}\label{SCk=Ck}
If $\K$ is a locally compact topological
field, then
a mapping from an open subset of a finite-dimensional
$\K$-vector space to a topological $\K$-vector
space is $C^k$ if and only if it is
an $SC^k$-map, by a simple induction based on
Lemma~\ref{C1impllud}
and Remark~\ref{remk1}.
The same conclusion holds for mappings
on locally compact subsets with dense interior.
\end{rem}
Compositions of composable $SC^k$-maps are $SC^k$.
\begin{prop}
Let $\K$ be a valued field, $E$, $F$
and~$H$ be topological
$\K$-vector spaces,
and $U\sub E$, $V\sub F$ be subsets with dense interior.
Let $k\in \N_0 \cup\{\infty\}$
and suppose that $f\colon  U\to V\sub F$
and $g\colon  V\to H$ are $SC^k$.
Then also $g\circ f\colon  U\to H$ is~$SC^k$.
\end{prop}
\begin{proof}
The case $k=0$ is trivial.
The case $k=1$
can be shown as follows:
Given $x\in U$, let $\gamma$ be a gauge
on~$H$.
There exists a gauge $\xi$ on~$H$ such that
$\|u+v\|_\gamma\leq\|u\|_\xi+\|v\|_\xi$
for all $u,v\in H$.
By strict differentiability of $g$ at~$f(x)$,
there exists a gauge $\zeta$ on~$F$
such that,
for each $\ve>0$, there exists
$\theta > 0$
such that
\begin{equation}\label{brrest1}
\|g(z)-g(y)-g'(f(x)).(z-y)\|_\xi \leq \ve \, \|z-y\|_\zeta\quad
\mbox{for all $y,z\in B_\theta^\zeta (f(x))\cap V$.}
\end{equation}
Let $\kappa$ be a gauge on~$F$
such that $\|u+v\|_\zeta\leq\max\{\|u\|_\kappa,\|v\|_\kappa\}$
for all $u,v\in F$.
After increasing $\kappa$ if necessary,
we may assume that $\|g'(f(x))\|_{\xi,\kappa}\leq 1$,
using the notation from Definition~\ref{likeopnorm}.
By strict differentiability of~$f$ at~$x$
and continuity of~$f$ at~$x$,
there exists a gauge $\eta$ on~$E$
such that,
for each $\theta>0$,
there exists $\delta> 0$
such that
\begin{equation}\label{brrest2}
\|f(z)-f(y)-f'(x).(z-y)\|_\kappa \leq \theta \, \|z-y\|_\eta\quad
\mbox{for all $z,y\in B_\delta^\eta(x)\cap U$.}
\end{equation}
After increasing $\eta$,
we may assume that $\|f'(x)\|_{\kappa,\eta}\leq 1$.
Given $\ve>0$,
choose $\theta\in \,]0, \min\{1,\ve\}]$
and $\delta\in \,]0,\theta [$ such that
(\ref{brrest1}) and (\ref{brrest2})
hold.
Given $y,z\in B_\delta^\eta(x)\cap U$,
we have
\begin{eqnarray*}
\|f(z)-f(x)\|_\zeta &\leq &
\max\{
\|f(z)-f(x)-f'(x).(z-x)\|_\kappa,
\|f'(x).(z-x)\|_\kappa\}\\
&\leq & \max\{\theta \|z-x\|_\eta,
\|z-x\|_\eta\}\; < \;
\theta
\end{eqnarray*}
and likewise $\|f(y)-f(x)\|_\zeta<\theta$.
Hence
\begin{eqnarray*}
\hspace*{-8mm}
\lefteqn{\|g(f(z))-g(f(y))
-g'(f(x)).f'(x).(z-y)\|_\gamma}\quad\quad\\
&\leq & \|g(f(z))-g(f(y))-g'(f(x)).(f(z)-f(y))\|_\xi\\
& &\quad\quad +\; \|g'(f(x)).(f(z)-f(y)-f'(x).(z-y))\|_\xi\\
&\leq& \ve\|f(z)-f(y)\|_\zeta+
\|g'(f(x))\|_{\xi,\kappa}\cdot
\|f(z)-f(y)-f'(x).(z-y)\|_\kappa \\
&\leq & \ve \|f(z)-f(y)-f'(x).(z-y)\|_\kappa+\ve\|f'(x).(z-y)\|_\kappa
+\theta\|z-y\|_\eta\\
&\leq & \ve \theta \|z-y\|_\eta
+\ve\|z-y\|_\eta+\theta\|z-y\|_\eta
\; \leq \; 3\ve\|z-y\|_\eta
\end{eqnarray*}
Hence indeed $g\circ f$ is
strictly differentiable at~$x$,
with differential $(g\circ f)'(x)=g'(f(x))\circ f'(x)$.

{\em The case $k\geq 2$.}
Let us call a map between subsets with dense
interior of Hausdorff
topological $\K$-vector spaces~$\cC^0$
if it is $SC^1$.
It is clear from the case $k=1$
that we obtain a $\cC^0$-concept
in the sense of \cite{Ber} in this way
(suitable adapted to non-open sets,
e.g.\ as in Appendix~\ref{appC0}).
Furthermore, $SC^{k+1}$-maps
are precisely the $\cC^k$-maps
for this $\cC^0$-concept,
for each $k\in \N_0\cup\{\infty\}$.
The assertion therefore reduces to the
Chain Rule for $\cC^k$-maps
(cf.\ \cite[Proposition~4.5]{Ber}).
\end{proof}
\section{Dependence of fixed points on parameters}\label{secfp}
%
%
We now study the dependence
of fixed points of contractions on parameters.
In particular, we shall establish
$C^k$-, $SC^k$- and $LC^k$-dependence under
natural hypotheses.
These results will be used in Section~\ref{sectimpl}
to prove generalizations of the inverse- and
implicit function theorems.\\[3mm]
We recall the notion of a contraction.
\begin{defn}
A mapping $f\colon X\to Y$
between metric spaces $(X,d_X)$ and $(Y,d_Y)$
is called a \emph{contraction}
if there exists $\theta\in [0,1[$
(a ``contraction constant'')
such that
\[
d_Y(f(x),f(y))\;\leq\; \theta\, d_X(x,y)\qquad\mbox{for all $x,y\in X$.}
\]
\end{defn}
Banach's Contraction
Theorem (see, e.g., \cite[Appendix~A]{Sch})
is a paradigmatic fixed point
theorem for contractions. We recall
it as a model for the slight generalizations
which we actually need for our purposes:
%
%
\begin{la}\label{banachfix}
Let $(X,d)$ be a $($non-empty$)$
complete metric
space and $f\colon X\to X$ be a contraction,
with contraction constant
$\theta\in [0,1[$. Thus
\[
d(f(x),f(y))\;<\; \theta\, d(x,y)\quad\mbox{for all $x,y\in X$.}
\]
Then $f(p)=p$ for a unique point $p\in X$.
Given any $x_0\in X$, we have $\lim_{n\to\infty}f^n(x_0)=p$.
Furthermore, the a priori estimate
\[
d(f^n(x_0),p)\leq \frac{\theta^n}{1-\theta}\, d(f(x_0),x_0)
\]
holds, for each $n\in \N_0$.\Punkt
\end{la}
Unfortunately, we are not always
in the situation of this theorem. But
the simple variants compiled
in the next proposition
are flexible enough for our purposes.
%
%
\begin{prop}\label{banfix2}
Let $(X,d)$ be a metric
space, $U\sub X$ be a subset and $f\colon U\to X$
be a contraction, with contraction constant~$\theta$.
Then the following holds:
\begin{itemize}
\item[\rm (a)]
$f$ has at most one fixed point.
\item[\rm (b)]
If $x_0\in U$ is a point and $n\in \N_0$
such that
$f^{n+1}(x_0)$ is defined, then
%
\begin{equation}\label{iterctr}
d(f^{k+1}(x_0), f^k(x_0))\;\leq\; \theta^k\, d(f(x_0),x_0)
\end{equation}
for all $k\in \{0,\ldots, n\}$,
and $\,d(f^{n+1}(x_0), x_0)\leq \frac{1-\theta^{n+1}}{1-\theta}\,
d(f(x_0),x_0)$.
\item[\rm (c)]
If $x_0\in U$ is a point
such that $f^n(x_0)$ is defined
for all $n\in \N$, then $(f^n(x_0))_{n\in \N}$
is a Cauchy sequence in~$U$,
and
%
\begin{equation}\label{makeCauex}
d(f^{n+k}(x_0), f^n(x_0))\;\leq\; \frac{\theta^n(1-\theta^k)}{1-\theta}\,
d(f(x_0),x_0)\quad\mbox{for all $\, n,k\in \N_0$.}
\end{equation}
If $(f^n(x_0))_{n\in \N}$ converges
to some $x\in U$, then $x$ is a fixed point
of~$f$, and
%
\begin{equation}\label{aprio2}
d(x, f^n(x_0))\;\leq\; \frac{\theta^n}{1-\theta}\,
d(f(x_0),x)\quad\mbox{for all $\,n \in \N_0$.}
\end{equation}
If $f^n(x_0)$ is defined for all $n\in \N$
and $f$ has a fixed point~$x$, then
$f^n(x)\to x$ as $n\to\infty$.
\item[\rm (d)]
Assume that $U=\wb{B}_r(x_0)$
is a closed ball of radius~$r$
around a point $x_0\in X$,
and
$d(f(x_0),x_0)\leq (1-\theta)r$.
Then $f^n(x_0)$ is defined for all
$n\in \N_0$.
Hence $f$ has a fixed point
inside~$\wb{B}_r(x_0)$, provided $X$
is complete.
Likewise, $f$ has a fixed point in the open
ball $B_r(x_0)$
if $X$ is complete,
$U=B_r(x_0)$,
and $d(f(x_0),x_0)<(1-\theta)r$.
\end{itemize}
\end{prop}
\begin{proof}
(a) If $x, y\in U$ are fixed points of~$f$,
then $d(x,y)=d(f(x),f(y))\leq \theta d(x,y)$,
entailing that $d(x,y)=0$ and thus $x=y$.

(b) For $k=0$, the formula (\ref{iterctr})
is trivial.
If $k<n$ and
$d(f^{k+1}(x_0),f^k(x_0))\leq \theta^k\, d(f(x_0),x_0)$,
then
$d(f^{k+2}(x_0),f^{k+1}(x_0))
=d(f(f^{k+1}(x_0)), f(f^k(x_0)))\leq $\linebreak
$\theta \, d(f^{k+1}(x_0), f^k(x_0))
\leq \theta^{k+1}\, d(f(x_0),x_0)$.
Thus (\ref{iterctr}) holds in general.

Using the triangle inequality and
the summation formula for the geometric series,
we obtain the estimates $d(f^{n+1}(x_0),x_0)
\leq\sum_{k=0}^n d(f^{k+1}(x_0),f^k(x_0))
\leq \sum_{k=0}^n \theta^k\, d(f(x_0),x_0)
=\frac{1-\theta^{n+1}}{1-\theta}\, d(f(x_0),x_0)$,
as asserted.

(c) Using both of the estimates from (b), obtain
\[
d(f^{n+k}(x_0),f^n(x_0)) \;\leq\;
\frac{1-\theta^k}{1-\theta}\, d(f^{n+1}(x_0),f^n(x_0))
\;\leq\;
\frac{1-\theta^k}{1-\theta}\, \theta^n d(f(x_0),x_0)\,.
\]
Thus (\ref{makeCauex}) holds,
and thus $(f^n(x_0))_{n\in \N}$ is a Cauchy sequence.
If $f^n(x_0)\to x$ for some $x\in U$,
then $x=\lim_{n\to\infty}f^{n+1}(x_0)=f(\lim_{n\to\infty}f^n(x_0))=
f(x)$ by continuity of~$f$, whence indeed~$x$
is a fixed point of~$f$.
Letting now $k\to\infty$ in (\ref{makeCauex}),
we obtain (\ref{aprio2}).

To prove the final assertion, assume that $f$ has a fixed
point $x$ and that $f^n(x_0)$
is defined for all~$n$.
We choose a completion $\wb{X}$
of~$X$ (with $X\sub \wb{X}$)
and let $\wb{U}$ be the closure of~$U$ in~$\wb{X}$.
Then $f$ extends to a contraction
$\wb{U}\to \wb{X}$, which we also denote
by~$f$. Since $(f^n(x_0))_{n\in \N}$
is a Cauchy sequence in $\wb{U}$ and $\wb{U}$
is complete, we deduce that $f^n(x_0)\to y$
for some $y\in \wb{U}$.
Then both $y$ and $x$ are fixed points of~$f$
and hence $x=y$.

(d) We show by induction that $f^n(x_0)$ is defined
for all $n\in \N$.
For $n=1$, this is trivial.
If $f^n(x_0)$ is defined,
then
\[
d(f^n(x_0),x_0)\, \leq\,
\frac{1-\theta^n}{1-\theta}\,d(f(x_0),x_0)
\,\leq\,
\frac{1-\theta^n}{1-\theta}\,(1-\theta)r\, \leq \, r
\]
and thus $f^n(x_0)\in \wb{B}_r(x_0)$,
whence also $f^{n+1}(x_0)=f(f^n(x_0))$ is defined.
Then $f^n(x_0)\in \wb{B}_r(x_0)$ for each $n\in \N$.
By (c), $(f^n(x_0))_{n\in \N}$ is a Cauchy
sequence. If $X$ is complete, then
so is $\wb{B}_r(x_0)$ and
thus $(f^n(x_0))_{n\in \N}$
converges to some point $x\in \wb{B}_r(x_0)$,
which is a fixed point of $f$ by~(c).
Finally, if $U=B_r(x_0)$ and $d(f(x_0),x_0)<(1-\theta)r$,
there exists $s\in \;]0,r[$ such that
$d(f(x_0),x_0)\leq (1-\theta)s$.
By the preceding, $f^n(x_0)\in \wb{B}_s(x_0)\sub B_r(x_0)$
for all $n\in \N$,
and $f$ has a fixed
point in $\wb{B}_s(x_0)\sub B_r(x_0)$.
\end{proof}
We are interested in uniform
families of contractions.
%
%
\begin{defn}\label{defunicon}
Let $F$ be a Banach space over
a valued field~$\K$,
and $U\sub F$ be a subset.
A family $(f_p)_{p\in P}$ of mappings
$f_p\colon U\to F$
is called a
\emph{uniform family of contractions}
if there exists $\theta\in [0,1[$
(a ``uniform contraction constant'')
such that
\[
\|f_p(x)-f_p(y)\|\,\leq\, \theta\|x-y\|\quad
\mbox{for all $x,y\in U$ and $p\in P$.}
\]
\end{defn}
If~$U$ is closed
and each $f_p$ is a self-map
of~$U$ here, then
Banach's Contraction Theorem
ensures
that, for each $p\in P$, the map
$f_p$ has a unique
fixed point $x_p$.
Our goal is to understand
the dependence of $x_p$ on the parameter~$p$.
In particular, for $P$ a subset
of a topological $\K$-vector space,
we want to find conditions
ensuring that the map $P\to F$,
$p\mto x_p$ is continuously differentiable.
We discuss dependence of fixed points
on parameters in two steps.
%
\begin{prop}\label{preparadep}
Let
$P$ be a topological space
and
$F$ be a Banach space over
a valued field~$\K$.
Let $U\sub F$ be a subset with dense interior and
$f\colon P\times U\to F$ be a map
such that $(f_p)_{p\in P}$
is a uniform family of contractions,
where $f_p:=f(p,\sbull)\colon U\to F$.
We assume that $f_p$ has a fixed point
$x_p$, for each $p\in P$.
Furthermore,
we assume that $U$ is open or
$f(P\times U)\sub U$
$($whence every $f_p$ is a self-map
of $U)$.
Then the following
holds:
\begin{itemize}
\item[\rm (a)]
If $f$ is continuous, then also the map
$\phi\colon P\to F$, $\phi(p):=x_p$
is continuous.
\item[\rm (b)]
If $P$ is a subset of a topological
$\K$-vector space~$E$
and $f$ is Lipschitz continuous,
then also $\phi$ is Lipschitz continuous.
\item[\rm (c)]
If $P$ is a subset with dense interior
of a topological
$\K$-vector space~$E$
and $f$ is $SC^1$,
then also $\phi$ is $SC^1$.
\item[\rm (d)]
If $P$ is a subset with dense interior
of a topological
$\K$-vector space~$E$
and $f$ is $C^1$,
then also $\phi$ is $C^1$.
\end{itemize}
\end{prop}
\begin{rem}
We shall see that
the differential of $\phi$ at~$p\in P$ is given by
%
\begin{equation}\label{difffp}
\phi'(p)\;=\; (\id_F -\beta_2)^{-1}\circ \beta_1
\end{equation}
in the situation of
Proposition~\ref{preparadep}\,(c) and (d),
where
$\beta_1:=d_1f(p,x_p,\sbull):=df(p,x_p,\sbull, 0)\in \cL(E,F)$
and $\beta_2:=d_2f(p,x_p,\sbull):=df(p,x_p,0,\sbull)\in \cL(F)$.
\end{rem}
\begin{proof}
Let $\theta\in [0,1[$ be a uniform contraction constant
for $(f_p)_{p\in P}$.

(a) If $f$ is continuous,
$p\in P$ and $\ve>0$,
we find a neighborhood $Q\sub P$ of~$p$
such that $\|x_p-f_q(x_p)\|\leq (1-\theta)\ve$
for all $q\in Q$.
If $f(P\times U)\sub U$,
then $\|x_p-x_q\|\leq \frac{1}{1-\theta}\|x_p-f_q(x_p)\|
\leq \ve$, by (\ref{aprio2})
in Proposition~\ref{banfix2}\,(c).
If $U$ is open, we may
assume that $\wb{B}_\ve(x_p)\sub U$
after shrinking~$\ve$.
Then Proposition~\ref{banfix2}\,(d)
applies to $f_q$
as a map $\wb{B}_\ve(x_p)\to E$
for each $q\in Q$,
showing that $f_q^n(x_p)$ is defined
for each $n\in \N$ and
$x_q=\lim_{n\to\infty}f_q^n(x_p)\in \wb{B}_\ve(x_p)$,
that is, $\|x_p-x_q\|\leq\ve$.

(b) Given $s\in P$, there exists a gauge $\gamma$ on~$E$ and
$R,L>0$ such that
%
\begin{equation}\label{bore1}
\|f_p(x)-f_q(y)\|\;\leq\; L\max\{\|p-q\|_\gamma,\|x-y\|\}
\end{equation}
for all $p,q\in P\cap B_R^\gamma(s)$ and $x,y\in U\cap
B_{2R}^F(x_s)$.
Let $\eta$ be a gauge on $E$ such that $\|u+v\|\gamma
\leq \max\{\|u\|_\eta,\|v\|_\eta\}$
for all $u,v\in E$.
If $U$ is open, we assume that
$\wb{B}_{2R}(x_s)\sub U$.
Set $r:=\min\{(1-\theta)R,(1-\theta)R/L\}$.
Given $p\in P\cap B_r^\eta(s)$, we have
$\|f_p(x_s)-x_s\|=
\|f_p(x_s)-f_s(x_s)\|
\leq L\|p-s\|_\gamma \leq (1-\theta)R$
by (\ref{bore1})
and hence $x_p\in \wb{B}_R(x_s)$,
by the proof of~(a).
Given $p,q\in P\cap B_r^\eta(s)$,
we either have $\|p-q\|_\gamma=0$;
then $\|f_q(x_p)-x_p\|=\|f_q(x_p)-f_p(x_p)\|=0$
by (\ref{bore1}), whence
$f_q(x_p)=x_p$. Hence $x_q=x_p$,
whence $\|x_q-x_p\|=0$ and thus
%
\begin{equation}\label{seelater}
\|x_p-x_q\|\;\leq\; \frac{L}{1-\theta}\, \|p-q\|_\gamma
\end{equation}
in particular.
Otherwise, $0<\ve:= \|p-q\|_\gamma \leq\max\{\|p-s\|_\eta,
\|q-s\|_\eta\}< r$
and
$\|f_q(x_p)-x_p\|=\|f_q(x_p)-f_p(x_p)\|\leq
L\|p-q\|_\gamma=\ve L$, by (\ref{bore1}).
If $U$ is open, then $\ve L <\ve L/(1-\theta)
\leq r L/(1-\theta)\leq R$
and thus $\wb{B}_{\ve L/(1-\theta)}(x_p)\sub U$.
Hence $\|x_p-x_q\|\leq \ve L/(1-\theta)=\|p-q\|_\gamma L/(1-\theta)$,
by the proof of\,(a).
Thus (\ref{seelater}) holds for all $p,q\in P\cap B^\eta_r(x_s)$.
We deduce that $\phi$ is Lipschitz continuous.

(c)
Assume that $f$ is strictly differentiable.
Given $s\in P$, by strict differentiability
of~$f$ at $(s,x_s)$, there exists a gauge $\gamma$ on~$E$
such that, for each $\ve>0$,
there exists $\delta>0$ such that
%
\begin{equation}\label{epdelt}
\|f(p,x)-f(q,y)-f'(s,x_s).(p-q,x-y)\|\;\leq\; \ve\,\max\{\|p-q\|_\gamma,
\|x-y\|\}
\end{equation}
for all $p,q\in P\cap B^\gamma_\delta(s)$ and
$x,y\in U\cap B_\delta^F(x_s)$.
Since $f$ is strictly differentiable, it is Lipschitz
continuous, whence also $\phi$ is Lipschitz continuous,
by (b). Hence, after replacing $\gamma$ by a larger
gauge, we may assume that 
%
\begin{equation}\label{hencesho}
\|x_p-x_q\|\;\leq\;  \|p-q\|_\gamma \quad\mbox{for all
$p,q\in P\cap B_1^\gamma (x_s)$.}
\end{equation}
Given $\ve>0$,
choose $\delta\in \;]0,1]$
such that (\ref{epdelt}) holds.
Taking $q=s$ in (\ref{hencesho}),
it follows that
$x_p\in U\cap B^F_\delta(x_s)$
for each $p\in B_\delta^\gamma(s)$.

It is useful to write
$f'(s,x_s)(u,v)=\beta_1(u)+\beta_2(v)$
in terms of the partial differentials
$\beta_1:=d_1f(x,x_s,\sbull)\colon E\to F$
and $\beta_2:=d_2f(x,x_s,\sbull)\colon F\to F$.
Then $\|\beta_2\|\leq\theta<1$
by Lemma~\ref{sweela}.
Abbreviate
\[
R(p,q):=f(p,x_p)-f(q,x_q)-f'(s,x_s).(p-q,x_p-x_q)
\]
for $p,q\in P\cap B_\delta^\gamma(s)$.
Combining (\ref{epdelt}) and (\ref{hencesho}), we find that
%
\begin{equation}\label{call44}
\|R(p,q)\|\, \leq \, \ve \max\{\|p-q\|_\gamma,\|x_p-x_q\|\}
\, \leq \, \ve \|p-q\|_\gamma
\end{equation}
for all $p,q\in P\cap B_\delta^\gamma(s)$.
Now
\begin{eqnarray*}
x_p-x_q &=& f(p,x_p)-f(q,x_q)
\,=\, f'(s,x_s).(p-q,x_p-x_q)+R(p,q)\\
&=& \beta_1(p-q)+\beta_2(x_p-x_q)+R(p,q)
\end{eqnarray*}
and therefore $x_p-x_q-(\id_F -\beta_2)^{-1}.\beta_1(p-q)=
(\id_F-\beta_2)^{-1}R(p,q)$,
where\linebreak
$\|(\id_F-\beta_2)^{-1}R(p,q)\|\leq \|(\id_F-\beta_2)^{-1}\|\,
\|R(p,q)\|
\leq \ve (1-\theta)^{-1}\|p-q\|_\gamma$,
by (\ref{geomse}) in Proposition~\ref{glislie}
and (\ref{call44}).
We have shown that $\phi$ is strictly
differentiable at~$s$,
with the desired differential.

(d) Being $C^1$, $f$ is Lipschitz continuous,
whence $\phi$ is Lipschitz continuous, by~(b).
Thus $\phi^{]1[}$ is continuous.
To see that $\phi$ is~$C^1$,
it only remains to show that,
for all $p_0\in P$ and $q_0\in E$,
there exists an open neighborhood
$W\sub P^{[1]}$ of $(p_0,q_0,0)$
and a continuous map $g\colon W\to F$
which extends the difference quotient map
$\phi^{]1[}|_{W\cap P^{]1[}}\colon W\cap P^{]1[}\to F$.
Then $\phi^{]1[}$ has a continuous
extension $\phi^{[1]}$ to all
of $P^{[1]}$, by Lemma~\ref{Dieudonn},
and thus $\phi$ will be $C^1$.
Our strategy is the following: We write
\begin{equation}
(f^{n+1}_{p+tq}(x_p)-x_p)/t\;=\;
\sum_{k=0}^n(f_{p+tq}^{k+1}(x_p)-f^k_{p+tq}(x_p))/t
\end{equation}
for
$(p,q,t)$ in a suitable neighborhood
$W$ of $(p_0,q_0,0)$.
For $W$ sufficiently small,
the left hand side converges to $\frac{x_{p+tq}-x_p}{t}
=\frac{\phi(p+tq)-\phi(p)}{t}$
as $n\to\infty$.
Furthermore, we can achieve that
each term on the right hand side extends continuously
to all of $W$, and
that the series
converges uniformly to a continuous
function on~$W$.
This will be our desired continuous extension~$g$.

Let us carry this out in detail.
Case~1. If $f(P\times U)\sub U$,
we set $W_0:=P^{[1]}$.
Case~2. Otherwise, $U$ is open,
whence there exists $\ve>0$ such that
$\wb{B}_{2\ve}^F(x_{p_0})\sub U$.
Since $f$ and $\phi$ are continuous
and $f_{p_0}(x_{p_0})=x_{p_0}$,
we find an open neighborhood $Q\sub P$ of $p_0$
such that $\|x_p-f_q(x_p)\|\leq (1-\theta)\ve$
for all $p,q\in Q$.
Then $f_q^k(x_p)$ is defined for all $k\in \N_0$,
$f_q^k(x_p)\in B_\ve^F(x_p)$,
and $x_q\in \wb{B}_\ve^F(x_p)$
(cf.\ proof of\,(a)).
In particular, $x_p\in \wb{B}_\ve^F(x_{p_0})$.
We now set $W_0:=Q^{[1]}$
and note that, if $(p,q,t)\in Q^{[1]}$,
then $p,p+tq\in Q$,
whence $f_{p+tq}^k(x_p)$
is defined
for all $k\in \N_0$
and $\lim_{k\to\infty}f_{p+tq}^k(x_p)=x_{p+tq}$
(by the preceding considerations).

In either case, we define
\[
h_0\colon W_0 \to F\,,\qquad
h_0(p,q,t)\,=\, f^{[1]}(p,x_p,q,0,t)
\]
and note that $h_0$ is a continuous map such that
%
%
\begin{equation}\label{dexter}
h_0(p,q,t)=(f_{p+tq}(x_p)-f_p(x_p))/t
=(f_{p+tq}(x_p)-x_p)/t\quad\mbox{if $t\not=0$.}
\end{equation}
For all $k\in \N$ and $(p,q,t)\in W_0$,
we have
%
\begin{eqnarray}
\lefteqn{\frac{f^{k+1}_{p+tq}(x_p)-f^k_{p+tq}(x_p)}{t}}\qquad\notag \\
&=&
\frac{f \bigl( p+tq, f^{k-1}_{p+tq} (x_p)+t
\, \frac{f^k_{p+tq}(x_p)-f^{k-1}_{p+tq}(x_p)}{t} \bigr)
- f(p+tq, f^{k-1}_{p+tq} (x_p))}{t} \notag \\
&=& f^{[1]} \Big(
p+tq, f^{k-1}_{p+tq} (x_p), 0, \frac{f^k_{p+tq}(x_p)-f^{k-1}_{p+tq} (x_p)}{t},
t \Big)\,. \label{dext2}
\end{eqnarray}
Recursively, we define
\[
h_k \colon W_0 \to F\,,\quad
h_k (p,q,t)\, :=\,
f^{[1]}\bigl(p+tq,f^{k-1}_{p+tq}(x_p),0, h_{k-1}(p,q,t),t\bigr)
\]
for $k\in \N$.
A simple induction based on (\ref{dexter})
and (\ref{dext2}) shows that the definition of
$h_k$ makes sense for each $k\in \N_0$,
and that
%
\begin{equation}\label{intph}
h_k(p,q,t)\,=\, \frac{f^{k+1}_{p+tq}(x_p)-f^k_{p+tq}(x_p)}{t}
\quad \mbox{for all $(p,q,t)\in W_0$ with $t\not=0$.}
\end{equation}
The function $h_0\colon W_0 \to F$,
$(p,q,t)\mto f^{[1]}(p,x_p,q,0,t)$
being continuous,
we find an open neighborhood
$W \sub W_0$ of $(p_0,q_0,0)$
and $C\in [0,\infty[$
such that
\[
\|f^{[1]}(p,x_p,q,0,t)\|\;\leq\; C\qquad \mbox{for all $(p,q,t)\in W$.}
\]
For all $(p,q,t)\in W$ such that $t\not=0$,
we have
\begin{eqnarray*}
\|f_{p+tq}^{k+1}(x_p)-f_{p+tq}^k(x_p)\|
&=& \|f_{p+tq}^k(f_{p+tq}(x_p))-f_{p+tq}^k(x_p)\|\\
&\leq & \theta^k\|f_{p+tq}(x_p)-x_p\|\\
& = &|t|\theta^k\left\|(f_{p+tq}(x_p)-f_p(x_p))/t\right\|\\
&= & |t|\theta^k \|f^{[1]}(p,x_p,q,0,t)\|
\; \leq \; |t|\theta^k C\,.
\end{eqnarray*}
Combining this with (\ref{intph}), for each $k\in \N_0$
we see that
$\|h_k(p,q,t)\|\leq \theta^k C$
for all $(p,q,t)\in W$ such that $t\not=0$,
and thus
%
\begin{equation}\label{estih}
\|h_k(p,q,t)\|\;\leq \; \theta^k C
\quad\mbox{for all $(p,q,t)\in W$,}
\end{equation}
because $h_k$ is continuous and $W\cap P^{]1[}$
is dense in~$W$.
As a consequence, $\sum_{k=0}^\infty \|h_k|_W\|_\infty\leq
\sum_{k=0}^\infty\theta^k C=\frac{1}{1-\theta}C<\infty$,
whence the series
$\sum_{k=0}^\infty h_k|_W$ of bounded continuous
functions converges
uniformly and absolutely. Thus
\[
g(p,q,t)\; :=\; \sum_{k=0}^\infty \, h_k (p,q,t)
\]
exists for all $(p,q,t)\in W$,
and $g\colon W\to F$ is continuous.
It only remains to observe that
\[
\frac{f_{p+tq}^{n+1}(x_p)-x_p}{t}
\; =\;  \sum_{k=0}^n \, \frac{f_{p+tq}^{k+1}(x_p)
-f_{p+tq}^k(x_p)}{t}
\; =\;
 \sum_{k=0}^n \, h_k(p,q,t)
\]
for all $(p,q,t)\in W$ such that $t\not=0$.
Since the left hand side converges to $\frac{x_{p+tq}-x_p}{t}$
and the right hand side converges to $g(p,q,t)$,
we obtain
\[
\frac{\phi(p+tq)-\phi(p)}{t}\;=\;\sum_{k=0}^\infty \, h_k(p,q,t)
\;=\; g(p,q,t)\,.
\]
Thus $g\colon W\to F$ is a continuous map
which extends $\phi^{]1[}|_{W\cap P^{]1[}}$,
as desired. This completes the proof.
\end{proof}
We now state and prove the main result of
this section.
%
%
%
\begin{thm}[Dependence of Fixed Points of Parameters]\label{fixpardep}
Let $\K$ be a\linebreak
valued field,
$E$ be a topological $\K$-vector space, and
$F$ be a Banach space over~$\K$.
Let $P\sub E$ be a subset with dense
interior, $U\sub F$ be open,
and
$f\colon P \times U\to F$
be a continuous map such that
$f_p:=f(p,\sbull)\colon U\to F$
defines a uniform family
$(f_p)_{p\in P}$
of contractions.
Then the following holds:
\begin{itemize}
\item[\rm (a)]
The set $Q$ of all $p\in P$ such that
$f_p$ has a fixed point $x_p$
is open in~$P$.
\item[\rm (b)]
If $f$ is $C^k$
for some $k\in \N_0\cup\{\infty\}$
$($resp., $SC^k$, resp., $LC^k)$,
then also
$\phi\colon Q \to U$, $\phi(p):=x_p$
is $C^k$, $SC^k$, resp., $LC^k$.
\end{itemize}
\end{thm}
\begin{proof}
Let $\theta\in [0,1[$ be a uniform contraction constant
for $(f_p)_{p\in P}$.

(a) If $p\in Q$, there is $r>0$
such that $\wb{B}_r(x_p)\sub U$.
There is a neighborhood
$S\sub Q$ of~$p$ such that
$\|f_q(x_p)-x_p\|=\|f_q(x_p)-f_p(x_p)\|\leq (1-\theta)r$
for all $q\in S$.
Now Proposition~\ref{banfix2}\,(d)
shows that $f_q$ has a fixed point $x_q$ in
$\wb{B}_r(x_p)$,
for each $q\in S$. Thus $S\sub Q$ and we deduce that
$Q$ is open.

(b) We may assume that $k<\infty$.
The proof is by induction on $k\in \N_0$.
The case $k=0$ is covered by
Proposition~\ref{preparadep}
(a) and (b).
Now assume that our assertion holds
for some $k\in \N_0$
and assume that $f$ is $C^{k+1}$
(resp., $SC^{k+1}$, resp., $LC^{k+1}$).
Then $\phi$ is $C^k$ (resp., $SC^k$, resp., $LC^k$) by the
induction hypothesis.
Furthermore, $\phi$ is $C^1$
by Proposition~\ref{preparadep}\,(d),
and if $f$ is $SC^1$, then $\phi$ is $SC^1$,
by Proposition~\ref{preparadep}\,(c).
If $k=0$, this already completes the
induction step for continuously differentiable
and for strictly
differentiable maps.
If $f$ is $C^{k+1}$ or $SC^{k+1}$
with $k\geq 1$ (whence $f$ is
$LC^1$ in particular)
or if $k\in \N_0$ and
$f$ is $LC^{k+1}$, it
remains to show
that $\phi^{[1]}$ is $C^k$,
$SC^k$ and
$LC^k$, respectively.
It suffices to prove that
$\phi^{[1]}|_Z$
is $C^k$ (resp., $SC^k$, resp., $LC^k$)
for $Z$ ranging through an
open cover of~$Q^{[1]}$
(cf.\ \cite[Lemma~4.9]{Ber}).
Because
%
\begin{equation}\label{onWop}
\phi^{[1]}(p,q,t)=\frac{\phi(p+tq)
-\phi(p)}{t}
\quad\mbox{for $(p,q,t)\in Q^{]1[}$,}
\end{equation}
we observe first that $\phi^{[1]}$ is $C^k$
(resp., $SC^k$, resp., $LC^k$)
on the open subset
$Q^{]1[}$ of $Q^{[1]}$.
Given $(p_0,q_0)\in Q\times E$,
let us find a description of
$\phi^{[1]}$ on a neighborhood of
$(p_0,q_0,0)$ in~$Q^{[1]}$.
Given $(p,q,t)\in Q^{]1[}$,
we have
\begin{eqnarray*}
\phi^{[1]}(p,q,t) &=& \frac{\phi(p+tq)-\phi(p)}{t}
\;=\;  \frac{f(p+tq ,\phi(p+tq))-f(p,\phi(p))}{t}\\
&=& \frac{f(p+tq ,\phi(p)+t\frac{\phi(p+tq)-\phi(p)}{t})
-f(p,\phi(p))}{t}
\end{eqnarray*}
and thus
%
\begin{equation}\label{aggfix}
\phi^{[1]}(p,q,t) \;=\;
f^{[1]}(p,\phi(p),q,\phi^{[1]}(p,q,t),t)
\end{equation}
for all $(p,q,t)\in Q^{]1[}$.
As both the left and right hand side
of (\ref{aggfix}) make sense
for all $(p,q,t)\in Q^{[1]}$
and are continuous there,
they coincide
on all of $Q^{[1]}$.\\[3mm]
Observe
that (\ref{aggfix}) means
that, for each parameter $(p,q,t)\in Q^{[1]}$,
the element
$\phi^{[1]}(p,q,t)\in F$
is a fixed
point of the function $f^{[1]}(p,\phi(p),q,\sbull, t)$.
We make this more precise
now and show that
we are dealing
with a $C^k$- (resp., $SC^k$-, resp., $LC^k$-)
family of uniform
contractions, which will enable us
to deduce from the induction hypothesis
that $\phi^{[1]}$ is $C^k$ (resp.,
$SC^k$, resp., $LC^k$).\\[3mm]
Abbreviate $v_0:=\phi(p_0)$
and $w_0:=\phi^{[1]}(p_0,q_0,0)$.
Because $U^{[1]}$ is open in $F\times F\times \K$
and $(v_0,w_0,0)\in U^{[1]}$,
there exist open neighborhoods
$V\sub U$ of $v_0$,
$W\sub F $ of $w_0$
and a $0$-neighborhood
$S\sub \K$ such that $V\times W\times S\sub U^{[1]}$
and thus $V+SW\sub U$.
Then
$P_1 :=
(Q^{[1]}\cap (Q\times E\times S))
\times V$
is a subset of $E \times E\times \K\times F$
with dense interior,
$W$ is an open subset
of the Banach space $F$, and
$(p,v,q,w,t)\in (Q\times U)^{[1]}$
for all $(p,q,t,v,w)\in P_1\times W$.
We can therefore define a
$C^k$-map (resp., $SC^k$-map, resp.,
$LC^k$-map)
%
\begin{equation}\label{toolg}
g\colon P_1\times W \to F\,,\quad
g(p,q,t,v,w):= f^{[1]}(p,v,q,w,t)\,.
\end{equation}
Choose $\ve\in \;]0,1-\theta[$.
Because $f$ is $LC^1$,
Lemma~\ref{monstrtech}
entails that
there
exists $s>0$ such that
$B_s^F(w_0)\sub W$
and an open neighborhood
$P_2\sub P_1$
of $(p_0,q_0,0,v_0)$
such that
%
\begin{equation}\label{seeq}
\|f^{[1]}(p,v,q,w_1,t)-f^{[1]}(p,v,q,w_2,t)
-f'_{p_0}(v_0).(w_1-w_2)\|\,\leq\,
\ve \|w_1-w_2\|
\end{equation}
whenever
$(p,q,t,v,w_1),
(p,q,t,v,w_2)\in P_2 \times B_s^F(w_0)$.
Since
$\|f'_{p_0}(v_0)\|\leq \theta$
by Lemma~\ref{sweela},
(\ref{seeq}) entails that
%
\begin{eqnarray}
\|f^{[1]}(p,v,q,w_1,t)-f^{[1]}(p,v,q,w_2,t)\|
& \leq &
(\|f'_{p_0}(v_0)\| +\ve)\cdot \|w_1-w_2\|\notag\\
&\leq  &(\theta+\ve)\cdot \|w_1-w_2\| \label{sefq}
\end{eqnarray}
whenever
$(p,q,t,v,w_1),
(p,q,t,v,w_2)\in P_2\times B_s^F(w_0)$.
Since $\theta+\ve<1$,
we see that the restriction of~$g$
to a map $g\colon P_2\times B_s^F(w_0)\to F$
is a uniform family of contractions.
The map $g$ being $C^k$ (resp., $SC^k$, resp., $LC^k$),
we know by induction that
the set $Q_2\sub P_2$
of all $(p,q,t,v)\in P_2$ such that $g(p,q,t,v,\sbull)\colon
B_s^F(w_0)\to F$ has a fixed point $\psi(p,q,t,v)$
is open in $P_2$,
and that the map $\psi\colon Q_2\to F$
so obtained is $C^k$ (resp., $SC^k$, resp., $LC^k$).
By continuity of $\phi$ and $\phi^{[1]}$,
there exists an open neighborhood
$Z$ of $(p_0,q_0,0)$ in $P^{[1]}$
such that $(p,,q,t,\phi(t))\in Q_2$ for all
$(p,q,t)\in Z$, and $\phi^{[1]}(p,q,t)\in B_s^F(w_0)$.
Then
(\ref{aggfix}) entails that
\[
\phi^{[1]}(p,q,t)=\psi(p,q,t,\phi(p))\qquad\mbox{for all
$(p,q,t)\in Z$,}
\]
whence $\phi^{[1]}|_Z$ is $C^k$ (resp., $SC^k$,
resp., $LC^k$).
This completes the proof.
\end{proof}
%
%
%
%
%
%
%
%
%
%
%
%
%
%
\section{Inverse and implicit function theorems}\label{sectimpl}
In this section, we prove
inverse- and implicit function
theorems for
mappings into Banach spaces over
valued fields.
In particular, we obtain the
following analog of
the classical Inverse Function Theorem:
%
%
%
\begin{thm}[Inverse Function Theorem]\label{basicinv}
Let $E$ be a Banach space
over a valued field~$\K$
and $f\colon U\to E$ be a mapping
on an open subset $U\sub E$
which is $C^k$, $SC^k$ and $LC^k$, respectively,
where $k\in \N \cup\{\infty\}$.
In the $C^k$-case, we
require that $k\geq 2$ or that
$\K$ is locally compact and $E$ has finite dimension.
If $f'(x)\in \GL(E)$ for some $x\in U$,
then there exists an open neighborhood
$V\sub U$ of~$x$ such that $f(V)$ is open
in~$E$ and $f|_V\colon V\to f(V)$
is a $C^k$-diffeomorphism
$($resp., an $SC^k$-diffeomorphism,
resp., an  $LC^k$-diffeomorphism$)$.
\end{thm}
While the classical implicit
function theorem deals with
functions between real Banach spaces,
we can discuss
implicit functions from arbitrary topological
vector spaces to Banach spaces:
%
%
\begin{thm}[Generalized Implicit Function Theorem]\label{genimp}
Let $\K$ be a valued field,
$E$ be a topological $\K$-vector space,
$F$ be a Banach space over~$\K$,
and $f\colon U\times V\to F$ be a
mapping, where $U\sub E$ is a subset with dense
interior and $V\sub F$ is open.
We assume that $f$ is $C^k$,
$SC^k$, respectively, $LC^k$
for some $k\in \N\cup\{\infty\}$.
In the $C^k$-case, we require 
that $k\geq 2$ or that $\K$ is locally compact and $F$ has finite
dimension.
Given $x\in U$, abbreviate $f_x:=f(x,\sbull)
\colon V\to F$.
If $f(x_0,y_0)=0$ for some $(x_0,y_0)\in U\times V$
and $f_{x_0}'(y_0)\in \GL(F)$,
then there exist open neighborhoods
$U_0\sub U$ of~$x_0$ and $V_0\sub V$ of~$y_0$
such that
\[
\{(x,y)\in U_0\times V_0\colon f(x,y)=0\}
\; =\; \graph\lambda
\]
for a map $\lambda \colon U_0\to V_0$
which is $C^k$, $SC^k$, resp., $LC^k$.
\end{thm}
We shall deduce the
Generalized Implicit
Function Theorem
from an ``Inverse Function Theorem
with Parameters.''
It is useful to get an idea
of the main steps of the proof
(for $SC^k$-maps, say)
before we carry them out in detail.\\[3mm]
\textbf{Strategy of the proof.}
Since
$f_{x_0}'(y_0) \in \GL(F)$,
where $\GL(F)$ is open in $\cL(F)$ and
the map $U \to \cL(F)$, $x \mto f'_x(y_0)$
is continuous, we see that
$f_x'(y_0)\in \GL(F)$ for $x$ close to $x_0$.
Then $f_x^{-1}$ exists locally
around $f_x(y_0)$,
by the Inverse Function Theorem (Theorem~\ref{basicinv}).
Now the essential point is that
the map
%
\begin{equation}\label{needdomain}
(x,z)\; \mto \; (f_x)^{-1}(z)
\end{equation}
actually makes sense on a whole
neighborhood $U_0\times W $ of $(x_0,0)$ in
$U\times F$,
and is $SC^k$ there.
This assertion is the main content
of the ``Inverse Function Theorem
with Parameters''
just announced.
Once we have this, the rest is easy:
The~map
\[
\lambda \colon U_0\to F \, ,\quad \lambda (x)\;:=\; (f_x)^{-1}(0)
\]
is $SC^k$,
and $f(x,\lambda(x))=0$.\vspace{3mm}\Punkt

\noindent
Motivated by these considerations,
we now proceed
as follows:
First, we formulate and
prove versions of the inverse
function theorem which provide
quantitative information on the ``size'' of the
domain of the local inverse $f^{-1}$,
and the size of the images of balls
(under $f$).
These studies can be carried out
under quite weak hypotheses
(for Lipschitz maps).
They will enable us to see that
the domain of definition in (\ref{needdomain})
is a neighborhood.
The second step,
then, is the exact formulation
and proof of the Inverse Function Theorem
with Parameters.

%
%
%
%
%
%
%
%
%
%
%
%
%
%
\subsection*{Lipschitz Inverse Function Theorems}\label{snultra}
We now prove a Lipschitz Inverse
Function Theorem for self-maps
of general Banach spaces,
and a variant
for ultrametric Banach spaces.
As a rule, inverse function theorems
are available in much stronger form
in the ultrametric case.
This is a general phenomenon
which we shall encounter repeatedly.
%
%
\begin{thm}[Lipschitz Inverse Function Theorem]\label{newton}
Let $(E,\|.\|)$ be a Banach space over a valued field
$(\K,|.|)$.
Let $r>0$, $x\in E$,
and $f\colon  B_r^E(x)\to E$ be a mapping.
We suppose that there exists $A\in \GL(E):=\cL(E)^\times$
such that
%
\begin{equation}\label{shrink}
\sigma:=\sup\left\{
\frac{\|f(z)-f(y)-A.(z-y)\|}{\|z-y\|}\colon 
\;\mbox{$y,z\in B_r(x)$, $y\not=z$}\right\}<\frac{1}{\|A^{-1}\|}\,.
\end{equation}
Then the following holds:
\begin{itemize}
\item[\rm (a)]
$f$ has open image and is a homeomorphism
onto its image.
\item[\rm (b)]
The inverse map $f^{-1}\colon f(B_r(x))\to B_r(x)$
is Lipschitz, with
%
\begin{equation}\label{lipofinv}
\Lip(f^{-1})\; \leq\; \frac{1}{\|A^{-1}\|^{-1}-\sigma}\,.
\end{equation}
\item[\rm (c)]
Abbreviating $a :=\|A^{-1}\|^{-1}-\sigma>0 $ and
$b :=\|A\|+\sigma$, we have
%
\begin{equation}\label{quasiisom5}
a \|z-y\|\;\leq\; \|f(z)-f(y)\|\; \leq\;
b \|z-y\|\quad
\mbox{for all $y,z\in B_r(x)$.}
\end{equation}
\item[\rm (d)]
The following estimates for the size of images
of balls are available:
For every $y\in B_r(x)$ and $s\in \;]0,r-\|y-x\|]$,
%
\begin{equation}\label{schachtel5}
B_{a s}(f(y))\; \sub \; f(B_s(y))\; \sub \; B_{b s}(f(y))
\end{equation}
holds.
In particular,
$B_{a r}(f(x)) \sub f(B_r(x)) \sub B_{b r}(f(x))$.
\end{itemize}
\end{thm}
%
%
\begin{rem}\label{againlp}
Note that the condition (\ref{shrink})
means that the
remainder term
\[
\tilde{f} \colon B_r(x)\to E,\quad \tilde{f}(y):=
f(y)-f(x)-A.(y-x)
\]
in the affine-linear approximation $f(y)=f(x)+A.(y-x)+\tilde{f}(y)$
is a Lipschitz map, with $\Lip(\tilde{f})=\sigma<\|A^{-1}\|^{-1}$.
\end{rem}
\begin{rem}
To understand the constants in Theorem~\ref{newton} better,
we recall that $\|A^{-1}\|^{-1}$
can be interpreted as a minimal
distortion factor, in the following sense:
For each $u\in E$, we have
$\|u\|=\|A^{-1}.(A.u)\|\leq \|A^{-1}\|\cdot
\|A.u\|$ and thus
%
\begin{equation}\label{leastdisto}
\|A.u\|\;\geq\; \|A^{-1}\|^{-1}\|u\|\quad
\mbox{for all $u\in E$.}
\end{equation}
Thus $A$ increases the norm of each element
by a factor of at least $\|A^{-1}\|^{-1}$.
Furthermore, $\|A^{-1}\|^{-1}$
is maximal among such factors,
as one verifies by going backwards
through the preceding lines.
Similarly, since $A^{-1}B_s(0)\sub B_{\|A^{-1}\|s}(0)$
and thus $B_s(0)\sub A.B_{\|A^{-1}\|s}(0)$
for each $s>0$, we find that
%
\begin{equation}\label{leastinfla}
A.B_s(0)\;\supseteq \; B_{\|A^{-1}\|^{-1}s}(0)\quad
\mbox{for all $s>0$.}
\end{equation}
\end{rem}
%
%
\begin{rem}\label{furthests}
The proof of Theorem~\ref{newton} also provides
the following information.
Set
$\alpha := a \|A^{-1}\|=
1-\sigma\|A^{-1}\|\in \;]0,1]$ and $\beta :=
1+\sigma\|A^{-1}\|\in [1,2[$.
Then $\beta\leq b\|A^{-1}\|$ and
%
\begin{equation}\label{quasiisom}
\alpha \|z-y\|\leq \|A^{-1}.f(z)-A^{-1}.f(y)\|\leq \beta \|z-y\|\quad
\mbox{for all
$y,z\in B_r(x)$.}
\end{equation}
For every $y\in B_r(x)$ and $s\in \;]0,r-\|y-x\|]$, we have
%
\begin{equation}\label{schachtel}
f(y)+A.B_{\alpha s}(0)\; \sub \; f(B_s(y))\; \sub \;
f(y)+A.B_{\beta s}(0)\,.
\end{equation}
Here $\alpha,\beta\to 1$ as $\sigma\to 0$.
\end{rem}
{\bf Proof of Theorem~\ref{newton}.}
(c) Given $y,z\in B_r(x)$, we have
\begin{eqnarray*}
\|f(z)- f(y)\| &= & \|f(z)-f(y)-A.(z-y)+A.(z-y)\|\\
&\leq & \|f(z)-f(y)-A.(z-y)\|+\|A.(z-y)\|\\
&\leq & (\sigma+\|A\|)\|z-y\|\; =\; b\|z-y\|
\end{eqnarray*}
and
\begin{eqnarray*}
\|z-y\|&= &\|A^{-1}.(f(z)-f(y)-A.(z-y))-(A^{-1}.f(z)-A^{-1}.f(y))\|\\
&\leq&
\|A^{-1}\|\cdot \|f(z)-f(y)-A.(z-y)\|+\|A^{-1}.f(z)-A^{-1}.f(y)\|\\
&\leq & \sigma\|A^{-1}\|\cdot \|z-y\|+\|A^{-1}\|\cdot\|f(z)-f(y)\|\,,
\end{eqnarray*}
whence (\ref{quasiisom5}) holds.
Likewise, (\ref{quasiisom}) from Remark~\ref{furthests}
follows from
\begin{eqnarray*}
\|A^{-1}.f(z)-A^{-1}.f(y)\| &= & \|A^{-1}.(f(z)-f(y)-A.(z-y))+z-y\|\\
&\leq &\|A^{-1}\|\cdot \|f(z)-f(y)-A.(z-y)\|+\|z-y\|\\
&\leq & (\sigma\|A^{-1}\|+1)\|z-y\|=\beta \|z-y\|
\end{eqnarray*}
and
\begin{eqnarray*}
\|z-y\|& \leq&
\|A^{-1}\|\cdot \|f(z)-f(y)-A.(z-y)\|+\|A^{-1}.f(z)-A^{-1}.f(y)\|\\
&\leq & \sigma\|A^{-1}\|\cdot \|z-y\|+\|A^{-1}.f(z)-A^{-1}.f(y)\|\,.
\end{eqnarray*}

(b) As a consequence of (\ref{quasiisom5}),
$f$ is injective,
a homeomorphism onto its image,
and $\Lip(f^{-1})\leq a^{-1}=(\|A^{-1}\|^{-1}-\sigma)^{-1}$.

(d) Suppose that $y\in B_r(x)$
and $s\in \;]0,r-\|y-x\|]$.
By (\ref{quasiisom5}),
we have $f(B_s(y))\sub B_{bs}(f(y))$,
proving the second half of (\ref{schachtel5}).
The second half of (\ref{schachtel})
can be shown similarly:
By (\ref{quasiisom}),
we have $A^{-1}.f(B_s(y))\sub A^{-1}.f(y)+B_{\beta s}(0)$
and thus
$f(B_s(y))\sub f(y)+A.B_{\beta s}(0)$.
We now show the first half of (\ref{schachtel}),
namely
%
\begin{equation}\label{sfhf}
f(y)+A.B_{\alpha s}(0)\; \sub \; f(B_s(y))\,.
\end{equation}
Then also the first half of (\ref{schachtel5})
will hold, as
$A.B_{\alpha s}(0)\supseteq B_{\|A^{-1}\|^{-1}\alpha s}(0)
=B_{as}(0)$
by (\ref{leastinfla}).
To prove (\ref{sfhf}),
let $c \in f(y)+A.B_{\alpha s}(0)$.
There exists $t\in \;]0,1[$ such that $c \in f(y)+A.\wb{B}_{t\alpha s}(0)$.
For $v\in \wb{B}_{st}(y)$,
we define
\[
g(v)\, :=\, v-A^{-1}.(f(v)-c)\,.
\]
Then $g(v)\in \wb{B}_{st}(y)$, because
\begin{eqnarray*}
\|g(v)-y\| &\leq &
\underbrace{\|v-y-A^{-1}.f(v)+A^{-1}.f(y)\|}_{\leq \|A^{-1}\|\sigma \|v-y\|
\leq \|A^{-1}\|\sigma st}
+ \underbrace{\|A^{-1}.c -A^{-1}.f(y)\|}_{\leq t\alpha s}\\
&\leq & (\|A^{-1}\|\sigma + \alpha )st=st\,.
\end{eqnarray*}
Thus $g(\wb{B}_{st}(y))\sub \wb{B}_{st}(y)$.
The map $g\colon  \wb{B}_{st}(y)\to \wb{B}_{st}(y)$
is a contraction, since
\begin{eqnarray}
\|g(v)-g(w)\| & = & \|v-w-A^{-1}.(f(v)-f(w))\|\nonumber\\
&\leq & \|A^{-1}\|\cdot \|f(v)-f(w)-A.(v-w)\|\nonumber\\
& \leq & \sigma \cdot \|A^{-1}\|\cdot \|v-w\| \label{hencecontr}
\end{eqnarray}
for all $v,w\in \wb{B}_{st}(y)$,
where $\sigma\|A^{-1}\|<1$.
By Banach's Contraction Theorem (Lemma~\ref{banachfix}),
there exists a unique element $v_0\in
\wb{B}_{st}(y)$ such that $g(v_0)=v_0$
and hence $f(v_0)= c$.

(a) We have already seen that $f$ is a homeomorphism onto
its image. As a consequence of (d), the image of $f$
is open.\vspace{3mm}\Punkt

\begin{rem}
Let $E$ be a Banach space, $U\sub E$ be an open subset,
$x\in U$ and $f\colon U\to E$
be a map which is strictly differentiable
at $x$ (for example,
an $SC^1$-map, $LC^1$-map
or $C^2$-map).
If $f'(x)\in \GL(E)$,
then the hypothesis of Theorem~\ref{newton}
is satisfied on $B_r(x)$
for some $r>0$,
with $A :=f'(x)$.
\end{rem}
Stronger results are available
for ultrametric Banach spaces.
In this case, $f$
behaves like an affine-linear
map, as far as the distortion of
balls is concerned:
%
%
%
\begin{thm}[Ultrametric Lipschitz Inverse Function Theorem]\label{newtonu}
\hfill Let\linebreak
$(E,\|.\|)$ be an ultrametric
Banach space over an ultrametric field
$(\K,|.|)$.
Let $r>0$, $x\in E$,
and $f\colon  B_r(x)\to E$ be a mapping.
We suppose that there exists $A\in \GL(E)$
such that
\begin{equation}\label{shrinku}
\sigma:=\sup\left\{
\frac{\|f(z)-f(y)-A.(z-y)\|}{\|z-y\|}\colon 
\;\mbox{$y,z\in B_r(x)$, $y\not=z$}\right\}<\frac{1}{\|A^{-1}\|}\,.
\end{equation}
Then the following holds:
\begin{itemize}
\item[\rm (a)]
$A^{-1}\circ f\colon  B_r(x)\to E$ is an isometry
onto an open subset of~$E$.
\item[\rm (b)]
$f$ is Lipschitz, with $\Lip(f)\leq \|A\|$,
and $f^{-1}\colon f(B_r(x))\to B_r(x)$
is Lipschitz, with $\Lip(f^{-1})\leq \|A^{-1}\|$.
\item[\rm (c)]
For all $y,z\in B_r(x)$, we have
%
\begin{equation}\label{updown6}
\|A^{-1}\|^{-1}\cdot \|z-y\| \;\leq\;
\|f(z)-f(y)\|\;\leq\;
\|A\|\cdot \|z-y\|\,.
\end{equation}
\item[\rm (d)]
For each $y\in B_r(x)$ and $s\in \,]0,r]$,
we have $B_s(y)\sub B_r(x)$ and
\begin{equation}\label{dobetter}
f(B_s(y)) \; = \; f(y)+A.B_s(0).
\end{equation}
\end{itemize}
\end{thm}
\begin{proof}
(a) For all $y,z\in B_r(x)$ with $y\not=z$,
we have
%
\begin{eqnarray*}
\|A^{-1}.f(z)-A^{-1}.f(y)-(z-y)\|
& \leq & \|A^{-1}\|\cdot \|f(z)-f(y)-A.(z-y)\|\\
& <& \|z-y\|\; ,
\end{eqnarray*}
using (\ref{shrinku})
to obtain the final inequality.
Hence, the norm~$\|.\|$ being ultrametric,
we must have $\|A^{-1}.f(z)-A^{-1}.f(y)\|=\|z-y\|$.
Thus $A^{-1}\circ f$ is in fact isometric.
As a consequence of Theorem~\ref{newton},
$A^{-1}\circ f$ has open image.

(b) Since $f=A\circ (A^{-1}\circ f)$
where $A^{-1}\circ f$ is an isometry,
$f$ is Lipschitz with $\Lip(f)\leq\Lip(A)=\|A\|$.
Likewise, $f^{-1}= (A^{-1}\circ f)^{-1}\circ A^{-1}|_{f(B_r(x))}$
is Lipschitz, with $\Lip(f^{-1})\leq \Lip(A^{-1})
=\|A^{-1}\|$.

(c) is a mere reformulation of (b).

(d) If $y\in B_r(x)$ and $s\in \,]0,r]$,
then $B_s(y)\sub B_r(y)=B_r(x)$,
as $\|.\|$ is ultrametric.
The map $A^{-1}\circ f$ being isometric,
we have $f(B_s(y))=A.(A^{-1}\circ f)(B_s(y))\sub
A.B_s(A^{-1}.f(y))=f(y)+A.B_s(0)$.
If $c \in f(y)+A.B_s(0)$ is given, define
\[
g(z)\, :=\, z-A^{-1}.(f(z)-c)\;\;\mbox{for
$z\in B_s(y)$.}
\]
Then
\begin{eqnarray*}
\|g(z)-y\|&=&\|(z-y)-(A^{-1}.f(z)-A^{-1}.f(y))+A^{-1}.(c-f(y))\|\\
&\leq & \max\{\|z-y\|,\|A^{-1}.f(z)-A^{-1}.f(y)\|,\|A^{-1}.(c-f(y))\|\}
<s
\end{eqnarray*}
for $z\in B_s(y)$, whence $g(z)\in B_s(y)$.
The map $g\colon  B_s(y)\to B_s(y)$
is a contraction, by the calculation from (\ref{hencecontr}).
Recall that, the norm on $E$ being ultrametric,
the open ball $B_s(y)$ is also closed and
therefore complete in the induced metric.
By Banach's Contraction Theorem
(Lemma~\ref{banachfix}),
there is a unique element $z_0\in B_s(y)$
such that $g(z_0)=z_0$ and thus $f(z_0)=c$.
\end{proof}
The following
consequence of Theorem~\ref{newtonu} (a) and (d)
is particularly useful.
\begin{cor}
If $A$ is an isometry
in the situation of Theorem~{\rm \ref{newtonu}}
$($for example, if $A=\id_E)$,
then $f(B_r(x))=B_r(f(x))$
and $f\colon B_r(x)\to B_r(f(x))$
is an isometry.\Punkt
\end{cor}
We are now in the position to formulate
the first version of an inverse function
theorem with parameters.
The result, and its proof,
can be re-used later to prove
the corresponding results for $C^k$-maps,
$SC^k$-maps, and $LC^k$-maps.
%
%
%
%
\begin{thm}[Lipschitz Inverse Function
Theorem with Parameters]\label{newtpara}
Let $(F,\|.\|)$ be a Banach space over a valued field
$\K$, and~$P$ be a topological space.
Let $r>0$, $x\in F$,
and $f\colon  P\times B \to F$ be a continuous
mapping, where $B:=B_r^F(x)$.
Given $p\in P$, we abbreviate
$f_p:=f(p,\sbull)\colon  B\to F$.
We suppose that there exists $A\in \GL(F)$
such that
%
\begin{equation}\label{shrinkpara}
\sigma:=\sup\left\{
\frac{\|f_p(z)-f_p(y)-A.(z-y)\|}{\|z-y\|}\colon 
\;\mbox{$p\in P$, $y,z\in B$, $y\not=z$}\right\}<\frac{1}{\|A^{-1}\|}\,.
\end{equation}
Then the following holds:
\begin{itemize}
\item[\rm (a)]
$f_p(B)$ is open in~$F$ and $f_p|_B$ is a homeomorphism
onto its image, for each $p\in P$.
\item[\rm (b)]
The set
$W:=\bigcup_{p\in P}\, \{p\}\times f_p(B)$
is open in $P\times F$, and the map
$\psi\colon  W\to F$,
$\psi(p,z):=(f_p|_B^{f_p(B)})^{-1}(z)$ is continuous.
\item[\rm (c)]
The map
$\xi \colon  P\times B\to W$,
$\xi (p,y):=(p,f(p,y))$
is a homeomorphism, with inverse given by
$\xi^{-1}(p,z)=(p,\psi(p,z))$.
\end{itemize}
If $P$ is a subset of a topological
$\K$-vector space~$E$ here
and $f$ is Lipschitz continuous,
then also $\psi$, $\xi$,
$\xi^{-1}$ and each of the maps
$(f_p|_B)^{-1}$ are Lipschitz continuous.
\end{thm}
\begin{proof}
(a)--(c): By Theorem~\ref{newton}, applied to~$f_p$,
the set $f_p(B)$ is open in~$F$ and $f_p|_B$
a homeomorphism onto its image.
Define $\alpha  :=1-\sigma\|A^{-1}\|$.
Let us show openness of $W$ and
continuity of $h$.
If $(p,z)\in W$, there exists $y\in B$
such that $f_p(y)=z$.
Let $\ve \in\,]0,r-\|y-x \|]$ be given.
There is an open neighborhood~$Q$ of~$p$ in~$P$
such that $f(q,y)\in f(p,y)+A.B_{\frac{\alpha \ve}{2}}(0)$
for all $q\in Q$, by continuity of~$f$.
Then, as a consequence of Remark~\ref{furthests},
Equation~(\ref{schachtel}),
\[
f_q(B_\ve(y))\supseteq
f(q,y)+A.B_{\alpha \ve}(0)\supseteq
f(p,y)+A.\wb{B}_{\frac{\alpha \ve}{2}}(0)
=z+A.\wb{B}_{\frac{\alpha \ve}{2}}(0)\,.
\]
By the preceding, $Q\times (z+A.\wb{B}_{\frac{\alpha \ve}{2}}(0))\sub W$,
whence $W$ is a neighborhood of $(p,z)$.
Furthermore, $\psi(q,z')=(f_q)^{-1}(z')\in B_\ve(y)=B_\ve((f_p)^{-1}(z))=
B_\ve(\psi(p,z))$
for all $(q,z')$ in the neighborhood
$Q\times (z+A.\wb{B}_{\frac{\alpha \ve}{2}}(0))$ of $(p,z)$.
Thus $W$ is open and $\psi$
is continuous. The assertions concerning~$\xi$ follow immediately.

\emph{Final assertions.}
Assume now that $P$
is a subset of a topological $\K$-vector space~$E$
and $f$ is Lipschitz continuous.
We only need
to show that $\psi$ is Lipschitz continuous
(the remaining assertions are then immediate).
To this end, let
$(p,z)\in W$.
Let $y$, $\ve$ and~$Q$ be as before.
By continuity of $f$,
after shrinking $Q$ and $\ve$,
we may assume that
%
\begin{equation}\label{nxtufl}
\|A^{-1}\|\cdot \|f_q(v)-f_p(v)\|\,\leq\, \frac{1-\theta}{8}\frac{\ve}{2}
\quad\mbox{for all $q\in Q$ and $v\in \wb{B}_{\frac{\ve}{2}}(y)$.}
\end{equation}
Abbreviate $\theta:=\sigma\|A^{-1}\|<1$.
Our goal is to see that $\psi$ is Lipschitz
continuous on the
neighborhood
$Q\times (z+A\wb{B}_{(1-\theta)\ve/4}(0))$
of $(p,z)$, where $1-\theta=\alpha$ from above.
To achieve this, for each $(q,c)$ in this neighborhood
we interpret
$\psi(q,c)$ as a fixed point
of a suitable contraction $g_{(q,c)}$
and then apply Proposition~\ref{preparadep}
on the Lipschitz continuous dependence
of fixed points on parameters.
The Lipschitz continuous map
\[
g\colon Q\times \big(z+A\wb{B}_{(1-\theta)\ve/4}(0)\big)
\times \wb{B}_{\frac{\ve}{2}}(y)\to
F \,,\quad
g(q,c,v)\, :=\,
v-A^{-1}.(f_q(v)-c)
\]
will serve our purpose.
Note first that
for all $(q,c)\in
Q\times (z+A\wb{B}_{(1-\theta)\ve/4}(0))$,
the map $g_{(q,c)}:=g(q,c,\sbull)\colon \wb{B}_{\frac{\ve}{2}}(y)\to F$
satisfies
%
\begin{eqnarray}
\|g_{(q,c)}(v)-g_{(q,c)}(w)\|
& = & \|g(q,c,v)-g(q,c,w)\|\notag\\
&=& \|A^{-1}.f_q(v)-A^{-1}.f_q(w) -(v-w)\|\notag\\
&\leq &\|A^{-1}\|\,\|f_q(v)-f_q(w)-A.(v-w)\|\notag \\
&\leq & \|A^{-1}\|\sigma\|v-w\|\;=\; \theta\|v-w\|
\,;\label{gitcnt}
\end{eqnarray}
we are thus dealing with a uniform family
of contractions.
Each $g_{(q,c)}$ is a self-map of $\wb{B}_{\frac{\ve}{2}}(y)$,
because, for each $v\in \wb{B}_{\frac{\ve}{2}}(y)$,
%
\begin{eqnarray}
\lefteqn{\|g(q,c,v)-y\|}\quad \notag\\
&=& \|g(q,c,v)-g(p,z,y)\|\notag\\
&\leq &
\|g(q,c,v)-g(p,c,v)\|
+\|g(p,c,v)-g(p,z,v)\|
+\|g(p,z,v)-g(p,z,y)\|\notag\\
&\leq & \|A^{-1}.(f_q(v)-f_p(v))\|
+\|A^{-1}.(c-z)\|+
\theta\|v-y\|\notag \\
&\leq & {\textstyle
\frac{1-\theta}{8}\frac{\ve}{2}
+ \frac{(1-\theta)\ve}{4} + \theta \frac{\ve}{2}
\;\leq \;\frac{\ve}{2}\,.}\label{callstr}
\end{eqnarray}
Here, we used (\ref{gitcnt}) to pass to the penultimate
line and then (\ref{nxtufl}).
By Banach's Contraction Theorem
(Lemma~\ref{banachfix}) and Proposition~\ref{preparadep},
$g_{(q,c)}$ has a unique fixed point
$\phi(q,c)$ for each $(q,c)$,
and the map $\phi\colon
Q\times (z+A\wb{B}_{(1-\theta)\ve/4}(0))\to F$
is Lipschitz continuous.
But $g_{(q,c)}(v)=v$ if and only if $f_q(v)=c$,
i.e., if and only if $v=\psi(q,c)$.
Thus $\phi(q,c)=\psi(q,c)$
and thus $\psi$ is Lipschitz continuous
on $Q\times (z+A\wb{B}_{(1-\theta)\ve/4}(0))$.
\end{proof}
%
%
\begin{rem}\label{facreus}
Note that if $\ve'\in [\frac{\ve}{2},\ve]$,
$q\in Q$ and $c\in z+A\wb{B}_{(1-\theta)\ve'/4}(0)$,
then the calculation (\ref{callstr})
shows that $\|g(q,c,v)-y\|\leq\frac{\ve'}{2}$.
This will be useful later.
\end{rem}
As an immediate consequence, we obtain an implicit function
theorem.
%
\begin{cor}[Lipschitz Implicit Function Theorem]\label{lipimpl}
In the situation\linebreak
of Theorem~{\rm\ref{newtpara}},
let $(p_0,y_0)\in P\times B$.
Then there exists an open neighborhood
$Q\sub P$ of~$p_0$ such that
$z_0:=f(p_0,y_0)\in f_p(B)$
for all $p\in Q$.
The mapping $\lambda \colon Q\to B$, $\lambda(p):=\psi(p,z_0)$
is continuous
$($resp., Lipschitz continuous if so is $f)$,
satisfies $\lambda (p_0)=y_0$, and
\[
\{(p,y)\in Q\times B \colon f(p,y)=z_0\}\;=\;
\graph\,(\lambda )\,.
\]
\end{cor}
\begin{proof}
Because $W$ is an open neighborhood
of $(p_0,z_0)$ in $P\times F$,
there exists an open neighborhood~$Q$
of $p_0$ in $P$
such that $Q\times \{z_0\}\sub W$.
Then $\lambda(p):=\psi(p,z_0)$ makes sense
for all $p\in Q$.
The rest is now obvious from Theorem~\ref{newtpara}.
\end{proof}
\subsection*{Inverse Function Theorem with Parameters}
We are now in the position
to formulate and prove
our main result,
an Inverse Function Theorem
with Parameters for various types
of differentiable mappings.
%
%
%
\begin{thm}[Inverse Function Theorem with Parameters]\label{advif}
Let $\K$ be a valued field,
$k\in \N\cup\{\infty\}$,
$E$ be a topological $\K$-vector
space, and $F$ be a Banach space over~$\K$.
Let $P_0 \sub E$ be a subset
with dense interior, $U\sub F$
be open, and $f\colon  P_0 \times U \to F$ be a map.
Assume that
\begin{itemize}
\item[\rm (i)]
$f$ is $LC^k$, respectively, $SC^k$; or:
\item[\rm(ii)]
$f$ is $C^k$ and $k\geq 2$,
or $f$ is $C^k$,
$\K$ is locally compact and $F$ has finite dimension.
\end{itemize}
Abbreviate $f_p:=f(p,\sbull)\colon  U \to F$ for $p\in P_0$.
Suppose that $(p_0,x_0)\in P_0\times U$ is given such that
$f_{p_0}'(x_0)
\in \GL(F)$.
Then there exists an open neighborhood $P\sub P_0$ of~$\,p_0$
and $r>0$ such that $B:=B_r(x_0)\sub U$ and the following holds:
\begin{itemize}
\item[\rm (a)]
$f_p(B)$ is open in~$F$, for each $p\in P$,
and $\phi_p\colon  B\to f_p(B)$,
$\phi_p(x):=f_p(x)=f(p,x)$ is an
$LC^k$-diffeomorphism $($resp.,
an $SC^k$-diffeomorphism; resp., a
$C^k$-diffeomorphism$)$.
\item[\rm (b)]
$W:=\bigcup_{p\in P} (\{p\}\times f_p(B))$ is open in $P_0\times F$,
and the map
\[
\psi\colon  W\to B, \quad
\psi(p,z)\, :=\, \phi_p^{-1}(z)
\]
is~$LC^k$ $($resp., $SC^k$, resp., $C^k)$.
Furthermore, the map
\[
\xi \colon  P \times B\to W, \quad
\xi (p,x):=(p, f(p,x))
\]
is an $LC^k$-diffeomorphism
$($resp., an $SC^k$-diffeomorphism, resp.,
a $C^k$-diffeomorphism$)$,
with inverse $\xi^{-1}(p,z)=(p,\psi(p,z))$.
\item[\rm (c)]
$P\times B_\delta(f_{p_0}(x_0)) \sub W$ for some $\delta>0$.
\end{itemize}
In particular, for each $p\in P$ there is a unique element
$\lambda(p)\in B$ such that $f(p,\lambda(p)) =f(p_0,x_0)$,
and the map $\lambda\colon  P \to B$
so obtained is~$LC^k$ $($resp., $SC^k$, resp., $C^k)$.
\end{thm}
%
%
\begin{rem}\label{additclaim}
Given $\alpha ,\beta \in \R$ such that $0<\alpha<1<\beta$
in the situation of Theorem~\ref{advif},
one can furthermore achieve that
%
\begin{equation}\label{eqaddclaim}
f_p(x)+ A.B_{\alpha s}(0) \; \sub \; f_p(B_s(x))\; \sub\;
f_p(x)+A.B_{\beta s}(0)
\end{equation}
for all $p\in Q$, $x\in B$
and $s\in \,]0,r-\|x-x_0\|]$.
\end{rem}
\begin{proof}
Given $\alpha , \beta \in \R$ such that $0<\alpha <1< \beta $,
define
\[
\tau \; :=\;  \min
\left\{ {\textstyle \frac{\beta -1}{\|A^{-1}\|},
\frac{1-\alpha}{\|A^{-1}\|}}\right\}
\, <\, {\textstyle \frac{1}{\|A^{-1}\|}}\, ,
\]
where $A:=f_{p_0}'(x_0)$.
Then $1-\tau \|A^{-1}\|\geq \alpha $ and $1+\tau \|A^{-1}\|\leq \beta$.
By strict differentiability of~$f$ at $(p_0,x_0)$
(resp., by Remark~\ref{varC1paraunif}),
there exists
an open neighborhood $P\sub P_0$ of~$p_0$ and $r>0$ such that
$B:=B_r(x)\sub U$ and
%
%
\begin{equation}\label{eqwithc}
\|f_p(z)-f_p(y)-f_{p_0}'(x).(z-y)\|
\; \leq \;
\tau \, \|z-y\|
\end{equation}
for all $p\in P$ and $y\not=z\in B$.
Hence
\begin{eqnarray}
\sigma & \! := \! & \sup\left\{
\frac{\|f_p(z)-f_p(y)-f_{p_0}'(x_0).(z-y)\|}{\|z-y\|}\colon 
\,\mbox{$p\in P$, $z\not=y\in B$} \right\}\notag\\
&\! \leq \! & \tau  \, < \, \frac{1}{\|A^{-1}\|}\,.\label{thirty1}
\end{eqnarray}
Thus Theorem~\ref{newtpara} applies to $f|_{P\times B}$
with $A:=f'_{p_0}(x_0)$,
whence $f_p(B)$ is open in~$F$ and $\phi_p:=f_p|_B^{f_q(B)}$ a
homeomorphism onto its image, for each $p\in P$;
the set
$W:=\bigcup_{p\in P}\{p\}\times f_p(B)$
is open in $P_0\times F$;
the map $\psi\colon  W\to B$, $\psi(p,z):=\phi_p^{-1}(z)$ is continuous;
and the mapping
$\xi \colon  P\times B\to W$, $\xi(p,y):=(p,f(p,y))$
is a homeomorphism, with inverse given by
$\xi^{-1}(p,z)=(p,\psi(p,z))$.
In view of (\ref{thirty1}),
Lemma~\ref{newton} applies to
$f_p|_B$, for all $p\in P$, whence~(\ref{eqaddclaim})
in Remark~\ref{additclaim}
holds.\\[3mm]
Also (c) is easily established:
we set $\delta:=\|A^{-1}\|^{-1} \frac{\alpha r}{2}$.
After shrinking $P$, we may assume that $\|f(p,x_0)-f(p_0,x_0)\|<\delta$
for all $p\in P$.
Then, using (\ref{eqaddclaim})
with $x:=x_0$ and $s:=r$,
we get $f_p(B)\supseteq f_p(x_0)+A.B_{\alpha r}(0)
\supseteq
B_{2\delta}(f_p(x_0))
\supseteq B_\delta(f_{p_0}(x_0))$,
for all $p\in P$.
Thus (c) holds.\\[3mm]
(a) and (b): If can show that
$\psi$ is $LC^k$ (resp., $SC^k$, resp., $C^k$),
then clearly all of the maps
$\psi$, $\xi$, $\lambda$ and
$\phi_q$ will have the desired
properties.
Given $(p,z)\in W$,
we let $A:=f_{p_0}'(x_0)$ as before and
define $y$, $\theta$, $Q$, $\ve$
and the map
\[
g\colon Q\times \big(z+A B_{(1-\theta)\ve/4}(0)\big)
\times B_{\frac{\ve}{2}}(y)\to F \,,\quad
g(q,c,v)\, :=\,
v-A^{-1}.(f_q(v)-c)
\]
as in the proof of Theorem~\ref{newtpara}
(using now open balls instead of closed balls).
Because the arguments from the proof of
Theorem~\ref{newtpara} apply to the restriction
of $g$ to
$Q\times \big(z+A \wb{B}_{(1-\theta)\ve'/4}(0)\big)
\times \wb{B}_{\frac{\ve'}{2}}(y)$
for each $\ve'\in [\frac{\ve}{2},\ve[$
(see Remark~\ref{facreus}),
we deduce that $g(q,c,\sbull)$
has a fixed point in
$B_{\frac{\ve}{2}}(y)$,
for each
$q\in Q$
and each\linebreak
$c\in z+A B_{(1-\theta)\ve/4}(0)$.
Repeating the arguments used in the proof of
Theorem~\ref{newtpara},
we see that
$(g(q,c,\sbull))_{q,c}$ is a uniform family
of contractions
for $q\in Q$, $c\in z+A B_{(1-\theta)\ve/4}(0)$,
and that
$\psi(q,c)$ is the unique fixed point
of the contraction $g(q,c,\sbull)\colon B_{\frac{\ve}{2}}(y)\to
B_{\frac{\ve}{2}}(y)$.
Since $g$ is $LC^k$
(resp., $SC^k$, resp., $C^k$),
Proposition~\ref{preparadep} shows that
$\psi$ is $LC^k$ (resp., $SC^k$, resp., $C^k$)
on the open neighborhood
$Q\times \big(z+A B_{(1-\theta)\ve/4}(0)\big)$
of $(p,z)$ in $W$.
This completes the proof.
\end{proof}
%
%
\begin{rem}\label{loosee}
If $f$ is $C^1$ in the situation of Theorem~\ref{advif}
but $\K$ is not locally compact or
$F$ is infinite-dimensional,
then the conclusions of the theorem
still remain intact if
we assume that (\ref{shrinkpara})
is satisfied by $f$ with $A:=f'_{p_0}(x_0)$.
This ensures that we are in the situation of
Theorem~\ref{newtpara},
and we can now complete the proof
(with $\alpha:=1-\sigma\|A^{-1}\|$ and $\beta=
1+\sigma\|A^{-1}\|$) as before,
noting that the $C^1$-dependence of fixed points
on parameters established in Theorem~\ref{fixpardep}
requires neither local compactness of~$\K$
nor finite-dimensionality of~$F$.

The same reasoning shows that if $f$ is $C^k$,
$SC^k$ or $LC^k$ for some $k\in \N\cup\{\infty\}$
in the situation of the Lipschitz Inverse Function
Theorem (Theorem~\ref{newton}),
then $f^{-1}\colon f(B_r(x))\to B_r(x)$
is $C^k$ (resp., $SC^k$, resp., $LC^k$).
\end{rem}
%
%
\begin{rem}\label{proofstanda}
Note that Theorem~\ref{advif}
subsumes as its final assertion
the Generalized Implicit Function Theorem
announced above (Theorem~\ref{genimp}).
Using a singleton set of parameters,
we also obtain the ordinary Inverse Function
Theorem (Theorem~\ref{basicinv})
as a special case.
\end{rem}
For ultrametric Banach spaces,
the Inverse Function Theorem with Parameters
attains a simpler form:
%
%
\begin{thm}\label{invparultra}
Let $(\K,|.|)$ be an ultrametric field,
$k\in \N\cup\{\infty\}$,
$E$ be a topological $\K$-vector
space, and $F$ be an ultrametric Banach space over~$\K$.
Let $P_0 \sub E$ be a subset
with dense interior, $U\sub F$
be open, and $f\colon  P_0 \times U \to F$ be a map.
Assume that
\begin{itemize}
\item[\rm (i)]
$f$ is $LC^k$, respectively, $SC^k$; or:
\item[\rm(ii)]
$f$ is $C^k$ and $k\geq 2$,
or $f$ is $C^k$,
$\K$ is locally compact and $F$ has finite dimension.
\end{itemize}
Abbreviate $f_p:=f(p,\sbull)\colon  U \to F$ for $p\in P_0$.
Suppose that $(p_0,x_0)\in P_0\times U$ is given such that
$A:=f_{p_0}'(x_0)
\in \GL(F)$.
Then there exists an open neighborhood $P\sub P_0$ of~$\,p_0$
and $r>0$ such that $B:=B_r(x_0)\sub U$ and the following holds:
\begin{itemize}
\item[\rm (a)]
$f_p(B)=f(p_0,x_0)+A.B_r(0)=:V$, for each $p\in P$,
and $\phi_p\colon  B\to V$, $\phi_p(y):=f(p,y)$ is an
$LC^k$-diffeomorphism
$($resp., an $SC^k$-diffeomorphism,
resp., a $C^k$-diffeomorphism$)$.
\item[\rm (b)]
$f_p(B_s(y))=f_p(y)+A.B_s(0)$
for all $p\in P$, $y\in B$ and $s\in \,]0,r]$.
\item[\rm (c)]
The map $\psi\colon  P\times V\to B$, $\psi(p,v):=\phi_p^{-1}(v)$
is $LC^k$, $SC^k$, resp., $C^k$.
\item[\rm (d)]
$\xi \colon  P\times B\to P\times V$,
$\xi(p,y):=(p,f(p,y))$
is an $LC^k$-diffeomorphism
$($resp., $SC^k$-diffeomorphism, resp., $C^k$-diffeomorphism$)$,
with inverse given by $\xi^{-1}(p,v)=(p,\psi(p,v))$.
\end{itemize}
\end{thm}
\begin{proof}
We let $P$ and $r$ be as in the proof of Theorem~\ref{advif},
choosing $P$ so small that
%
\begin{equation}\label{lttlhlp}
\|A^{-1}.f(p,x_0)-A^{-1}.f(p_0,x_0)\|\; <\; r
\end{equation}
for all $p\in P$.
Since (b) holds by Theorem~\ref{newtonu}\,(d),
taking $s:=r$ we deduce that
\begin{eqnarray*}
f_p(B_r(x_0)) & = & f_p(x_0)+A.B_r(0)\;=\;
f_{p_0}(x_0)+(f_p(x_0)-f_{p_0}(x_0))+A.B_r(0)\\
& = &
f_{p_0}(x_0)+A.B_r(0)\; =:\; V\,.
\end{eqnarray*}
The final equality holds because
$f_p(x_0)-f_{p_0}(x_0)\in A.B_r(0)$
by (\ref{lttlhlp}) and $B_r(0)$ is an
additive subgroup of~$F$.
Thus (a) holds by Theorem~\ref{advif}\,(a).
Furthermore, $W:=\bigcup_{p\in P}\{p\}\times f_p(B)=P\times V$
by (a), whence (c) and (d) hold by Theorem~\ref{advif}\,(b).
\end{proof}
\appendix
%
\section{Appendix: {\boldmath$\cC^0$}-concepts}\label{appC0}
In this appendix, we describe a version
of the notion of ``$\cC^0$-concept''
introduced in \cite{Ber}.
While only mappings between open sets
were considered in~\cite{Ber},
we now define $\cC^0$-concepts
for mappings between subsets of topological
vector spaces with dense interior.
In other respects,
our $\cC^0$-concepts
are more restrictive than those from~\cite{Ber}.
In particular, we are working only with
Hausdorff topological vector spaces
over Hausdorff topological fields,
whereas suitable topologized modules
over suitable topologized rings
provided the general framework
in \cite{Ber}.
\begin{defn}
Let $\K$ be a (non-discrete, Hausdorff) topological
field and $\cE$ be a class of (Hausdorff)
topological $\K$-vector spaces
satisfying the following axioms:
\begin{itemize}
\item[(E1)]
$\K\in \cE$ and $\{0\}\in \cE$ hold:
\item[(E2)]
If $E\in \cE$ and $F$ is a topological
$\K$-vector space isomorphic to~$E$,
then $F\in \cE$;
\item[(E3)]
If $E_1,E_2 \in \cE$, then also $E_1\times E_2\in \cE$
(when equipped with the product topology).
\end{itemize}
A \emph{$\cC^0$-concept}
over~$\K$
(with underlying class of topological
vector spaces $\cE$)
assigns
a set $\cC^0(U,V)\sub C(U,V)$
of continuous maps
to all $E,F\in \cE$
and subsets $U\sub E$ and $V\sub F$
with dense interior,
such that the following axioms
are satisfied:
\begin{itemize}
\item[(C01)]
If $E, F, H\in \cE$
and $U\sub E$, $V\sub F$, $W\sub H$
are subsets with dense
interior, then $g\circ f \in \cC^0(U, W)$
for all $f\in \cC^0(U,V )$ and
$g\in \cC^0(V, W)$.
Furthermore, $\id_U\in \cC^0(U,U)$.
\item[(C02)]
For each $E\in \cE$ and subset $U\sub E$
with dense interior,
the inclusion map $i_U\colon U\to E$
is $\cC^0$.
\item[(C03)] 
For all $E, F\in \cE$
and subsets $U\sub E$, $V\sub F$
with dense interior,
a map $f\colon U\to V$ is $\cC^0$
if and only if it is $\cC^0$ as a map into
$F$, i.e., if and only if $i_V\circ f\in \cC^0(U,F)$.
\item[(C04)]
If $n\in \N$,
$E_1,\ldots, E_n, F\in \cE$ and $\beta\colon E_1\times\cdots\times
E_n \to F$ is a continuous
$n$-linear map, then $\beta\in \cC^0(E_1\times\cdots\times E_n,F)$.
\item[(C05)]
Given $F, E_1, E_2\in \cE$
and a subset $U\sub F$ with dense interior,
a mapping $f=(f_1,f_2)\colon U\to E_1\times E_2$
is $\cC^0$ if and only if both components
$f_1$ and $f_2$ are $\cC^0$.
\item[(C06)]
(Locality).
If $E,F\in \cE$, $U\sub E$ is a subset with dense
interior and $f\colon U\to F$ a mapping such that
$f|_{U_i}\in \cC^0(U_i,F)$
for an open cover $(U_i)_{i\in I}$ of~$U$,
then $f\in \cC^0(U,F)$.
\end{itemize}
\end{defn}
\begin{rem}
Whenever a $\cC^0$-concept is used in the present article,
$\cE$ simply is the class of all Hausdorff
topological $\K$-vector spaces.
But, of course,
there are other interesting
classes of topological $\K$-vector spaces,
for example the classes of locally convex
(polynormed, normable, resp., complete normable)
spaces over valued fields,
or classes of spaces with certain
completeness properties, etc.
\end{rem}
\begin{rem} Here are some consequences of the axioms.
\begin {itemize}
\item[(a)]
Property (C01) means that the pairs $(E,U)$,
where $E\in \cE$ and $U\sub E$ is a subset with dense
interior, form a category with $\Hom((U,E), (V,F))=\cC^0(U,V)$
as respective set of
morphisms.\footnote{We prefer to suppress $E$ and $F$ in the notation
$\cC^0(U,V)$, with little risk
of misunderstanding.}
\item[(b)]
Properties (C01) and (C02) guarantee in particular that
$f|_U=f\circ i_U \in \cC^0(U,F)$ for each $\cC^0$-map
$f\colon E\to F$ and subset $U\sub E$
with dense interior.
\item[(c)]
(C04) ensures
that continuous linear maps and continuous
bilinear maps are $\cC^0$.
Hence the addition map $E\times E\to E$
and the scalar multiplication map
$\K \times E\to E$ are $\cC^0$, for each $E\in \cE$.
\item[(d)]
Note that if $U\sub \K$ has dense
interior and $t\in \K$,
then $U\setminus \{t\}$
is dense in~$U$, whence
a $\cC^0$-map $f\colon U\to F$
is uniquely determined by its restriction
to $U\setminus\{t\}$.
This property (``determination axiom'')
plays an essential role in the discussions
of \cite{Ber} (in the case of open subsets).
\item[(e)]
Given a $\cC^0$-concept in the present sense,
restricting attention to mappings
between open subsets we obtain a $\cC^0$-concept
in the sense of \cite{Ber}.
\item[(f)]
In \cite{Ber},
Property (C04) is not required in full,
but it is satisfied by all interesting examples
of $\cC^0$-concepts based on topological vector spaces,
and makes it unnecessary to distinguish
between continuous linear maps and linear maps which are
$\cC^0$ (and similar nuisances).
\item[(g)]
Usually, one only specifies the $\cC^0$-maps
$U\to F$ on subsets $U\sub E$ with dense interior,
for all $E,F\in \cE$.
One then tacitly uses (C03)
as the \emph{definition} of
$\cC^0$-maps to a subset $V\sub F$ with dense
interior.
\end{itemize}
\end{rem}
\begin{defn}
Given a $\cC^0$-concept over a topological
field~$\K$ based on a class~$\cE$ of
topological $\K$-vector spaces,
we define $\cC^1$-maps as follows:
Let $E,F\in \cE$
and $f\colon U\to F$ be a $\cC^0$-map
on a subset $U\sub E$ with dense
interior.
We say that $f$ is $\cC^1$
if there exists a $\cC^0$-map
$f^{[1]}\colon U^{[1]}\to F$
which extends $f^{]1[}\colon U^{]1[}\to F$
(where $U^{]1[}$ and $f^{]1[}$
are as in Section~\ref{secprel}).
Recursively, we say that
$f$ is $\cC^k$ if $f$ is $\cC^1$
and $f^{[1]}$ is $\cC^{k-1}$.
\end{defn}
Then all relevant results from \cite{Ber}
(and their proofs) remain valid.
In particular,
the Chain Rule holds for $\cC^k$-maps;
being $\cC^k$ is a local property;
finite-order Taylor expansions
are available, etc.
\section{Appendix: A variant of Lemma~\ref{monstrtech}}\label{apptechni}
%
%
%
In this appendix, we prove the following variant
of Lemma~\ref{monstrtech}
for $C^1$-maps (which is not needed in the main text).
%
%
%
\begin{la}\label{monstrtech2}
Let $\K$ be a locally compact field,
$E$ and $H$ be topological $\K$-vector spaces
and $F$ be a finite-dimensional
normed $\K$-vector space.
Let $U\sub E$ and $V\sub F$
be subsets with dense interior
and
$f\colon U \times V \to H$ be a $C^1$-map.
Let $u_0\in U$, $v_0\in V$, 
$x_0\in E$, $y_0\in F$,
$\gamma$ be a gauge on~$H$,
and $\ve>0$.
Then there exist neighborhoods
$U_0 \sub U$ of~$u_0$,
$V_0\sub V$ of~$v_0$,
$X_0\sub E$ of $x_0$,
$Y_0\sub F$ of~$y_0$
and a $0$-neighborhood
$S_0\sub \K$
such that
%
\begin{equation}\label{eqmonstrt2}
\|f^{[1]}(u,v,w,y_1,t)-f^{[1]}(u,v,w,y_2,t)-f'_{u_0}(v_0).(y_1-y_2)\|_\gamma
\;\leq\; \ve\, \|y_1-y_2\|
\end{equation}
for all elements $u\in U_0$, $v\in V_0$,
$x\in X_0$, $y_1,y_2\in Y_0$ and $t\in S_0$
such that 
$(u,v,x,y_1,t), (u,v,x,y_2,t)\in (U\times V)^{[1]}$.
\end{la}
The proof of Lemma~\ref{monstrtech2}
uses a variant of Lemma~\ref{symmetries}
for $f^{[2]}$:
%
\begin{la}\label{symmetries2}
Let $E$ and $F$ be topological vector
spaces over a topological field~$\K$, and $f\colon  U\to F$
be a $C^2$-map, defined on a subset of~$E$ with dense
interior.
If $t\in \K^\times$, $x,x_1,y,y_1\in E$ and $s,s_1,s_2\in \K$ such that
\[
((x,y,ts),\, (x_1,y_1,ts_1),\, ts_2)\in U^{[2]}\,,
\]
then also
$((x,t^2y,\frac{s}{t}),\, (tx_1,t^3y_1,s_1),\, s_2)\in U^{[2]}$, and
%
\begin{equation}\label{bracket2}
{\textstyle t^3 \, f^{[2]}((x,y,ts),\, (x_1,y_1,ts_1),\, ts_2)
=f^{[2]}((x,t^2y,\frac{s}{t}),\, (tx_1,t^3y_1,s_1),\,
s_2)\,.}
\end{equation}
\end{la}
\begin{proof}
See \cite[Lemma~3.3\,(b)]{IMP}
for the case of open domains.
The proof carries over verbatim.
\end{proof}
{\bf Proof of Lemma~\ref{monstrtech2}.}
We may assume that $\|F\|\sub |\K|$
(because all norms on~$F$
are equivalent, and we can take a maximum
norm with respect to some basis).
Let
$\Omega$ be the set of all
$(u,v,x,t,y_1,y_2)\in
U\times V\times E\times\K\times F\times F$
such that
$(u,v,x,y_1,t), (u,v,x,y_2,t)\in (U\times V)^{[1]}$.
We choose a gauge $\|.\|_\zeta$
on~$H$ such that $\|a+b\|_\gamma\leq\max\{\|a\|_\zeta,\|b\|_\zeta\}$
for all $a,b\in H$.
Given $(u,v,x,t,y_1,y_2)\in\Omega$,
we have
%
%
%
\begin{eqnarray}
\lefteqn{\|f^{[1]}(u,v,x,y_1,t)-f^{[1]}(u,v,x,y_2,t)
-f'_{u_0}(v_0).(y_1-y_2)\|_\gamma}\quad\notag\\
&=&
\|f^{[1]}(u+tx,v+ty_2,0,y_1-y_2,t)
-f'_{u_0}(v_0).(y_1-y_2)\|_\gamma\notag\\
&=&
\|f^{[1]}(u+tx,v+ty_2,0,y_1-y_2,t)
-f'_{u+tx}(v).(y_1-y_2)\notag\\
& & \qquad\qquad
+(f'_{u+tx}(v)-f'_{u_0}(v_0)).(y_1-y_2)\|_\gamma\notag\\
&\leq &\max\big\{
\| f^{[1]}(u+tx,v+ty_2,0,y_1-y_2,t)-f_{u+tx}'(v).(y_1-y_2)\|_\zeta,
\notag\\
& &
\|(f'_{u+tx}(v)-f'_{u_0}(v_0)).(y_1-y_2)\|_\zeta
\big\}\,,\label{spltprb2}
\end{eqnarray}
using Lemma~\ref{firsiden} to obtain
the first equality.
Abbreviate
\[
\Omega_1 := 
\{(u,v,x,t,y_1,y_2,w,r)\in
\Omega\times F\times \K \colon
(u+tx,v+ty_2,0,w,r)\in (U\times V)^{[1]}\}.
\]
The function
\[
g\colon \Omega_1\to H\,,\quad
g(u,v,x,t,y_1,y_2, w,r ):=
f^{[1]}(u+tx,v+ty_2,0,w,r)
-f_{u+tx}'(v).w
\]
is continuous and vanishes
on the compact set
$\{(u_0,v_0,x_0,0, y_0,y_0)\}\times K\times \{0\}$,
where $K:=\wb{B}^F_1(0)$.
Hence, there exist
neighborhoods
$U_0 \sub U$ of~$u_0$,
$V_0\sub V$ of~$v_0$,
$X_0\sub E$ of $x_0$,
$Y_0\sub F$ of~$y_0$
and a balanced $0$-neighborhood
$S_0\sub \K$
such that
\[
\|g(u,v,x,t,y_1,y_2, w,r )\|_\zeta\, \leq \, \ve
\]
for all $(u,v,x,t,y_1,y_2,w,r)\in
\Omega_1\cap (U_0\times V_0\times
X_0\times S_0\times Y_0\times Y_0\times K\times S_0)=:\Omega_2$.
Given
$(u,v,x,t,y_1,y_2,w,r)\in \Omega_2$ and
$s\in \K^\times $ such that
$\|w\|\leq |s|\leq 1$,
we have
$(u,v,x,t,y_1,y_2,s^{-1}w, sr)\in
\Omega_2$
and
\begin{eqnarray}
\lefteqn{\|g(u,v,x,t,y_1,y_2,w,r)\|_\zeta}\qquad\notag\\
&=& \|f^{[1]}(u+tx,v+ty_2,0,w,r)
-f_{u+tx}'(v).w\|_\zeta\notag\\
&=&
|s|\cdot \|f^{[1]}(u+tx,v+ty_2,0,s^{-1}w,sr)
-f_{u+tx}'(v).s^{-1}w\|_\zeta \notag\\
&=&
|s|\cdot \|g(u,v,x,t,y_1,y_2,s^{-1}w, sr )\|_\zeta\leq
|s|\cdot \ve\,.\label{wtg}
\end{eqnarray}
If $\|w\|>0$, we can choose $s$ such that
$|s|=\|w\|$; if $\|w\|=0$, we can let $s$ pass
to~$0$.
In either case, (\ref{wtg}) entails that
$\|g(u,v,x,t,y_1,y_2,w,r)\|_\zeta\leq \ve \|w\|$.
After shrinking $Y_0$, we may assume
that $Y_0-Y_0\sub K$.
For all
$(u,v,x,t,y_1,y_2)\in \Omega\cap (U_0\times V_0\times X_0\times
S_0\times Y_0\times Y_0)$,
we then have
$(u,v,x,t,y_1,y_2,y_1-y_2,t)\in
\Omega_2$
and
\begin{eqnarray}
\hspace*{-18mm}\lefteqn{\|f^{[1]}(u+tx,v+ty_2,0,y_1-y_2,t)
-f_{u+tx}'(v).(y_1-y_2)\|_\zeta}\qquad \notag\\
&=&
\|g(u,v,x,t,y_1,y_2,y_1-y_2,t)\|_\zeta\, \leq \, \ve \|y_1-y_2\|\,.
\label{ztg}
\end{eqnarray}
By Lemma~\ref{operctsfin},
after shrinking $U_0, V_0, X_0, Y_0$ and $S_0$,
we can achieve that also
%
\begin{equation}\label{lstp}
\|(f'_{u+tx}(v)-f'_{u_0}(v_0)).(y_1-y_2)\|_\zeta\;\leq\;
\ve\, \|y_1-y_2\|
\end{equation}
for all $(u,v,x,t,y_1,y_2)\in \Omega\cap
(U_0\times V_0\times X_0\times S_0\times Y_0\times Y_0)$.
Combining (\ref{spltprb2}), (\ref{ztg})
and (\ref{lstp}), we now see that
(\ref{eqmonstrt2}) holds.\Punkt

%
%
%
%
%
%
\section{Strictly differentiable mappings of several\\
variables over complete valued fields}\label{newppx}
In this manuscript, finite-dimensional
vector spaces over locally compact
fields played a special role,
and it was frequently possible
to obtain extra results for $C^1$-maps
on open subsets of such spaces
(even when only $C^k$-maps with
$k\geq 2$ could be treated otherwise).
After the manuscript was completed,
the author realized that the local
compactness is inessential,
and that most of the
extra results
can be proved just as well for $C^1$-maps on
open subsets of finite-dimensional
vector spaces over complete valued
fields, using different arguments.
In this section, we explain
these further generalizations.\footnote{At a
later stage, this section
may be merged with the main text.}
We begin with an appropriate
replacement for Lemma~\ref{C1impllud}.
Recall that if $(\K,|.|)$
is a \emph{complete}
valued field, then every
finite-dimensional
(Hausdorff) topological $\K$-vector space
automatically
carries the canonical vector
topology (see Theorem~2
in \cite[Chapter~I, \S2, no.\,3]{BTV}).
%
\begin{prop}\label{newprop}
Let $(\K,|.|)$ be a
valued field,
$E$ be a finite-dimensional
$\K$-vector space,
equipped with its canonical
vector topology,
and $F$ be a topological $\K$-vector space.
Let $U\sub E$ be an open subset,
$f\colon  U\to F$ be a map,
and $k\in \N\cup\{\infty\}$.
Then $f$ is~$C^k$
if and only if $f$ is $SC^k$.
\end{prop}
\begin{proof}
We may assume
$k\in \N$, and reduce to $k=1$ by an obvious
induction.
Hence Proposition~\ref{newprop}
will hold if
we can show that $f$ is $C^1$
if and only if $f$ is strictly
differentiable.
But this equivalence
follows from Lemma~\ref{strictC1}
and the
following generalization
of Lemma~\ref{C1paraunif}
(using a singleton set
of parameters).
\end{proof}
%
%
\begin{la}\label{C1paraunifV}
Let $(\K,|.|)$ be a valued
field,
$E$ be a finite-dimensional
$\K$-vector space
equipped with the canonical
vector topology,
$U\sub E$ be an open subset,
and $\|.\|$ be a norm on~$E$
defining its topology.
Let $F$ be a topological $\K$-vector space,
$P$ be a topological space,
and $f\colon  P\times U\to F$ be a continuous
map such that $f_p:=f(p,\sbull)\colon  U\to F$
is $C^1$ for all $p\in P$, and such
that the map
\[
P\times U^{[1]}\to F,\quad
(p,y)\mto (f_p)^{[1]}(y)
\]
is continuous.
Let $p\in P$ and $u\in U$ be given.
Then, for every $\ve>0$ and gauge
$\gamma$ on~$F$, there is a neighborhood
$Q$ of~$p$ in~$P$ and $\delta>0$ such that
%
\begin{equation}\label{wannashow}
\|f_q(z)-f_q(y)-f_q'(u).(z-y)\|_\gamma\, \leq \, \ve \|z-y\|
\end{equation}
for all $q\in Q$ and $y,z\in B_\delta^E(u)\cap U$,
where $f_q'(u):=d(f_q)(u,\sbull)$.
\end{la}
\begin{proof}
We may assume that $E=\K^n$ for some $n\in \N$.
Also, we may assume that
$\|.\|=\|.\|_\infty$
is the maximum norm (as all norms
giving rise to the canonical
vector topology are equivalent).
Let $p\in P$, $u=(u_1,\ldots, u_n)\in U$,
$\ve>0$ and a gauge
$\gamma$ on~$F$ be given.
After shrinking~$U$, we may
assume that $U=U_1\times\cdots\times U_n$
for certain open neighbourhoods
$U_1,\ldots, U_n\sub \K$
of $u_1,\ldots, u_n$,
respectively.
As a consequence of Lemma~\ref{substitut},
there exists a gauge~$\eta$ on~$F$
such that
\begin{equation}\label{itertriv}
\|x_1+\cdots+x_{2n}\|_\gamma
\;\leq\;
\max_{j=1,\ldots,2n} \|x_j\|_\eta
\quad\mbox{for all $\, x_1,\ldots, x_{2n}\in F$.}
\end{equation}
Let $e_1,\ldots, e_n$ be the
standard basis for $\K^n$,
with components $(e_j)_i=\delta_{ij}$.
We shall use the
partial derivatives
$D_jf(q,x):=d(f_q)(x,e_j)=(f_q)^{[1]}(x,e_j,0)$;
then
$D_jf\colon P\times U\to F$ is continuous.
Define $R\colon P\times U^{[1]}\to F$,
\[
R(q,x,y,t)\; :=\; f_q^{[1]}(x,y,t)-f_q^{[1]}(x,y,0)
\]
for $q\in P$, $(x,y,t)\in U^{[1]}$.
Then $R$ is continuous
and
\begin{equation}\label{agtr}
R(q,x,y,0)\;=\; 0\quad
\mbox{for all $\,q\in P$,
$x\in U$ and $y\in E$.}
\end{equation}
Furthermore,
%
\begin{equation}\label{nxttrv}
f(q,x+ty)-f(q,x)- f_q'(x).ty
\;=\;tR(q,x,y,t)\quad
\mbox{for $q\in P$, $(x,y,t)\in U^{[1]}$.}
\end{equation}
Since
$D_jf$ and $R$
are continuous
and $R(p,u,e_j,0)=0$
for $j\in \{1,\ldots, n\}$,
we find a neighbourhood
$Q\sub P$ of~$p$
and $\delta>0$
such that
%
\begin{equation}\label{whatn1}
\|D_jf(q,v)-D_jf(q,w)\|_\eta\;\leq\;\ve
\end{equation}
for all
$j\in\{1,\ldots, n\}$, $q\in Q$
and $v,w\in B_\delta^E(u)\cap U$;
and such that
\begin{equation}\label{whatn2}
\|R(q,v,e_j,t)\|_\eta\;\leq \;\ve
\end{equation}
for all $j\in\{1,\ldots, n\}$, $q\in Q$, $v \in B_\delta^E(u)\cap U$
and $t\in B_{2\delta}^\K(0)$
with $(v,e_j,t)\in U^{[1]}$.\\[2.5mm]
Let $y=(y_1,\ldots, y_n)$ and $z=(z_1,\ldots, z_n)$
be in $B_\delta(u)\cap U$,
and $q\in Q$. Then
\[
f_q(z)-f_q(y)=\sum_{j=1}^n \left(
f_q(z_1,\ldots, z_j, y_{j+1},\ldots, y_n)
-f_q(z_1,\ldots, z_{j-1},y_j,\ldots, y_n)\right)\,,
\]
where
\begin{eqnarray*}
\hspace*{-7mm}\lefteqn{f_q(z_1,\ldots, z_j, y_{j+1},\ldots, y_n)
-f_q(z_1,\ldots, z_{j-1},y_j,\ldots, y_n)}\quad \\
&= & D_jf(q,z_1,\ldots, z_{j-1},y_j,\ldots, y_n).(z_j-y_j)\\
& &  +\;\;
(z_j-y_j)R(q;z_1,\ldots, z_{j-1},y_j,\ldots, y_n;e_j;
z_j-y_j)\,.
\end{eqnarray*}
Hence
\begin{eqnarray*}
\hspace*{-7mm}\lefteqn{f_q(z)-f_q(y)-f'_q(u).(z-y)}\quad\\
&=&
\;\;\sum_{j=1}^n (z_j-y_j)\cdot
\left(D_jf(q,z_1,\ldots, z_{j-1},y_j,\ldots, y_n)
-D_jf(q,u)\right)\\
& &  \! + \sum_{j=1}^n (z_j-y_j)\cdot R(q;z_1,\ldots, z_{j-1},y_j,\ldots, y_n;e_j;
z_j-y_j)\,.
\end{eqnarray*}
Since $|z_j-y_j|\leq \|z-y\|$,
we see with
(\ref{whatn1}) and (\ref{whatn2})
that the $\eta$-gauge of each
of the $2n$ summands on the right hand side
of the previous equation is
bounded by $\ve\|z-y\|$.
Thus (\ref{wannashow})
follows, using~(\ref{itertriv}). 
\end{proof}
%
%
%
\begin{rem}\label{finremnonop}
If $E=\K^n$, then Lemma~\ref{C1paraunifV}
and its proof work just as well
if $U$ is not open,
but of the form $U=U_1\times\cdots\times U_n$,
where $U_1,\ldots, U_n\sub \K$
are subsets with dense interior.
\end{rem}
%
%
\begin{rem}\label{esspnt}
Since
Lemma~\ref{C1paraunifV}
can be used as a replacement for
Lemma~\ref{C1paraunif},
we see that
Theorem~\ref{genimp}
(the implicit function theorem)
and Theorem~\ref{advif}
(the inverse function theorem with parameters)
may be
extended to the
case of $C^1$-maps
over an arbitrary
valued field,
if~$F$ is finite-dimensional
and equipped with the canonical
vector topology.
Theorem~\ref{invparultra}
(the ultrametric inverse
function theorem with parameters)
extends
to $C^1$-maps
over an arbitrary
ultrametric field,
if~$F$ is finite-dimensional
and equipped with the canonical
vector topology.
Also, Theorem~\ref{basicinv}
(the inverse function theorem)
extends to $C^1$-maps
over valued fields,
if~$E$ is
finite-dimensional
and equipped with the canonical
vector topology.
\end{rem}
{\footnotesize
\noindent
{\bf Helge Gl\"{o}ckner}, TU~Darmstadt, FB~Mathematik~AG~5,
Schlossgartenstr.\,7,\\
64289~Darmstadt, Germany.
\,E-Mail: gloeckner@mathematik.tu-darmstadt.de}
\end{document}